\documentclass[USenglish]{article}
\usepackage{fullpage}
\usepackage{amsmath,amsfonts,amsthm,amssymb}
\usepackage{color,xcolor}
\usepackage{ifpdf}
\usepackage{psfrag}
\usepackage{graphicx,graphics}
\usepackage{hhtensor}
\usepackage{comment}
\usepackage[small]{caption}
\usepackage{subfigure}
\usepackage{tikz}
\usepackage{mathtools}
\usepackage{bbm}
\usepackage{algorithm,algorithmic}
\usepackage{pgfplots}
\usepackage{pgfplotstable}
\pgfplotsset{compat = 1.3}

\usetikzlibrary{fillbetween}
\usepgfplotslibrary{fillbetween}
\usepackage{xparse}
\usepackage{bigints}

\usepackage{url}
\usepackage{hyperref}
\usepackage{todonotes}
\usepackage{ulem} 
\newtheorem{theorem}{Theorem}[section]

\newtheorem{remark}[theorem]{Remark}
\newtheorem{definition}[theorem]{Definition}

\newtheorem{conclusion}[theorem]{Conclusion}
\usepackage{cancel}
\usepackage{algorithm,algorithmic}
\usepackage{enumitem}

\newcommand{\R}{\mathbb R}

\newcommand{\diff}[1]{{\mathrm{d}{#1}}}

\newcommand{\ddt}[1]   {\frac{\partial{#1}}{\partial{t}}}
\newcommand{\ddx}[1]   {\frac{\partial{#1}}{\partial{x}}}
\newcommand{\ddy}[1]   {\frac{\partial{#1}}{\partial{y}}}

\newcommand{\ddo}[2]   { \frac{\diff #1}{\diff #2}}

\newcommand{\bu}{\mathbf{u}}

\newcommand{\bF}{\mathbf{F}}
\newcommand{\bG}{\mathbf{G}}
\newcommand{\bbS}{\mathbf{S}}
\newcommand{\bFF}{\mathcal{F}}

\newcommand{\ww}[1]{\underline{#1}}

\newcommand{\bU}{\mathbf{U}}

\newcommand{\LL}{\mathcal{L}}

\newcommand{\xip}{x_{i+1/2}}
\newcommand{\xin}{x_{i-1/2}}

\newcommand{\yjp}{y_{j+1/2}}
\newcommand{\yjn}{y_{j-1/2}}

\newcommand{\iip}{i+1/2}
\newcommand{\iin}{i-1/2}

\newcommand{\jjp}{j+1/2}
\newcommand{\jjn}{j-1/2}

\newcommand{\bS}{\mathcal{S}}

\renewcommand{\vec}[1]{\ww{#1}}

\usepackage{bbm}
\def\L{\mathcal{L}}

\def\R{\mathbb{R}}
\def\bbc{\underline{\boldsymbol{y}}}
\def\bc{\boldsymbol{y}}

\newcommand{\mass}{\mathbb M}
\def\dt{\Delta t}

\definecolor{darkspringgreen}{rgb}{0., 0.55, 0.3}
\definecolor{dartmouthgreen}{rgb}{0.05, 0.5, 0.06}
\definecolor{etonblue}{rgb}{0.59, 0.78, 0.64}
\definecolor{airforceblue}{rgb}{0., 0.4, 0.66}
\definecolor{arylideyellow}{rgb}{0.91, 0.84, 0.42}
\definecolor{emerald}{rgb}{0.31, 0.78, 0.47}
\definecolor{uclagold}{rgb}{1.0, 0.7, 0.0}
\definecolor{cadmiumorange}{rgb}{0.93, 0.53, 0.18}

\newsavebox{\DelimiterBox}
\newlength{\DelimiterHeight}
\newlength{\DelimiterDepth}
\newsavebox{\ArgumentBox}
\newlength{\ArgumentHeight}
\newlength{\ArgumentDepth}
\newlength{\ResizedDelimiterHeight}
\newlength{\ResizedDelimiterDepth}

\ifluatex

\else

\fi

\newcommand{\Prod}{p}
\newcommand{\dest}{d}

\newcommand{\bby}{\mathbf{y}}

\hypersetup{
pdftitle={An Arbitrary High Order and Positivity Preserving Method for the Shallow Water Equations},
pdfauthor={M. Ciallella, L. Micalizzi, P. \"Offner, D. Torlo},
pdfpagemode=UseOutlines,
linkbordercolor=0 0 0,
linkcolor=red,
citecolor=blue,
colorlinks=true
}

\begin{document}
\title{An Arbitrary High Order and Positivity Preserving Method for the Shallow Water Equations}

\author{M. Ciallella$^{(1)}$, L. Micalizzi$^{(2)}$, P. \"Offner$^{(3)}$ and D. Torlo$^{(4)}$\footnote{Corresponding authors: \texttt{\ttfamily mirco.ciallella@inria.fr} (M. Ciallella), \texttt{\ttfamily lorenzo.micalizzi@math.uzh.ch} (L. Micalizzi), \texttt{\ttfamily poeffner@uni-mainz.de} (P. \"Offner), \texttt{\ttfamily davide.torlo@sissa.it} (D. Torlo) } \\ 
{\small (1): Team CARDAMOM, INRIA, Univ. Bordeaux, CNRS, Bordeaux INP, IMB, UMR 5251, France}\\
{\small (2): Institute of Mathematics,
University of Zurich, Switzerland,} \\
{\small (3): Institute of Mathematics, Johannes Gutenberg-University Mainz, Germany, }\\
{\small (4): SISSA mathLab, SISSA, via Bonomea 265, 34136, Trieste. }\\
{\small All authors contributed equally to this project.} }
\date{\today}
\maketitle

\begin{abstract}

In this paper, we develop and present an arbitrary high order well-balanced finite volume WENO method 
combined with the modified Patankar Deferred Correction (mPDeC) time integration method for the shallow water equations.
Due to the positivity-preserving property of mPDeC, the resulting scheme 
is unconditionally positivity preserving for the water height. To apply the mPDeC approach, we have to interpret the spatial semi-discretization in terms of production-destruction systems. 
Only small modifications inside the classical WENO implementation are necessary and we explain how it can be done. 
In numerical simulations, focusing on a fifth order method, we demonstrate the good performance of the new method and verify the theoretical properties. 

\end{abstract}
\vspace*{2mm}


\section{Introduction}
In the last years, the development of structure preserving  high-order methods for hyperbolic conservation/balance laws have been an active field of research \cite{abgrall2019reinterpretation,abgrall2020high, chenreview, kuzmin2020entropycg, gaburro2021unified, thomann2020all, veiga2021arbitrary}. In the context of the shallow water equations, one is mainly interested in maintaining positive levels of water height, in conserving the equilibrium/stationary states  and in entropy conservation/dissipation methods. There exists various ways to obtain these desired results, e.g.\ the applications of limiters for the positivity is only one example,
cf.\ \cite{berberich2021high, ricchiuto2011c, mantri2021well,noelle2007high,ranocha2017shallow,cheng2019new, xing2014survey} and references therein.

In this paper, we also deal with these issues and we present a new high-order, well-balanced,  positivity preserving method for the shallow water equation starting from a classical WENO scheme. In order to obtain well-balanced (WB) solutions, we subtract the residual of the a priori known stationary solution from our numerical scheme as shown in~\cite{berberich2021high}.\\
Then, to ensure the positivity of the water height, the modified Patankar Deferred Correction (mPDeC) method  is used for the time-integration of this variable. Even if (modified) Patankar (mP) methods have been already used inside a numerical method for fluid simulations~\cite{huang2019positivity, huang2018third, meister2016positivity}, those mP schemes have been based on extensions of classical RK methods and they are of maximum order three\footnote{Mostly, they have been only used for the source terms, multicomponent terms or in a post-processing process.}. By applying the modified Patankar trick inside the Deferred Correction (DeC)  framework, the authors of~\cite{offner2020arbitrary} were able to construct a conservative, arbitrarily high-order and positivity preserving method for production-destruction systems (PDS) of ordinary differential equations.\\
To obtain a positive WENO spatial reconstruction, a positive limiter must be used~\cite{zhang2010positivity,perthame1996positivity}.
In this work the mPDeC is applied for the first time to a PDE problem, with a finite volume WENO spatial discretization. This lead to a high order accurate method enjoying all the previously cited properties (WB, positivity preservation). \\
Finally, in order to apply mPDeC on the semi-discretized problem, the finite volume method must be rewritten into a PDS and we explain in details how this can be done. It should be stressed out that only small modifications inside the classical finite volume implementation are necessary as it can be seen in the reproducibility repository \cite{ourrepo}. 
To our opinion the modifications can be adapted to most WENO codes in a straightforward manner and the approach is a good alternative to already existing methods. \\ 
The paper is structured as follows:
In Section~\ref{se:SW_Equation}, we introduce the considered model, the classical shallow water equations, and we repeat the basic properties focusing on steady state solutions and the positivity of the water height. 
Next, in Section~\ref{sec_Space} we describe the used classical finite volume WENO approach from~\cite{shu1998essentially, shu1988efficient}
and its well-balanced modification from \cite{berberich2021high}.
We focus on the fifth order WENO method. However, all the ingredients to go to arbitrary high order can be found in the related repository \cite{ourrepo}. 
 In  Section~\ref{se_time_discretization}, 
the time-integration is considered focusing on Deferred Correction (DeC),  the modified Patankar approach
and its combination to the modified Patankar DeC (mPDeC) method developed in ~\cite{offner2020arbitrary}.
In Section~\ref{se_implemenation}, we describe how mPDeC can be combined with the WENO approach. It is important to interpret the 
semi-discretization in terms of production-destruction systems. Details of the implementation are given with additional algorithms. Then, in Section~\ref{se:numierics}, we verify the theoretical properties of the scheme with numerical simulations with WENO5 focusing on the high-order accuracy, the well-balanced and the positivity preserving properties. In addition, we demonstrate also the excellent performance for more challenging test cases. 
Finally, in Section~\ref{se:summary} we summarize the obtained results and perspectives for future works.

In Appendixes~\ref{sec:PrimitiveVarRecon} and~\ref{sec:WENO5GP4} we describe in details the WENO reconstruction and apply it for WENO5 with 4-points Gaussian quadrature rule, which, up to our knowledge, is not available in literature.

\section{Shallow Water Equations}\label{se:SW_Equation}
\subsection{Model}

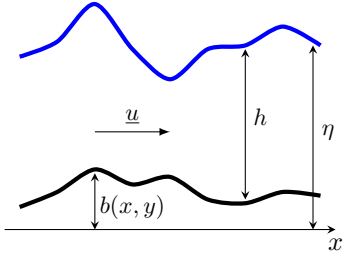
\begin{figure}
\centering
\begin{tikzpicture}

\draw [-stealth] (-0.2,-0.3) -- (4.2,-0.3);
\node [black,scale=1] at (4.2,-0.5) {$x$};

\draw [stealth-stealth] (1,-0.3) -- (1,0.45);
\node [black,scale=0.9] at (1.5,0) {$b(x,y)$};
\draw [stealth-stealth] (3,0.1) -- (3,2.1);
\node [black,scale=0.9] at (3.2,1.2) {$h$};
\draw [stealth-stealth] (3.9,-0.3) -- (3.9,2.15);
\node [black,scale=0.9] at (4.1,1.0) {$\eta$};

\draw [-latex] (1,1) -- (2,1);
\node [black,scale=0.9] at (1.5,1.2) {$\vec{u}$};

\draw [black,ultra thick] plot [smooth] coordinates {(0,0) (0.5,0.2) (1,0.5) (1.5,0.3) (2,0.4) (2.5,0.1) (3,0.05) (3.5,0.2) (4,0.15)};
\draw [blue, ultra thick]  plot [smooth] coordinates {(0,2) (0.5,2.2) (1,2.7) (1.5,2.1) (2,1.7) (2.5,2.1) (3,2.15) (3.5,2.4) (4,2.15)};

\end{tikzpicture}
\caption{Shallow Water Equations: definition of the variables.}\label{fig:SWE}
\end{figure}

The shallow water equations (SWE) model the behaviour of shallow free surface flows under the action of gravity.
They are used to simulate the  flows in rivers and coastal areas, and can be applied to predict tides, storm surge levels and coastline changes from hurricanes and ocean currents. They are also used in 
atmospheric flows, debris flows, and certain hydraulic structures like open channels and sedimentation tanks. 
SWEs take the form of non-homogeneous hyperbolic conservation laws with source terms modeling the effects of bathymetry and viscous friction. 
In this paper, we will consider the effect of the bathymetry as the only source term. If the bottom topography is assumed to be constant with respect to time, the SWEs can be recast in balance law form as:
\begin{equation}\label{eq:CL}
\ddt\bu + \nabla\cdot\bFF(\bu) = \bS(\bu,x,y) \;\;\;\;\;\text{on}\;\;\;\;\;\boldsymbol{\Omega}_T = \boldsymbol{\Omega}\times[0,T]\subset\R^2\times\R^+  
\end{equation}
with conserved variables, flux and source terms given by
\begin{equation}\label{eq:SWE}
\bu=\begin{bmatrix} h \\ hu  \\ hv \end{bmatrix}\;,\;\;
\bFF(\bu)=\begin{bmatrix}\bF &\bG \end{bmatrix}=\begin{bmatrix} hu & hv\\ hu^2+g\frac{h^2}{2}& huv \\ huv &hv^2+g\frac{h^2}{2} \end{bmatrix}\;,\;\;
\bS(\bu,x,y)=-gh\begin{bmatrix} 0 \\ \ddx b(x,y)  \\ \ddy b(x,y) \end{bmatrix}
\end{equation}
where $h$ represents the relative water height, $\vec{u}=(u,v)$ are the flow speed components, 
$g$ is the gravity acceleration and  $b(x,y)$ is the local bathymetry. 
The source term helps modeling the effects induced on the flow caused by the
bathymetry changes in space. Finally, it is also convenient to introduce the free surface water level $\eta:=h+b$. 
All the aforementioned variables can be better interpreted by looking at Figure~\ref{fig:SWE}. 

\subsection{Properties of the model}

As it is described \textit{inter alia} in \cite{bollermann2013well, berthon2016fully, lukavcova2007well, ricchiuto2015explicit, xing2014survey}, the construction and development of effective and accurate numerical methods for the shallow water equations have received much interest in the last decades and it is still ongoing. In particular, one is interested in schemes that preserve physical quantities or structures from the continuous level. 
In this paper, we are targeting two types of difficulties which are often encountered when working with SWE: the preservation of steady state solutions and water height positivity\footnote{We do not focus on entropy/energy preservation, see \cite{berthon2016fully, ranocha2017shallow}  and reference therein for this topic.}.

\subsubsection{Steady state solutions}

The SWE system~\eqref{eq:SWE} is known to admit some steady state solutions whose form depends on the equilibrium between the source terms $\bS$ and the remaining terms of the equations. The numerical simulations should be able to capture these behaviors even on coarse grids. 
Without additional techniques, many methods fail at balancing the source terms and the flux,  resulting in small pertubations of the steady state. The perturbations could be then amplified by the method causing a bad approximation of the exact behaviour. 
This situation is sometimes called \emph{numerical storm} in such context. 
To prevent it, one is interested in schemes that are capable of exactly balancing the flux and the source terms 
to obtain the desired steady-state solution.  
Numerical schemes enjoying this property are called \textbf{well-balanced schemes}. \\
The still water surface is often the first equilibrium taken in consideration.
It is given by 
\begin{equation}\label{lake_at_rest}
u=v=0; \qquad \eta(x,y,t) =h(x,y,t) +b(x,y)  \equiv\eta_0 \in \R^+_0, \quad \forall\, (x,y)\in \boldsymbol{\Omega},\, t \in [0,T].
\end{equation}
It represents a steady-state solution,  and is referred to as \textbf{lake at rest}. 
However, that is only a special case of the \textbf{moving water equilibrium}.
 Provided we verify the compatibility condition for the bathymetry which is $
(-v, u) \cdot  \nabla b=0, 
$
the steady state solutions are characterized by the invariants
\begin{equation}\label{moving_equilibrium} 
h(x,y,t)\bu(x,y,t)=Const. \quad  \text{ and } E(x,y,t)=\frac{1}{2}  ||{\bu(x,y,t)}||^2+g (h(x,y,t)+b(x,y,t))=Const. 
\end{equation}
Here, $E$ is the specific total energy (moving water equilibrium variable), cf. \cite{ricchiuto2015explicit} for more details. \\
Obviously, the lake at rest~\eqref{lake_at_rest} is a special case of \eqref{moving_equilibrium} when the velocity reduces to zero. 
Because of these reasons, it is desirable to solve~\eqref{eq:SWE}  with well-balanced schemes and also our described method has this property as we shall describe in Section~\ref{subsec_Well_balance}. 

\begin{remark}[C-Property]\label{C-Property}
As described in~\cite{cea2012unstructured, ricchiuto2011c, ricchiuto2015explicit}, instead of speaking of well-balanced schemes, one could 
alternatively say that a scheme enjoys the \textbf{C-property} if it preserves exactly the steady state~\eqref{lake_at_rest}.
However one still speaks of C-property when referring to other steady states~\eqref{moving_equilibrium}.
When the conservation of the steady state is not exact but is obtained within error rates below the formal accuracy of the scheme, one often speaks of \textbf{generalized C-property}~\cite{ricchiuto2011c}. 
\end{remark}

\subsubsection{Positivity of the solution}

The other major difficulty which pops up in many simulations of the SWE is the appearance of dry regions in many real applications as for example dam break problems, flood waves and run-up phenomena at  a coast with tsunamis. Here, the water height $(h=0)$ will be zero. As a result, all eigenvalues of the Jacobian of the flux coincide,  cf.~\cite{ricchiuto2015explicit} and the SWE model will be not strictly hyperbolic anymore. 
If, by numerical oscillations $h$ becomes negative, the problem is also not well-posed and the calculations will simply break down. It is thus essential for a good scheme to preserve the positivity of $h$ at any time and any point. 
Especially, in situations when we have dry and wet areas, the scheme has to be constructed such that it can handle these numerical challenges. In this paper, we shall target the issue of positivity of the solutions by combining a positivity preserving time-integration method together with the WENO approach. More about this will follow in Section~\ref{sub_modified}. 
\begin{remark}
Apart from the wet-dry areas in which the value of $h$ approaches $0$, we may have some regions of the space domain in which the bottom topography $b$ overcomes the water level resulting in areas which are completely dry. 
In this context, the given definition~\eqref{lake_at_rest} of the lake at rest steady state is unsatisfactory because in the mentioned regions we have $\eta_0< b$ which, if we strictly stick to~\eqref{lake_at_rest}, would give $h<0$ which does not correspond to the physics of the phenomen that must be modeled.
Therefore in such cases we need to modify the definition of the lake at rest steady state introducing the so-called \textbf{dry lake at rest}
\begin{equation}\label{dry_lake_at_rest}
u=v=0; \qquad h(x,y,t)=\begin{cases} \eta_0-b(x,y) , & \text{ if } b\leq  \eta_0 , \\
0, & \text{ else}, \end{cases} \quad \forall\, (x,y)\in \boldsymbol{\Omega},\, t \in [0,T].
\end{equation}
\end{remark}

\section{Space discretization: Finite Volume method}\label{sec_Space}
The system of PDEs considered herein will be solved by means of the Method Of Lines (MOL).
Thus, space and time are going to be treated separately. This section has the goal of presenting a standard high order
finite volume framework. In particular, we will consider a Cartesian setting but the method can be easily applied to a general unstructured framework in a straightforward way.

The computational domain $\boldsymbol{\Omega}$ is covered with $N_x\times N_y$ non-overlapping control volumes 
\begin{equation}
\Omega_{i,j} = [\xin,\xip]\times[\yjn,\yjp]
\end{equation}
with the dimensions given by $\Delta x = \xip-\xin$ and $\Delta y = \yjp-\yjn$.\\
Considering the system of hyperbolic balance laws described by~\eqref{eq:CL} and~\eqref{eq:SWE}, 
for the control volume $\Omega_{i,j}$ we can define the cell average at time $t$:
\begin{equation}\label{eq:cell average}
\bU_{i,j}(t):=\frac{1}{\Delta x \Delta y}\int_{\xin}^{\xip} \int_{\yjn}^{\yjp} \bu(x,y,t)\;\diff{x}\diff{y}.
\end{equation}
Integrating \eqref{eq:CL} over $\Omega_{i,j}$  provides the semi-discrete evolution formula with respect to $\bU_{i,j}$
\begin{equation}\label{eq:evol cell average}
\ddo{\bU_{i,j}(t)}{t} + \frac{1}{\Delta x}(\bF_{\iip,j}(t)-\bF_{\iin,j}(t)) + \frac{1}{\Delta y}(\bG_{i,\jjp}(t)-\bG_{i,\jjn}(t)) = \bbS_{i,j}(t) 
\end{equation}
where $\bbS_{i,j}$ is the cell average of the source terms over cell $\Omega_{i,j}$ at time $t$
\begin{equation}\label{eq:source average}
\bbS_{i,j}(t):=\frac{1}{\Delta x \Delta y}\int_{\xin}^{\xip} \int_{\yjn}^{\yjp} \bS(x,y,t)\;\diff{x}\diff{y}
\end{equation}
and $\bF_{\iip,j}$ and $\bG_{i,\jjp}$ are the cell-averages of the physical fluxes over cell boundaries at time $t$:
%
\begin{align}
\bF_{\iip,j}(t) &= \frac{1}{\Delta y}\int_{\yjn}^{\yjp}\bF(\bu(\xip,y,t))\;\diff{y},   \label{eq:flux F} \\
\bG_{i,\jjp}(t) &= \frac{1}{\Delta x}\int_{\xin}^{\xip}\bG(\bu(x,\yjp,t))\;\diff{x}  .  \label{eq:flux G}
\end{align}
%
Equation~\eqref{eq:evol cell average}, along with all the previous definitions, is so far exact. \\
From now until the end of this section, by an abuse of notation, we will not explicitly write the dependence on $t$ in all variables. $\bU_{i,j}$, $\bF_{\iip,j}$, $\bG_{i,\jjp}$ and $\bbS_{i,j}$ can be then approximated to the desired order of accuracy using appropriate quadrature 
formulae and reconstruction techniques.
From now on, we shall focus only on $\bF_{\iip,j}$ ($\bG_{i,\jjp}$ is obtained in a similar manner).
Many low order and high order reconstruction techniques have been developed in the last decades. 
Few of them can be found in~\cite{godunov1959finite,van1979towards,colella1984piecewise,shu1988efficient,gaburro2020high,bassi1997high,dumbser2010arbitrary,nunez2018hybrid}.

Once the reconstruction has been performed, at each face we have two sets of values of $\bU$, corresponding to $\xip^L$ and $\xip^R$,
which will be referred to as the left and right extrapolated values:
\begin{align}
\bu^L_{\iip,\theta} = \bu(\xip^L,y_\theta)\;,\;\;\;\bu^R_{\iip,\theta} = \bu(\xip^R,y_\theta).
\end{align}
In our case, the weighted essentially non-oscillatory (WENO) \cite{shu1998essentially, shu1988efficient} reconstruction has been considered to avoid severe oscillations at discontinuities.
By applying a consistent quadrature rule, and dropping the time dependence, the flux in the x-direction reads,
\begin{equation}\label{eq:flux with GP}
\bF_{\iip,j} = \frac{1}{\Delta y} \sum_{\theta=1}^{N_\theta} w_\theta \bF(\bu(\xip,y_\theta)) = \frac{1}{\Delta y} \sum_{\theta=1}^{N_\theta} w_\theta \hat\bF(\bu_{\iip,\theta}^L,\bu_{\iip,\theta}^R). 
\end{equation}
where the subscript $\theta=1,\ldots,N_\theta$ corresponds to different Gaussian points $y_\theta \in [y_{j-1/2},y_{j+1/2}]$ and weights $w_\theta \in [0,1]$.
The last step in the evaluation of the fluxes replaces $\bF({\bu}(\xip,y_\theta))$ in~\eqref{eq:flux with GP}
with a monotone and consistent numerical flux $\hat\bF({\bu}^L,{\bu}^R)$. The consistency is proved only
if $\hat\bF({\bu},{\bu}) = \bF({\bu})$. \\
In our scheme we employ a Rusanov-type Riemann solver:
\begin{equation}
\hat\bF({\bu}^L,{\bu}^R) = \frac{1}{2}\left(\bF({\bu}^R) + \bF({\bu}^L)\right) - \frac{1}{2}s_{max}\left({\bu}^R - {\bu}^L\right),
\end{equation}
where $s_{max}$ is the maximum eigenvalue of the normal flux-Jacobian of the system~\eqref{eq:CL}.\\
%
%
In order to express the dependence of some quantities on several cell averages $\bU_{i,j}$, it is useful to collect them all in the vector $\bU$.

\subsection{Weighted Essentially Non-Oscillatory (WENO) method}
As it is described in Shu's seminal paper \cite{shu1998essentially} the classical WENO (as well as ENO) approach contains three major steps: 
\begin{enumerate}
	\item Use the WENO reconstruction procedure which will be described in the following to obtain the values at the Gaussian points. This step involves two one-dimensional reconstructions in the two directions.
	\item Compute the numerical flux at quadrature points and integrate them at each cell interface.
	\item Form a classical semidiscrete FV method and apply a time-integration method to update the cell averaged values. 
\end{enumerate}

\subsubsection{Scalar reconstruction}
The goal of the WENO method is to compute point-wise value of variable of interest $u(x,y)$ at Gaussian quadrature points $(\xip,y_\theta)$,
%
%
in order to have a conservative and high order accurate procedure.
In general, two ways can be followed to obtain the same result: genuine multidimensional reconstruction~\cite{shi2002technique} and 
dimension-by-dimension reconstruction~\cite{titarev2004finite,casper1993finite}. 
The latter is a procedure made up by successive one-dimensional reconstruction sweeps and it is much simpler and less computationally expensive than the genuine multidimensional one.
For this reason, we shall only focus on this one. 

The high order reconstructed variables we are looking for will be referred to as $u^L_{\iip,\theta}$ and $u^R_{\iip,\theta}$.
For the left values, we need to reconstruct the variable inside the cell $\Omega_{i,j}$, while, for the right values, similar arguments apply on the cell $\Omega_{i+1,j}$.  We aim at reconstructing the variables with an accuracy of order $p$ ($p$ odd). So, we define a stencil of $p$ cells, $\lbrace\Omega_{l_x,l_y},\quad l_x = i-r+1, \dots, i+r-1, \; l_y = j-r+1,\dots, j+r-1 \rbrace,$ 
%
%
where $2r-1=p$. For instance, WENO5 has accuracy $p=5$, with $r=3$, and uses a $5$-cells stencil from $i-2$ to $i+2$.\\
In the first step of the two-dimensional reconstruction, a one-dimensional WENO reconstruction along the x-direction
is performed obtaining the averages at cell interface $x_{i+1/2}$ with respect to the y-direction for $l_y=j-r+1,\dots, j+r-1$
\begin{align}
v^R_{l_y} = \frac{1}{\Delta y}\int_{y_{l_y-1/2}}^{y_{l_y+1/2}}u(\xip^R,y)\;\diff{y},\qquad 
v^L_{l_y} = \frac{1}{\Delta y}\int_{y_{l_y-1/2}}^{y_{l_y+1/2}}u(\xip^L,y)\;\diff{y} .
\end{align}
In the second sweep we perform another one-dimensional reconstruction along the y-direction 
in the Gaussian integration points on the y-axis $(x=\xip,y=y_\theta)$, with $y_\theta \in [y_{j-1/2},y_{j+1/2}]$.
The reconstructed values can be, more generally, defined for each WENO sweep as the one-dimensional averages $q_i=\frac{1}{\Delta\xi}\int_{\xi_{i-1/2}}^{\xi_{i+1/2}} q(\xi)\;\diff{\xi}$
of a function $q(\xi)$
%
%
where $\Delta\xi=\xi_{i+1/2}-\xi_{i-1/2}$ is the cell size.\\
For each one-dimensional step of the procedure, there are $r$ candidate stencils for reconstruction. 
For each of these stencils, made up by $r$ cells, there is a correspoding polynomial of degree $(r-1)$ referred as $p_m(\xi)$ $m = 0,\ldots,r-1$.
%
%
The goal of the WENO reconstruction is that of using all information coming from the $r$ stencils employed for the reconstruction. 
For this reason, the WENO approach defines the reconstructed value as a convex combination of the $r$ values of all polynomials in each quadrature point,
weighted with positive nonlinear weights. The weights are chosen in order to achieve $(2r-1)$th order of accuracy when the 
solution is smooth and prefer the smoother stencils when discontinuities occur in the field. 
For a given (quadrature) point $\tilde\xi$ the design of weights consists of three steps.
Firstly, the optimal linear weights $d_m$ are sought so that the combination of all polynomials with these weights produces the 
polynomial of degree $(2r-2)$ corresponding to the large stencil.
Then, the nonlinear weights $\omega_m$ can be defined as $\omega_m = \frac{\alpha_m}{\sum^{r-1}_{k=0}\alpha_k}$ with $\alpha_m = \frac{d_m}{(\beta_m+\epsilon)^2}$,
%
%
where $\epsilon$ is a small constant introduced to avoid division by zero (we use $\epsilon=10^{-6}$ in the simulations) and $\beta_m$ are the smoothness indicators
%
\begin{equation}
\beta_m = \sum_{k=1}^{r-1} \int_{\xi_{i-1/2}}^{\xi_{i+1/2}} \left(\frac{\diff{}^k}{\diff{x}^k} p_m(\xi)\right)^2 \Delta\xi^{2k-1}\diff{\xi}\;,\;\;\;m = 0,\ldots,r-1.
\end{equation}
%
%
If some of $d_m$ are negative then a special procedure must be used to tackle the reconstruction problem~\cite{shi2002technique}.\\ 
The final WENO reconstructed quantity is given by $q(\tilde\xi) = \sum_{k=0}^{r-1} p_k(\tilde\xi)\omega_k$.\\ 
%
%
The numerical experiments presented herein have been performed through a piece-wise parabolic WENO5 reconstruction ($r=3$), which
formally corresponds to fifth order accurate approximation for smooth solutions. However, in order to actually retain the fifth-order accuracy the quadrature formulae must be consistent with the WENO reconstruction.  
As Titarev and Toro stated in~\cite{titarev2004finite}, the best results in terms of accuracy and computational cost for $r=3$ are
obtained if the two-point Gaussian quadrature rule is used. 
However, 
this leads eventually to a formal fourth order of accuracy. 
For this reason, we implemented the four-point Gaussian quadrature rule with positive optimal weights. 
%
%
%
Up to our knowledge, the two-point Gaussian quadrature 
has already been thoroughly discussed in 
many references cited above. However, the case with 4-point Gaussian quadrature 
has not been fully described in literature, hence, we are going to introduce all the coefficients needed to use such formula in Appendix~\ref{sec:WENO5GP4} and a Matlab script to compute the weights and coefficients for all orders is provided in~\cite{ourrepo}.

\begin{remark}[Positivity limiter]\label{rem:positivity_limiter}
We aim at a positive solution and during the reconstruction procedure, it might happen that $h(x^L_{i+1/2})$ or $h(x^R_{i+1/2})$ become negative. 
In order to ensure that positive cell averages lead to positive reconstructions at the cell interfaces, we use the positivity limiter introduced by Perthame and Shu~\cite{perthame1996positivity} and developed for two dimensional problems in~\cite{zhang2010positivity}. 
The limiter is used in the simulation section with a parameter $\varepsilon=10^{-6}$ as minimum water height, if not otherwise specified. 
We refer to~\cite{xing2014survey,zhang2010positivity} for details on the implementation.
This limiter when used in combination with the forward Euler (FE) method restricts the CFL conditions to $\text{CFL}^{\text{FE}}:=w^{\text{Lobatto}}_{1}$ the weight of the Gauss--Lobatto quadrature rule of the corresponding space accuracy. For instance, with WENO5, $\text{CFL}^{\text{FE}}=1/12$. The restriction slightly improves for high order SSPRK methods, for example we have $\text{CFL}^{\text{SSPRK(5,4)}}\approx 1.508 \cdot\text{CFL}^{\text{FE}}$. Unfortunately, explicit SSPRK methods are at most fourth order accurate, so for fifth order schemes (as DeC5), there is no warranty that the solution stays nonnegative under any CFL condition. Instead, we highlight that the new presented approach is unconditionally positive and thus not subjected to any CFL restriction.
\end{remark}

\subsection{Well-Balanced modification of the standard Finite Volume method}
\label{subsec_Well_balance}
%
In order to achieve Well-Balancing with respect to the (eventually dry) lake at rest steady state, in this work we coupled the WENO formulation with a simple modification firstly introduced in \cite{berberich2021high}. The modification consists in recasting the original problem into an equivalent one in terms of the deviation of the seeked solution $\bU$ from the reference solution $\tilde\bU$ which must be preserved. In the particular case in which a steady solution ($\ddt{\tilde{\bU}}=0$) must be preserved, the modification leads to the new problem
\begin{align}
\begin{aligned}
\frac{\diff{}}{\diff{t}} \bU_{i,j} + 
&\frac{1}{\Delta x}(\bF_{\iip,j}(\bU)-\bF_{\iin,j}(\bU)) - \frac{1}{\Delta x}(\bF_{\iip,j}(\tilde\bU)-\bF_{\iin,j}(\tilde\bU)) +\\
&\frac{1}{\Delta y}(\bG_{i,\jjp}(\bU)-\bG_{i,\jjn}(\bU)) - \frac{1}{\Delta y}(\bG_{i,\jjp}(\tilde\bU)-\bG_{i,\jjn}(\tilde\bU)) =\\
&\bbS_{i,j}(\bU) - \bbS_{i,j}(\tilde\bU), \label{eq:WB2} 
\end{aligned}
\end{align}
which can be interpreted as a classical finite volume formulation with modified fluxes and source:
\begin{equation}\label{eq:numerical_flux1}
\begin{aligned}
\overline{\bF}_{\iip,j} (\bU)&=\bF_{\iip,j}(\bU)-\bF_{\iip,j}(\tilde\bU) , \\
\overline{\bG}_{i,\jjp} (\bU)&=\bG_{i,\jjp}(\bU)-\bG_{i,\jjp}(\tilde\bU),  \\
\overline{\bbS}_{i,j}(\bU) &=\bbS_{i,j}(\bU) - \bbS_{i,j}(\tilde\bU).
\end{aligned}
\end{equation}
This approach is very easy to code and, further, the structures related to the steady reference solution can be computed in advance once and then used for every timestep without affecting the computational time.
It must be underlined that, with this technique, all cell average computations, WENO reconstruction and source terms of the reference solution are performed following the same procedures and quadrature rules carried out for solving the balance law. Hence, all the terms always match when at the equilibrium.

\section{Time discretization}\label{se_time_discretization}

After the description of the space discretization, we will introduce the time discretization in the this section.
Due to the MOL approach, it is enough to focus here on the simple ODE case. 
The time discretization is one of the major point of this paper since we want to apply for the first time the arbitrary high-order, conservative and positivity preserving  modified Patankar 
Deferred Correction method (mPDeC) together with the described WENO approach to the shallow water equations resulting in a 
high-order, conservative, unconditionally positivity preserving,  non-oscillatory and well-balanced scheme. 
In agreement with \cite{meister2014unconditionally}, we refer to the numerical methods which are provably positive with respect to the water height without any CFL restriction, i.e., for all time steps $\dt>0$, as \textbf{unconditionally} positivity preserving, in the context of the numerical solution of the shallow water equations.

Before describing the combined algorithm, we introduce the Deferred Correction (DeC) method~\cite{dutt2000dec} as described in~\cite{abgrall2017high, abgrall2020high, torlo2020hyperbolic}
and we repeat its modification using the Patankar trick from~\cite{offner2020arbitrary, offner2020stability, torlo2021stability}.

\subsection{Deferred Correction method}

The general DeC approach has been introduced in~\cite{dutt2000dec} and has been further 
developed and applied in~\cite{christlieb2010integral, liu2008strong, minion2003dec} in different contexts,
whereas a simplified version  was presented in~\cite{abgrall2017high}. In~\cite{abgrall2017high}, a compact operator 
notation was also introduced and we shall follow this framework herein.
However, even if the notation changed, the main idea of DeC is always the same. DeC is based on the 
 Picard-Lindel\"of Theorem in the continuous setting and the classical proof makes use of Picard iterations 
 to minimize the error and to obtain convergence. 
 DeC is constructed to mimic these Picard iterations at the discrete level and decreases the approximation error
 in several iterative steps.
 To explain the method, we consider the following time-dependent initial 
value problem 
 \begin{equation}\label{eq:initial_prob}
\begin{aligned}
y'(t) = f(y(t)),  \qquad y(t_0)= y_0,
\end{aligned}
\end{equation}
where $y:\R\to \R^S$ and $f: \R^S\to\R^S$ resulting from our MOL approach.
 
  For our description, two operators are introduced: $\LL^1$ and $\LL^2$.
Here,  the $\LL^1$ operator represents a low-order easy-to-solve numerical scheme,
e.g.\ the explicit Euler method, 
and $\LL^2$ is a high order operator that can present difficulties in
its practical solution, e.g.\ an implicit RK scheme.
The DeC method can be written as a combination of these two operators.
\begin{figure}[ht]
	\centering
	\begin{tikzpicture}
		\draw [thick]   (0,0) -- (10,0) node [right=2mm]{};
		\fill[black]    (0,0) circle (1mm) node[below=2mm] {$t^n=t^{n,0}=t^0 \,\, \quad$} node[above=2mm] {$y^0$}
		(2,0) circle (0.7mm) node[below=2mm] {$t^{n,1}=t^1$} node[above=2mm] {$y^1$}
		(4,0) circle (0.7mm) node[below=2mm] {}
		(6,0) circle (0.7mm) node[below=2mm] {$t^{n,m}=t^m$} node[above=2mm] {$y^m$}
		(8,0) circle (0.7mm) node[below=2mm] {}
		(10,0) circle (1mm) node[below=2mm] {$\qquad t^{n,M}=t^M=t^{n+1}$} node[above=2mm] {$y^M$}; 
	\end{tikzpicture} \caption{Time interval divided into subintervals}\label{Fig:Time_interval}
\end{figure}
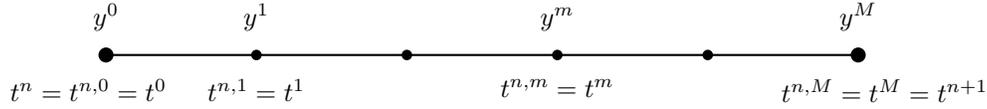
Given a time interval $[t^n, t^{n+1}]$, we subdivide 
it into $M$ subintervals  $\lbrace [t^{n,m-1},t^{n,m}]\rbrace_{m=1}^M$,
where $t^{n,0} = t^n$ and $t^{n,M} = t^{n+1}$. Therefore, we  mimic for every 
subinterval $[t^0, t^m]$ the Picard--Lindel\"of Theorem for both operators $\LL^1$
and $\LL^2$. With an abuse of notation, we drop the dependency on the timestep $n$ for subtimesteps $t^{n,m}$ 
and substates $y^{n,m}$ as denoted in Figure~\ref{Fig:Time_interval}.  

Then, the $\LL^2$ operator is given by
\begin{equation}\label{eq:L2operator}
\LL^2(y^0, \dots, y^M) :=
\begin{cases}
y^M-y^0 - \Delta t \sum_{r=0}^M \theta_r^M f(y^r)\\
\vdots\\
y^1-y^0 - \Delta t \sum_{r=0}^M \theta_r^1 f(y^r)
\end{cases}.
\end{equation}
where $\theta_r^m$ are the weights of a high order quadrature rule in $\lbrace t^m \rbrace_{m=0}^M$.

The $\LL^2$ operator represents a high order numerical scheme if set equal
to zero, i.\,e., $\LL^2(y^0, \dots, y^M)=0$. The order depends on the distribution of the subtimesteps, for instance, with $M$ equispaced subtimesteps one obtains $(M+1)$th order, while with $M$ Gauss--Lobatto quadrature subtimesteps one has $(2M)$th order. Unfortunately, the resulting scheme is implicit and, further, the terms $f$ may be nonlinear. 
The $\LL^1$ operator is given by the forward Euler discretization for each state $y^m$ in the time interval, i.\,e.,  
\begin{equation}\label{eq:L1}
\LL^1(y^0, \dots, y^M) :=
\begin{cases}
 y^M-y^0 - \beta^M \Delta t f(y^0) \\
\vdots\\
y^1- y^0 - \beta^1 \Delta t f(y^0)
\end{cases}
\end{equation}
with coefficients $\beta^m:=\frac{t^m-t^0}{t^M-t^0}$.\\
To simplify the notation and to describe  DeC, we  introduce the matrix of states for the variable $y$ at all subtimesteps.
\begin{align}\label{eq:definition_bbc}
&\bby :=  (y^0, \dots, y ^M) \in \R^{(M+1)\times S}, \text{ such that }\\
&\LL^1(\bby) := \LL^1(y^0, \dots, y^M) \text{ and } \LL^2(\bby) := \LL^2(y^0, \dots, y^M) .
\end{align}
The DeC algorithm uses a combination of the $\LL^1$ and $\LL^2$ operators
to provide an iterative procedure. The aim is to recursively approximate $\bby^*$, the numerical solution of
the $\LL^2(\bby^*)=0$ scheme, similarly to the Picard iterations in the 
continuous setting. The successive states of the iteration process will be denoted 
by the superscript $(k)$, where $k$ is the iteration index, e.g. $\bby^{(k)}\in \R^{(M+1)\times S}$.
The total number of iterations (also called correction steps in the following) is denoted by $K$.
To describe the procedure, we have to refer to both 
 the $m$-th subtimestep and the $k$-th iteration of the DeC algorithm. We will indicate the variable by $y^{m,(k)} \in \R^S$.
 Finally, the DeC method can be written as \\

 \centerline{\textbf{DeC Algorithm}}

\begin{equation}\label{DeC_method}
\begin{split}
&y^{0,(k)}:=y(t^n), \quad k=0,\dots, K,\\
&y^{m,(0)}:=y(t^n),\quad m=1,\dots, M,\\
&\LL^1(\bby^{(k)})=\LL^1(\bby^{(k-1)})-\LL^2(\bby^{(k-1)}) \text{ with }k=1,\dots,K,
\end{split}
\end{equation}
where $K$ is the number of iterations that we  want to compute. 

Using the procedure~\eqref{DeC_method}, we  need, in particular, 
as many iterations as the desired  order of accuracy $p$, i.\,e., $K=p$. This means that we choose the number of subtimesteps in a way that the order of the $\L^2$ operator is itself equal to $p$.
In practice, for each correction and each subtimestep, $\L^{1,m}(\bbc^{(k)})=\L^{1,m}(\bbc^{(k-1)})-\L^{2,m}(\bbc^{(k)})$ reduces to solve 
\begin{equation}
	\label{eq:explicit_dec_explicit}
	y_\alpha^{m,(k)}-y^0_\alpha -\Delta t\sum_{r=0}^M \theta_r^m  f_{\alpha}(y^{r,(k-1)})=0, \qquad \forall \alpha = 1, \dots, I.
\end{equation}

For more information and properties of the DeC approach, we refer to~\cite{abgrall2021relaxation, veiga2021dec} and references therein. 
In the following, we explain how to adapt the DeC approach to obtain a conservative and positivity preserving time integration scheme.

\subsection{Patankar method for production-destruction systems}

Many problems  \eqref{eq:initial_prob} in nature can be written as a 
production destruction system (PDS) for the unknown $y\in \R^S$
\begin{equation}
\label{eq:pd-system}
  f_\alpha(y)
  =
  \ \sum_{\beta=1}^S (\Prod_{\alpha,\beta}(y) - \dest_{\alpha,\beta}(y)),
\end{equation}
where $\Prod_{\alpha,\beta}, \dest_{\alpha,\beta} \geq 0$ are the production and destruction
terms, respectively. The production and destruction terms are conveniently written as
matrices.
Applications for PDS are for example the biological and/or chemical reactions such as algal bloom \cite{burchard2003high}.
Also parts (or all) of  the semi discretization of hyperbolic conservation/balance laws can be interpreted in such PDS system as described in~\cite{huang2019positivity, huang2018third, meister2014unconditionally}
and also later in this work.  The calculated solutions are often describing physical quantities  that enjoy some properties, for instance concentrations of chemicals or water height in the context of SWE should be nonnegative.  
The following definition may be introduced for ODE systems:  
\begin{definition}\label{def:POSCONSPDS} 
  An ODE~\eqref{eq:initial_prob} is called \emph{positive}, if positive initial data $y_0 > 0$
  result in positive solutions $y(t) > 0, \forall t$. Here, inequalities for
  vectors are interpreted componentwise, i.e., $y(t) > 0$ means
  $\forall \alpha\colon y_\alpha(t) > 0$.
  A PDS~\eqref{eq:pd-system} is \emph{conservative},
  if $\Prod_{\alpha,\beta}(y) = \dest_{\beta,\alpha}(y), \, \forall \, \alpha,\,\beta,\,y$.
\end{definition}
These properties should be preserved by the numerical scheme as well. Thus, we introduce
the following discrete counterpart.
\begin{definition}
  A numerical method computing $y^{n+1} \approx y(t_{n+1})$ given
  $y^{n} \approx y(t_n)$ is called \emph{conservative}, if
  $\sum_\alpha y^{n+1}_\alpha = \sum_\alpha y^{n}_\alpha$.
  It is called \emph{unconditionally positive}, if $y^{n} > 0$ implies $y^{n+1} > 0$ for any time step $\dt>0$.
\end{definition}
From literature~\cite{burchard2003high}, it is well-known that the implicit Euler method is conservative and \textit{unconditionally }positive preserving whereas the explicit Euler method is only
conservative (it might be positive under time step restrictions). To avoid solving a fully nonlinear system of equations, the so--called Patankar modifications have been applied to the explicit Euler method. 
To build an unconditionally positive numerical scheme,
Patankar had the idea~\cite{patankar1980numerical} of firstly 
weighting the destruction term in the original explicit Euler method
with a coefficient as follows
\begin{equation}\label{eq:patankar}
y_\alpha^{n+1}= y_\alpha^n+\Delta t  \left( \sum_{\beta=1}^S \Prod_{\alpha,\beta}(y^n) - 
\sum_{\beta=1}^S \dest_{\alpha,\beta}(y^n) \frac{y^{n+1}_\alpha}{y_\alpha^n} \right), \quad \alpha=1,\dots, S.
\end{equation}
Indeed, the resulting scheme~\eqref{eq:patankar} is unconditionally positive and the implicit terms can be collected on the left hand side, but the conservation relation 
is violated. 
Burchard et al. had the idea~\cite{burchard2003high} not only to weight the destruction term but also the production term:
\begin{equation}\label{eq:mod_patankar}
y_\alpha^{n+1}=y_\alpha^n+\Delta t  \left( \sum_{\beta=1}^S \Prod_{\alpha,\beta}(y^n) \frac{y^{n+1}_\beta}{y_\beta^n} - \sum_{\beta=1}^S \dest_{\alpha,\beta}(y^n) \frac{y^{n+1}_\alpha}{y_\alpha^n} \right), \quad \alpha=1,\dots, S.
\end{equation}
They called their constructed scheme~\eqref{eq:mod_patankar} \textbf{modified Patankar scheme} and proved that it is 
 is unconditionally positive and conservative. The resulting scheme is linearly implicit, meaning that collecting all the implicit terms on the left hand side, we obtain a linear system at each time iteration.
 Based on this technique, extensions to second and third order modified Patankar Runge--Kutta (MPRK) methods have been made by several researchers in such context, cf.~\cite{huang2019positivity, huang2018third, kopecz2018unconditionally, kopecz2019existence}. 
Also the semi implicit RK methods proposed in~\cite{chertock2015steady} can be interpreted as Patankar methods as they weight only the destruction terms~\cite{torlo2021stability}.
Finally, 
in~\cite{offner2020arbitrary} an arbitrarily high-order, conservative and positivity preserving scheme based on the DeC framework has been constructed.
Herein we describe the main idea. For the details of the properties and proofs, we refer again to~\cite{offner2020arbitrary}.

\subsection{Modified Patankar Deferred Correction method}
\label{sub_modified}
The modified Patankar Deferred Correction (mPDeC) is based on the DeC algorithm \eqref{eq:explicit_dec_explicit} and it consists in a modification of the $\LL^2 $ operator through the modified Patankar trick.
This amounts to weight the production-destruction terms with respect to the intermediate approximations. 

Using the fact that initial states $y_\alpha^{0,(k)}$
are identical for any correction $k$, the mPDeC correction steps
can be rewritten \cite{offner2020arbitrary} for $k=1,\dots,K$, $m =1,\dots, M$ and $\forall \alpha =1,\dots, S$ as
\begin{equation}
\label{eq:explicit_dec_correction}
y_\alpha^{m,(k)}-y^0_\alpha -\sum_{r=0}^M \theta_r^m \Delta t   \sum_{\beta=1}^S \left( \Prod_{\alpha,\beta}(y^{r,(k-1)}) \frac{y^{m,(k)}_{\gamma(\beta,\alpha, \theta_r^m)}}{y_{\gamma(\beta,\alpha, \theta_r^m)}^{m,(k-1)}} - \dest_{\alpha,\beta}(y^{r,(k-1)})  \frac{y^{m,(k)}_{\gamma(\alpha,\beta, \theta_r^m)}}{y_{\gamma(\alpha,\beta, \theta_r^m)}^{m,(k-1)}} \right)=0,
\end{equation}
where $\theta_r^m$ are the DeC quadrature weights in time and 
$$\gamma(\alpha,\beta,\theta):=\begin{cases}
	\alpha & \text{if }\theta \geq 0,\\
	\beta & \text{if }\theta <0.
\end{cases}$$
Finally, the new numerical solution
is $y^{n+1}=y^{M,(K)}$.
As in the classical DeC framework, the choice of the distribution, the number of subtimesteps $M$ and the number of iterations $K$ determines the order of accuracy of the scheme. 
To reach order $p$, classically $M=p-1$ equispaced subintervals and $K=p$ corrections should be used. As proven in~\cite{offner2020arbitrary}, the scheme is conservative, positivity preserving and can reach arbitrary high order.   \\
A description of the assembly of the mass matrix of the system \eqref{eq:explicit_dec_correction} is described in Algorithm \ref{algo:mass} and one timestep of the mPDeC algorithm is sketched in Algorithm~\eqref{DeC_method} where the evolution formula is given by Algorithm \ref{algo:mPDeC}.
\begin{remark}[Subtimestep distribution]
In our numerical simulations, we apply Gauss-Lobatto nodes in every timestep. They have the advantage of requiring less subtimesteps to reach $p$th order of accuracy. In the following we will use $M=3$ Gauss-Lobatto subtimesteps, which guarantee 6th order of accuracy for the operator $\L^2$ and $K=5$ iterations aiming at a 5th order scheme to match the spatial discretization accuracy of WENO5. 
\end{remark}

\begin{remark}[Solution of the linear system]
	At each subtimestep $m$ and iteration $(k)$ we need to solve the linear system given by \eqref{eq:explicit_dec_correction}. The mass matrix obtained has the following form:
	\begin{equation}\label{eq:matrixMPDeC}
		\mathbb M (y^{m,(k-1)})_{\alpha,\beta} = 
		\begin{cases} 
			1+\Delta t \displaystyle\sum_{r=0}^M \displaystyle\sum_{\beta=1}^S \frac{\theta_r^m}{y_\alpha^{m,(k-1)}} \left( d_{\alpha,\beta}(y^{r,(k-1)})\mathbbm{1}_{\{\theta_r^m>0\}} - p_{\alpha,\beta}(y^{r,(k-1)})\mathbbm{1}_{\{\theta_r^m<0\}} \right), &\text{for } \alpha=\beta, \\ 
			-\Delta t \displaystyle\sum_{r=0}^M \frac{\theta_r^m}{y_\beta^{m,(k-1)}} \left( p_{\alpha,\beta}(y^{r,(k-1)})\mathbbm{1}_{\{\theta_r^m>0\}} - d_{\alpha,\beta}(y^{r,(k-1)})\mathbbm{1}_{\{\theta_r^m<0\}} \right), &\text{for } \alpha \neq \beta,
		\end{cases}
	\end{equation}
	where $\mathbbm{1}$ is the indicator function. The mass matrix assembly algorithm is described in Algorithm \ref{algo:mass}. The linear system will then read
	\begin{equation}
		\mathbb M(y^{m,(k-1)}) y^{m,(k)} = y(t^n).
	\end{equation}
\end{remark}

\begin{remark}[Division on almost wet areas]\label{rem:division_zero}
	When the water height is low, we might encounter troubles in computing the divisions in \eqref{eq:matrixMPDeC} as the denominator might be very small. They hypothesis behind the production and destruction system that says that as $h_\alpha\to 0$ also $d_\alpha\to 0$, can be difficult to be obtain at a numerical level. Hence, to be sure that those divisions do not lead to extremely high values when they should go to 0, we slightly modify the way we implement the division as suggested in \cite{meister2016positivity}.
	Given any numerator $n$ and denominator $d$ of \eqref{eq:matrixMPDeC}, we approximate the division by
	\begin{equation}
		\frac{n}{d}\approx \begin{cases}
			0 & d < \varepsilon,\\
			\frac{2d\cdot n}{d^2+\max\lbrace d^2,\varepsilon\rbrace} & d \geq  \varepsilon,
		\end{cases}
	\end{equation}
with $\varepsilon$ a small tolerance value. Along the computations, if not specified, we will use $\varepsilon:= 10^{-6}$.
This formulation allow to smoothly pass from $\frac{n}{d}$ to $0$ as $d\to 0$. Moreover, when $d^2\geq \varepsilon$ the division will be exact. 
\end{remark}

In the following, we apply the mPDeC time marching algorithm \eqref{eq:explicit_dec_correction} to the WENO finite volume semidiscretization to solve 
the shallow water equations. Below we describe the actual implementation procedure.

\begin{algorithm}
	\fontsize{10pt}{10pt}\selectfont
	\caption{Mass} 
	\begin{algorithmic}[1]
		{\REQUIRE Production-destruction functions $p_{\alpha,\beta}(\cdot),\,d_{\alpha,\beta}(\cdot)$, $\Delta t$, previous correction variables $\bbc^{(k-1)}$, current subtimestep $m$.
			\STATE $\mass\gets \mathbb I$
			\FOR{$\alpha=1 $ \TO $ S$}
			\FOR{$\beta=1$ \TO $S$}
			\FOR{$r=0$ \TO $M$}
			\IF{$\theta_r^m\geq 0$}
			\STATE $\mass_{\alpha,\beta} \gets  \mass_{\alpha,\beta} -\Delta t \theta_r^m \frac{p_{\alpha,\beta}(\bc^{r,(k-1)})}{y_\beta^{m,(k-1)}}$
			\STATE $\mass_{\alpha,\alpha} \gets  \mass_{\alpha,\alpha} +\Delta t \theta_r^m \frac{d_{\alpha,\beta}(\bc^{r,(k-1)})}{y_\alpha^{m,(k-1)}}$
			\ELSE
			\STATE $\mass_{\alpha,\beta} \gets  \mass_{\alpha,\beta} +\Delta t \theta_r^m \frac{d_{\alpha,\beta}(\bc^{r,(k-1)})}{y_\beta^{m,(k-1)}}$
			\STATE $\mass_{\alpha,\alpha} \gets  \mass_{\alpha,\alpha} -\Delta t \theta_r^m \frac{p_{\alpha,\beta}(\bc^{r,(k-1)})}{y_\alpha^{m,(k-1)}}$	
			\ENDIF
			\ENDFOR
			\ENDFOR
			\ENDFOR
			\RETURN $\mass$
		}
	\end{algorithmic}\label{algo:mass}
\end{algorithm}

\begin{algorithm}
	\fontsize{10pt}{10pt}\selectfont
	\caption{mPDeC Update formula} 
	\begin{algorithmic}[1]
		{\REQUIRE $\bbc^{(k-1)}$, $\Delta t$, production-destruction functions $p_{\alpha,\beta}(\cdot),\,d_{\alpha,\beta}(\cdot)$, m.
			\STATE Compute the mass matrix $\mass(\bc^{m,(k-1)})\gets $Mass($p_{\alpha,\beta}(\cdot),\,d_{\alpha,\beta}(\cdot)$,$\Delta t,\bbc^{(k-1)}, m)$ using Algorithm \ref{algo:mass}
			\STATE Compute $\bc^{m,(k)}$ solving the linear system $\mass(\bc^{m,(k-1)}) \bc^{m,(k)} = \bc^{n}$  given by \eqref{eq:explicit_dec_correction} with Jacobi Algorithm \ref{algo:jacobi}
			\RETURN $\bc^{m,(k)}$
		}
	\end{algorithmic}\label{algo:mPDeC}
\end{algorithm}

\begin{algorithm}
	\fontsize{10pt}{10pt}\selectfont
	\caption{Jacobi iterative method} 
	\begin{algorithmic}[1]
		{\REQUIRE $\mathbb D$ diagonal of the matrix, $\mathbb L$ off-diagonal terms of the matrix, $\mathtt r$ right hand side of the system, $\mathtt{ tol}$ tolerance.
			\STATE $\mathtt{err}\gets 2\cdot \mathtt{tol}$,\quad  $k\gets0$,\quad $\mathtt x^k\gets \mathtt r$
			\WHILE{$\mathtt{err} > \mathtt{tol}$} 
			\STATE $k\gets k+1$
			\STATE $\mathtt x^{k+1}\gets \mathbb D^{-1} (\mathtt r - \mathbb L \mathtt x^{k})$
			\STATE $\mathtt{err}\gets||\mathtt x^{k}-\mathtt x^{k-1}||$
			\ENDWHILE
			\RETURN $\mathtt x^{k+1}$
		}
	\end{algorithmic}\label{algo:jacobi}
\end{algorithm}

\section{Implementation of well-balanced mPDeC-WENO Scheme}\label{se_implemenation}

In the following part, we will describe how the semi-discretization, using the WENO approach, can be written and interpreted as a PDS in order to apply the mPDeC scheme. 
We will see that only small modifications in the WENO code are necessary. 
Furthermore, this interpretation is not only restricted to the WENO procedure but can also applied to other high-order FV/FD discretizations. 
Our WENO approach is working only as a generic example. Finally, the complete algorithm will be presented.

\subsection{WENO and production-destruction systems}
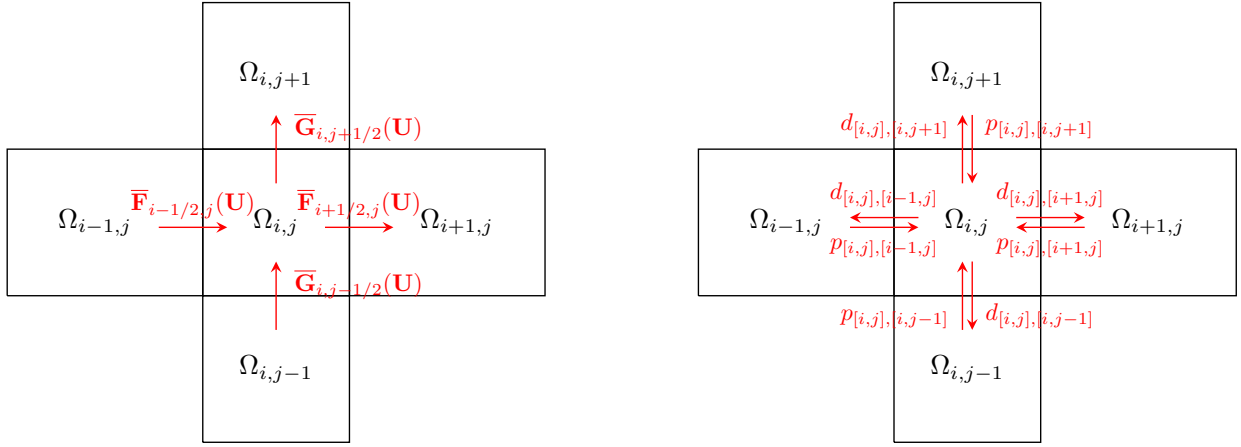
\begin{figure}
\centering
\scalebox{1.3}{
\begin{tikzpicture}
%
%

\draw (0.,2) -- (2,2) -- (2,3.5) -- (0.,3.5) -- (0.,2);
\draw (2,2) -- (3.5,2) -- (3.5,3.5) -- (2,3.5) -- (2,2);
\draw (3.5,2) -- (5.5,2) -- (5.5,3.5) -- (3.5,3.5) -- (3.5,2);
\draw (2,0.5) -- (3.5,0.5) -- (3.5,2) -- (2,2) -- (2,0.5);
\draw (2,3.5) -- (3.5,3.5) -- (3.5,5) -- (2,5) -- (2,3.5);

\node [black,scale=0.8] at (0.9,2.75) {$\Omega_{i-1,j}$};
\node [black,scale=0.8] at (2.75,2.75) {$\Omega_{i,j}$};
\node [black,scale=0.8] at (4.6,2.75) {$\Omega_{i+1,j}$};
\node [black,scale=0.8] at (2.75,4.25) {$\Omega_{i,j+1}$};
\node [black,scale=0.8] at (2.75,1.25) {$\Omega_{i,j-1}$};

\node [red, scale=0.7, anchor=south] at (1.9,2.7) {$\overline{\bF}_{i-1/2,j}(\bU)$};

\draw [-stealth,red] (1.55,2.7) -- (2.25,2.7);

\node [red, scale=0.7, anchor=south] at (3.6,2.7) {$\overline{\bF}_{i+1/2,j}(\bU)$};

\draw [stealth-,red] (3.95,2.7) -- (3.25,2.7);

\node [red, scale=0.7, anchor=base west] at (2.85,2.05) {$\overline{\bG}_{i,j-1/2}(\bU)$};

\draw [-stealth,red] (2.75,1.65) -- (2.75,2.35);

\node [red, scale=0.7, anchor=base west] at (2.85,3.65) {$\overline{\bG}_{i,j+1/2}(\bU)$};

\draw [-stealth,red] (2.75,3.15) -- (2.75,3.85);

\end{tikzpicture}}\hfill
\scalebox{1.3}{
	\begin{tikzpicture}
		%
		%
		
		\draw (0.,2) -- (2,2) -- (2,3.5) -- (0.0,3.5) -- (0.,2);
		\draw (2,2) -- (3.5,2) -- (3.5,3.5) -- (2,3.5) -- (2,2);
		\draw (3.5,2) -- (5.5,2) -- (5.5,3.5) -- (3.5,3.5) -- (3.5,2);
		\draw (2,0.5) -- (3.5,0.5) -- (3.5,2) -- (2,2) -- (2,0.5);
		\draw (2,3.5) -- (3.5,3.5) -- (3.5,5) -- (2,5) -- (2,3.5);
		
		\node [black,scale=0.8] at (0.9,2.75) {$\Omega_{i-1,j}$};
		\node [black,scale=0.8] at (2.75,2.75) {$\Omega_{i,j}$};
		\node [black,scale=0.8] at (4.6,2.75) {$\Omega_{i+1,j}$};
		\node [black,scale=0.8] at (2.75,4.25) {$\Omega_{i,j+1}$};
		\node [black,scale=0.8] at (2.75,1.25) {$\Omega_{i,j-1}$};
		
		\node [red, scale=0.7, anchor=base] at (1.9,2.48) {$p_{[i,j],[i-1,j]}$};
		\node [red, scale=0.7, anchor=south] at (1.9,2.82) {$d_{[i,j],[i-1,j]}$};
		
		\draw [-stealth,red] (1.55,2.7) -- (2.25,2.7);
		\draw [stealth-,red] (1.55,2.8) -- (2.25,2.8);

		\node [red, scale=0.7, anchor=base] at (3.6,2.48) {$p_{[i,j],[i+1,j]}$};
		\node [red, scale=0.7, anchor=south] at (3.6,2.82) {$d_{[i,j],[i+1,j]}$};
		
		\draw [-stealth,red] (3.95,2.7) -- (3.25,2.7);
		\draw [stealth-,red] (3.95,2.8) -- (3.25,2.8);

		\node [red, scale=0.7, anchor=base east] at (2.65,1.75) {$p_{[i,j],[i,j-1]}$};
		\node [red, scale=0.7, anchor=base west] at (2.85,1.75) {$d_{[i,j],[i,j-1]}$};
		
		\draw [-stealth,red] (2.7,1.65) -- (2.7,2.35);
		\draw [stealth-,red] (2.8,1.65) -- (2.8,2.35);

		\node [red, scale=0.7, anchor=base east] at (2.65,3.7) {$d_{[i,j],[i,j+1]}$};
		\node [red, scale=0.7, anchor=base west] at (2.85,3.7) {$p_{[i,j],[i,j+1]}$};
		
		\draw [-stealth,red] (2.7,3.15) -- (2.7,3.85);
		\draw [stealth-,red] (2.8,3.15) -- (2.8,3.85);
		
\end{tikzpicture}}
\caption{Cell $\Omega_{i,j}$, its four neighbors and its production and destruction terms.}\label{fig:omega_ij}
\end{figure}

%
First of all we must underline the fact that in order to preserve the positivity of the water height $h$, 
the mPDeC scheme is going to be applied only to the first equation of system~\eqref{eq:CL} as $h$ is the only variable that must stay nonnegative.
Thus, a simple DeC approach is going to be used to evolve the momentum equations.
So from now on, we shall only talk about the modifications introduced for the first equation to turn it into a production-destruction system.
Given the foundations of finite volume schemes, each control volume has fluxes entering and exiting its boundary and, for each boundary
face, the flux going from element $\alpha=[i,j]$ to element $\beta=[l,r]$ is going to be equal in module and opposite in sign to the flux from element $\beta$ to element $\alpha$. This is the key feature for turning a finite volume schemes into a PDS. Note that the source is zero on $h$ component.\\
Therefore, we can define the production and destruction terms for a general $\alpha=[i,j]$ and the neighboring $\beta=[l,r]$ as 
\begin{equation}\label{eq:PDSforWENO}
	\begin{split}
		p_{[i,j],[i-1,j]}({\bU})= +\frac{1}{\Delta x}\overline{\bF}_{i-1/2,j} (\bU)^+, \quad d_{[i,j],[i-1,j]}({\bU})=- \frac{1}{\Delta x}\overline{\bF}_{i-1/2,j} (\bU)^-,\\
		p_{[i,j],[i+1,j]}({\bU})= -\frac{1}{\Delta x}\overline{\bF}_{\iip,j} (\bU)^-, \quad d_{[i,j],[i+1,j]}({\bU})=+ \frac{1}{\Delta x}\overline{\bF}_{\iip,j} (\bU)^+,\\
		p_{[i,j],[i,j-1]}({\bU})= + \frac{1}{\Delta y}\overline{\bG}_{i,j-1/2} (\bU)^+, \quad d_{[i,j],[i,j-1]}({\bU})=-\frac{1}{\Delta y} \overline{\bG}_{i,j-1/2} (\bU)^-,\\
		p_{[i,j],[i,j+1]}({\bU})= - \frac{1}{\Delta y} \overline{\bG}_{i,j+1/2} (\bU)^-, \quad d_{[i,j],[i,j+1]}({\bU})=+ \frac{1}{\Delta y}\overline{\bG}_{i,j+1/2} (\bU)^+,
	\end{split}
\end{equation}
where with the superscript $\,^+$ and $\,^-$ we denote the positive and the negative part respectively. All the other $p_{\alpha,\beta}$ and $d_{\alpha,\beta}$ not defined here are set to 0. Clearly, this define a conservative and positive PDS, as the properties in Definition \ref{def:POSCONSPDS} are verified. The visualization of the production and destruction terms in Figure~\ref{fig:omega_ij} may help the reader.
We clearly observe that the matrices $(p_{\alpha,\beta})$ and $(d_{\alpha,\beta})$ are $S\times S$ sparse matrices, with, at most, 4 nonzero entries per row and $S=N_x\cdot  N_y$.

%
Once the production and destruction matrices have been assembled, the next step consist in running the mPDeC algorithm in Eq.~\eqref{eq:explicit_dec_correction}. Note that the matrix built in \eqref{eq:matrixMPDeC} is sparse as well with at most 5 nonzero entries for each row (4 nonzero entries of the production/destruction terms and the diagonal term).
Hence, in the numerical computations we will use the classical Jacobi iterative method, see Algorithm \ref{algo:jacobi}, to obtain the solution of system \eqref{eq:explicit_dec_correction} at each iteration. 
Indeed, it is provable \cite{offner2020arbitrary,burchard2003high} that the Jacobi iteration algorithm converges when applied on the matrix defined in \eqref{eq:matrixMPDeC}. 

In all calculations, we set the tolerance to machine precision and the algorithm converge towards the solution in few iterations. Experimentally, we have seen that usually 10-20 iterations suffice, in the worst cases 40 iterations are needed and the number of iterations do not depend much on the mesh size. 
Overall, considering the assembly the mass matrix, the inversion of the system with Jacobi iterations, the whole mPDeC procedure increases the computational costs of around 18\% with respect to the original DeC algorithm using the same CFL. The code has not been construct with the goal of optimizing all the procedures, so it might well be that this extra computational cost can be decreased with better implementations. The greatest advantage is anyway that no CFL restrictions are required in order to guarantee the positivity of the solution. So, a CFL$=1$ suffices to guarantee the stability of the method.
In the numerical test section, we will also detail the number of Jacobi iterations and the overhead computational cost that the method brings into the system with respect to the explicit method.

\begin{remark}[Efficiency and properties of Jacobi iterative method]
	Though being a very simple algorithm, Jacobi iterative method is particularly effective for this application. Indeed, the mass matrix in $\R^{S\times S}$ is very sparse, i.e., only $5S$ non-zero elements, hence, at each iterations, only $5S$ multiplications are computed. Moreover, experimentally, we have seen that the overhead computational cost of the Jacobi solver is around 10\% of the whole computational cost and that the iterations needed are usually below 20, see section~\ref{se:numierics}. Moreover, Jacobi guarantees the positivity of the solution at every iteration of the procedure.\\		
	Clearly, there are other many iterative methods to solve linear system \cite{saad2003iterative}, e.g. Krylov preconditioning methods. All of these methods have larger complexity per iteration, but may converge faster to the solution of the system. Nevertheless, for most of these methods it is not possible to guarantee the positivity of the solution of the iterative solver along the process.
	That is why, we will use the Jacobi iterative method.
\end{remark}
%
%
\subsection{Full Algorithm}

After describing how the WENO (FV) procedure can be interpreted as a PDS, we give a 
more precise description of the full algorithm which is used to calculate the numerical solution.
%
%

In our version of this approach, we have to adapt the steps taking into account the re-interpretation of the WENO approach as a production destruction system and the well balancing approach and this is described in Algorithm \ref{algo:WENO}.

Finally, we can combine all the ingredients described above in a full algorithm as in Algorithm \ref{algo:fullAlgo}. There, we simply use the mPDeC to evolve in time and the production and destruction functions are given by the WENO description from above.
%
%
%

\begin{algorithm}
	\begin{algorithmic}[1] 
		\REQUIRE $\bU_{i,j}$, well balanced fluxes
		\STATE Reconstruct on the quadrature points on cell interfaces the variable $\bU$ in a high order fashion
		\STATE Compute the numerical fluxes at quadrature points on cell interfaces $\hat \bF (\bu^L, \bu^R)$ 
		\STATE Subtract the correction for well balanced problems and obtain $\overline{\bF}$ and $\overline{\bG}$
		\STATE Integrate over the cell interface to obtain the numerical fluxes $\bF_{\iip,j},\bF_{\iin,j}$ and $\bG_{i,\jjn}, \bG_{i,\jjp}$
		\STATE Compute $p_{\alpha,\beta}(\bU),\,d_{\alpha,\beta}(\bU)$ as in \eqref{eq:PDSforWENO}
		\RETURN $p_{\alpha,\beta}(\bU),\,d_{\alpha,\beta}(\bU)$
	\end{algorithmic}
	\caption{WENO FV with PDS structure}\label{algo:WENO}
\end{algorithm}

\begin{algorithm}
	\begin{algorithmic}[1]
		\REQUIRE $\bU_{i,j}^0$, $T$
		\STATE $t=0$
		\WHILE{$t<T$}
		\STATE Compute $\Delta t$ by CFL restrictions
		\STATE $\bU^{n+1}$=DeC($\bU^n,\Delta t$, WENOPDS) with DeC Algorithm \eqref{DeC_method} where update formula \eqref{eq:explicit_dec_correction} is used for $h$ and \eqref{eq:explicit_dec_explicit} is used for $hu$ and $hv$ and the WENO PSD function are given by Algorithm \ref{algo:WENO}
		\STATE $t=t+\Delta t$
		\ENDWHILE
	\end{algorithmic}
\caption{Full algorithm}\label{algo:fullAlgo}
\end{algorithm}

 \begin{conclusion}
 	The described method is high order accurate, positivity preserving, conservative, non-oscillatory and well-balanced for the shallow water equation. 
 \end{conclusion}

\begin{remark}[Difference with respect to classical WENO]
	We want to highlight the differences between a classical WENO and the proposed algorithm are minimal. Indeed, once the spatial discretization is performed with a simple WENO step we need to apply two easy modifications. The first one consists in subtracting the flux related to the steady state variables. The second one consists in the definition of the production and destruction terms. Then the mPDeC can be applied as a simple time integration scheme.
	The code and these modifications are available at the reproducibility repository \cite{ourrepo}.
\end{remark} 

\begin{remark}[Advantages of mPDeC]
	We shall remark that the presented method does not require any CFL constraint to obtain positive solutions for $h$, while the classical positivity limiter for WENO5 requires a CFL number of $1/12\approx0.083$. Clearly, a CFL number for a classical explicit method must be anyway used (between 1 and 1.5), but this allow to run the simulation with much less time steps than a classical explicit WENO scheme, with the extra computational cost of the Jacobi iterative method, which is negligible with respect the cost of decreasing the CFL number of factor of 12.
	Moreover, the procedure allows to have a provably positive method with arbitrarily high order of accuracy, while with positive limiters applied to WENO schemes only SSPRK methods guarantee the positivity of the solutions and they exist only up to order 4 \cite{gottlieb2011strong}. 
\end{remark}
 
 \begin{remark}[Stability and accuracy of combined time-integration methods]
We combined two time-integration methods, the mPDeC and the DeC approach. How can we assess the stability and accuracy properties the combined method. For the accuracy of the combined scheme, we can exclude a loss of accuracy
as long as the two methods are consistent up to the same order of accuracy, even if applied separately to each equation.
For instance, one could perform a simple Taylor expansion to prove the error behavior for each component and also in the combination step-after-step. In practice, as done in \cite{offner2020arbitrary}, it is sufficient to develop the Taylor expansion on the modified Patankar weighting coefficients to obtain a high order approximation of the DeC.
This is also verified in our numerical tests. 
Stability is more problematic since it is not clear how to analyze it even for modified Patankar schemes in the ODE case, where only few preliminary works are available \cite{izgin2022lyapunov, torlo2021stability}. This is also due to the nonlinear character of the scheme itself.
Therefore, a stability investigation for PDEs is far beyond the goals of  this paper. 
However, in general explicit methods have a more restrictive stability region than implicit methods.  Therefore, we will simple assume  that our stability region is determined by the 
underlying explicit method (DeC). 
\end{remark}

\section{Numerical Simulations}\label{se:numierics}

The goal of this section is to present the results obtained with the fifth order positivity-preserving 
mPDeC scheme, compared to that given by the classical fifth order DeC time integration method.
The first test case consists in assessing the convergence properties of the spatial and temporal discretization
on an unsteady vortex-type solution~\cite{ricchiuto2021analytical}. Afterwards, we focus on testing the well-balanced
implementation for the lake at rest solution by showing its impact on perturbation analysis. 
Finally, three challenging simulations are performed to prove its capabilities to cope with wet-dry fronts. 
It is important to underline that in all simulations, especially those involving shocks, when there is no comparison with DeC5 (dam break problems),
the mPDeC5 method obtains consistent results with respect to classical methods.
For all the simulations carried out herein periodic boundary conditions have been considered together with the local Lax-Friedrichs (Rusanov) numerical flux.

\subsection{Unsteady vortex}

\begin{figure}
	\centering
	\subfigure[Computational mesh]{\includegraphics[width=0.45\textwidth]{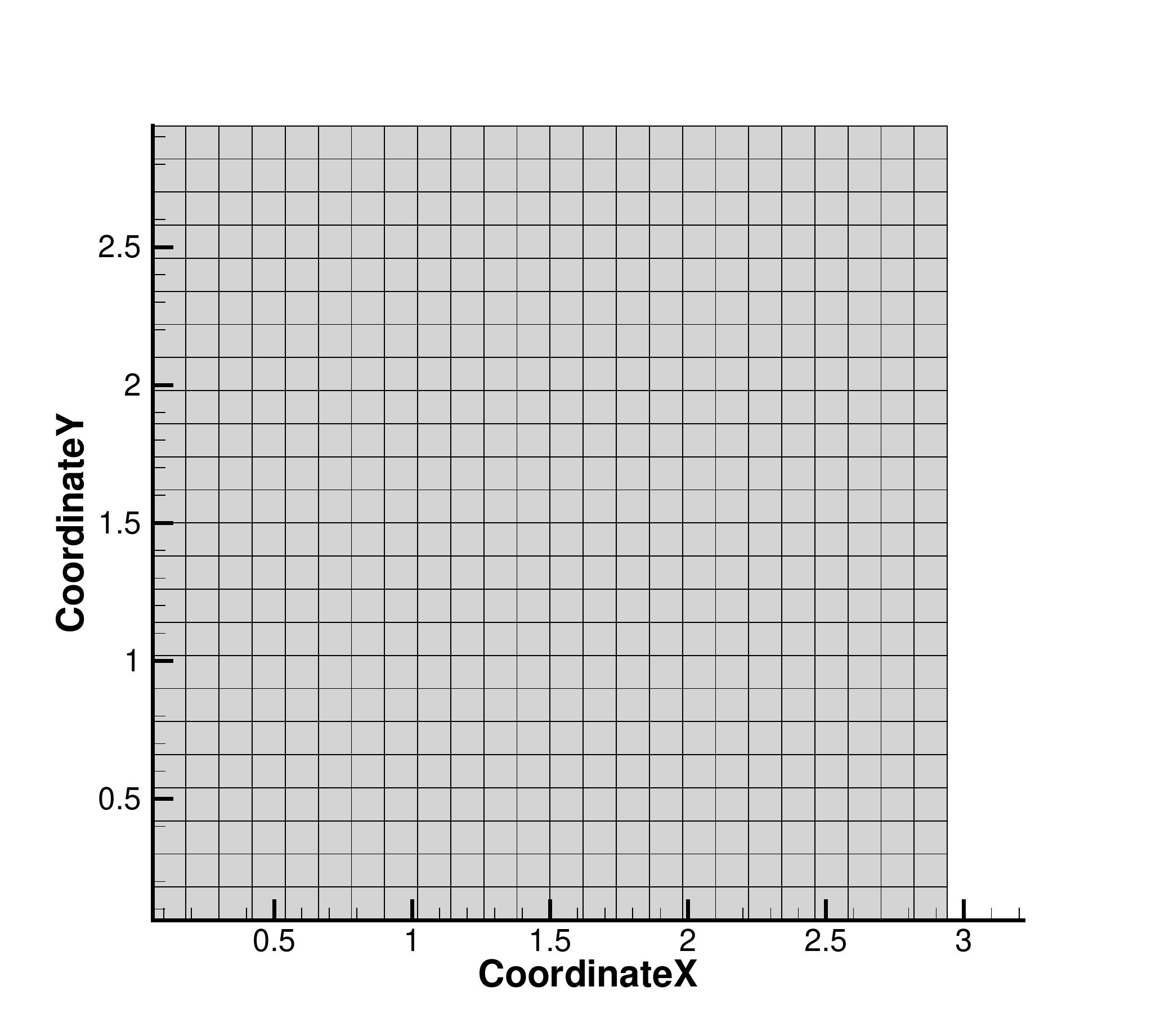}}
	\subfigure[Exact solution]{\includegraphics[width=0.45\textwidth]{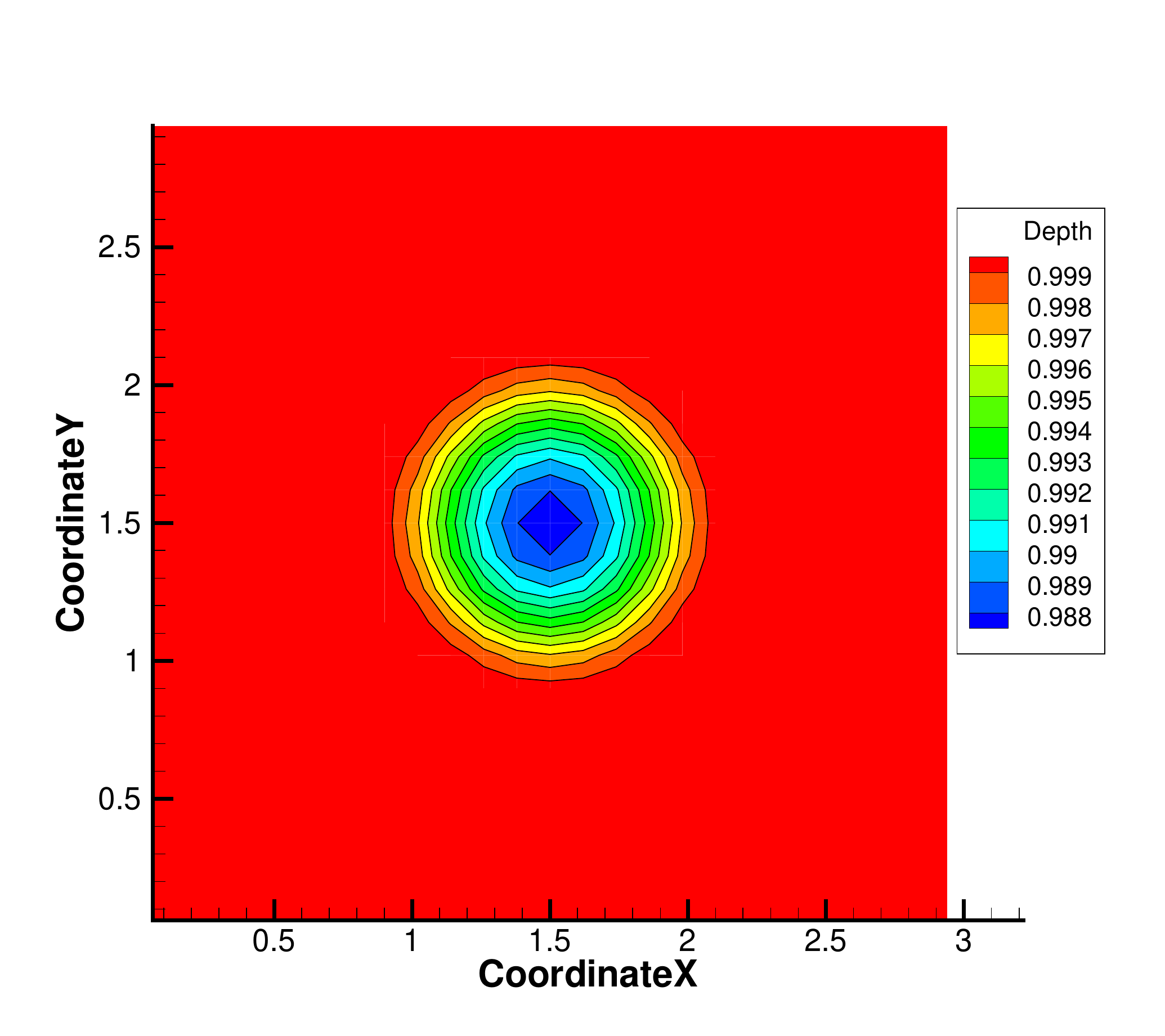}}
	\caption{Unsteady vortex: test case setting.}\label{fig:vortex_set}
\end{figure}

In order to verify the order of accuracy we consider a moving smooth vortex. The computational domain is the square
$[0,3]\times[0,3]$. The initial condition is 
given by some perturbations $\delta$ applied on a homogeneous background field $(h_0,u_0,v_0)=(1,2,3)$. Hence,
the perturbation for the depth variable $h$ is
\begin{equation}\label{eq:movingvortexh}
h(r)=h_0 - \delta h(r) = h_0 - \gamma \begin{cases} e^{-\frac{1}{\arctan^3(1-r^2)}}, & \text{ if } r < 1, \\
0, & \text{ else}, \end{cases}\qquad \text{with }r=\sqrt{(x-1.5)^2+(y-1.5)^2}
\end{equation}
and the vortex amplitude is $\gamma=0.1$. The velocity field is affected by the following perturbations
\begin{equation}\label{eq:movingvortexu}
\left(\begin{array}{c} \delta u \\ \delta v \end{array}\right) = \sqrt{2\,g\,\partial_r h} \left(\begin{array}{c} (y-1.5) \\ -(x-1.5) \end{array}\right),
\end{equation}
where $\partial_r h=\partial_r h(r)$ is a function of the radial distance from the center of the vortex
\begin{equation}
\partial_r h(r) = \,\frac{3\,\gamma\, e^{-\frac{1}{\arctan^3(1-r^2)}}}{\arctan^4(r^2-1)((r^2-1)^2+1)}.
\end{equation}
It is important to highlight the fact that this solution is $\mathcal{C}^\infty$, which is a fundamental property 
for testing arbitrarily high order schemes. Many vortex-type solutions
can be found available online but most of them can only be used to test lower order schemes. The exact solution of this problem is given by
\begin{equation}
	h(x,y,t) = h(x-u_0t,y-v_0t,0), \quad u(x,y,t) = u(x-u_0t,y-v_0t,0), \quad v(x,y,t) = v(x-u_0t,y-v_0t,0).
\end{equation}

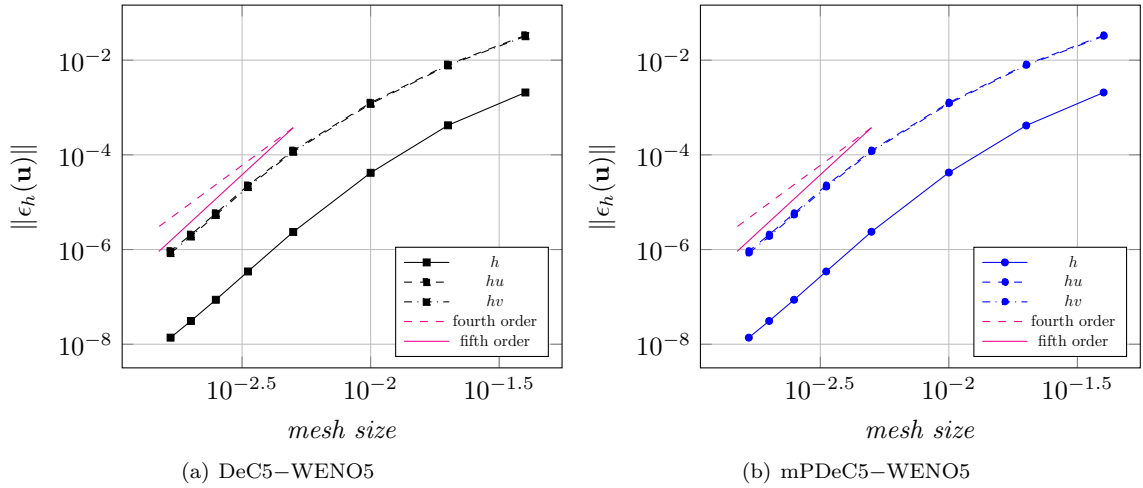
\begin{figure}
\centering
\subfigure[DeC5$-$WENO5]{
\begin{tikzpicture}
    \begin{axis}[
     xmode=log, ymode=log,
     grid=major,
     xlabel={\textit{mesh size}},
     ylabel={$\|\epsilon_h(\mathbf{u})\|$},
     xlabel shift = 1 pt,
     ylabel shift = 1 pt,
     legend pos= south east,
     legend style={nodes={scale=0.6, transform shape}},
     width=.45\textwidth
     ]
     \addplot[mark=square*,mark size=1.3pt,black] table [y=h, x=size]{vortex_Dec5.dat};
     \addlegendentry{$h$}

     \addplot[mark=square*,mark size=1.3pt,black,dashed] table [y=hu, x=size]{vortex_Dec5.dat};
     \addlegendentry{$hu$}

     \addplot[mark=square*,mark size=1.3pt,black,dashdotted] table [y=hv, x=size]{vortex_Dec5.dat};
     \addlegendentry{$hv$}
     \addplot[magenta,domain=0.005:0.0015, dashed]{x*x*x*x*(600000)};
     \addlegendentry{fourth order}
     \addplot[magenta,domain=0.005:0.0015]{x*x*x*x*x*(120000000)};
     \addlegendentry{fifth order}

     \end{axis}
\end{tikzpicture}}
\subfigure[mPDeC5$-$WENO5]{
\begin{tikzpicture}
    \begin{axis}[
     xmode=log, ymode=log,
     grid=major,
     xlabel={\textit{mesh size}},
     ylabel={$\|\epsilon_h(\mathbf{u})\|$},
     xlabel shift = 1 pt,
     ylabel shift = 1 pt,
     legend pos= south east,
     legend style={nodes={scale=0.6, transform shape}},
     width=.45\textwidth
     ]
     \addplot[mark=otimes*,mark size=1.3pt,blue] table [y=h, x=size]{vortex_mPDec5.dat};
     \addlegendentry{$h$}

     \addplot[mark=otimes*,mark size=1.3pt,blue,dashed] table [y=hu, x=size]{vortex_mPDec5.dat};
     \addlegendentry{$hu$}

     \addplot[mark=otimes*,mark size=1.3pt,blue,dashdotted] table [y=hv, x=size]{vortex_mPDec5.dat};
     \addlegendentry{$hv$}
     \addplot[magenta,domain=0.005:0.0015, dashed]{x*x*x*x*(600000)};
     \addlegendentry{fourth order}
     \addplot[magenta,domain=0.005:0.0015]{x*x*x*x*x*(120000000)};
     \addlegendentry{fifth order}

     \end{axis}
\end{tikzpicture}}
\caption{Unsteady vortex: convergence tests.}\label{fig:vortex_conv}
\end{figure}

For the unsteady vortex, two convergence tests are run for WENO5 coupled with both the time integration schemes DeC5 and mPDeC
to corroborate the fact that, for smooth flows, the results should be almost identical. We used CFL=0.7.
The convergence tests are run on cartesian meshes of size $25\times25$, $50\times50$, $100\times100$, $200\times200$, $300\times300$, $400\times400$, $500\times500$ and $600\times600$.
The computational mesh of the coarsest grid and initial condition for this test case are shown in Figure~\ref{fig:vortex_set}. 
For these convergence tests we used a tolerance $\varepsilon=10^{-30}$ both for the positivity limiter and the mPDeC divisions, since the errors that we obtain are of the order of $10^{-8}$.
The error $||\epsilon_h(\bu)||$ is the $\mathbb L^1$ norm of the difference between the exact solution and the approximated one.
Figure~\ref{fig:vortex_conv} points out the predicted fifth order behavior for both time integration schemes.

\begin{figure}
		\centering
		\subfigure[Ratio of computational time of mPDeC over DeC\label{fig:vortex_time_ratio}]{
		\begin{tikzpicture}
	\begin{axis}[
		xmode=log,
		grid=major,
		xlabel={$N_x$},
		xlabel shift = 1 pt,
		width=.4\textwidth
		]
		\addplot[mark=o,mark size=1.3pt,black]             table [y=time_ratio   , x=N]{vortex_time/CFL01_comparison.txt};
		\addplot[mark=square,mark size=1.3pt,black]             table [y=time_ratio   , x=N]{vortex_time/CFL02_comparison.txt};
		\addplot[mark=triangle,mark size=1.3pt,black]             table [y=time_ratio   , x=N]{vortex_time/CFL04_comparison.txt};
		\addplot[mark=diamond,mark size=1.3pt,black]             table [y=time_ratio   , x=N]{vortex_time/CFL07_comparison.txt};
	\end{axis}
\end{tikzpicture}}
\subfigure[Computational time and error\label{fig:vortex_time_error}]{\begin{tikzpicture}
			\begin{axis}[
				xmode=log, ymode=log,
				grid=major,
				xlabel={Time},
				ylabel={$\|\epsilon(h)\|$},
				xlabel shift = 1 pt,
				ylabel shift = 1 pt,
				legend pos= outer north east,
				legend style={nodes={scale=0.6, transform shape}},
				width=.4\textwidth
				]
				\addplot[mark=o,mark size=1.3pt,black]             table [y=errorhDeC   , x=timeDeC]{vortex_time/CFL01_comparison.txt};
				\addlegendentry{DeC5-WENO5 CFL=0.1}
				\addplot[mark=o*, dashed, mark size=1.3pt,red]             table [y=errorhmPDeC   , x=timemPDeC]{vortex_time/CFL01_comparison.txt};
				\addlegendentry{mPDeC5-WENO5 CFL=0.1}
				\addplot[mark=square,mark size=1.3pt,black]             table [y=errorhDeC   , x=timeDeC]{vortex_time/CFL02_comparison.txt};
				\addlegendentry{DeC5-WENO5 CFL=0.2}
				\addplot[dashed,mark=square*, mark size=1.3pt,red]             table [y=errorhmPDeC   , x=timemPDeC]{vortex_time/CFL02_comparison.txt};
				\addlegendentry{mPDeC5-WENO5 CFL=0.2}
				\addplot[mark=triangle,mark size=1.3pt,black]             table [y=errorhDeC   , x=timeDeC]{vortex_time/CFL04_comparison.txt};
				\addlegendentry{DeC5-WENO5 CFL=0.4}
				\addplot[mark=triangle*,dashed, mark size=1.3pt,red]             table [y=errorhmPDeC   , x=timemPDeC]{vortex_time/CFL04_comparison.txt};
				\addlegendentry{mPDeC5-WENO5 CFL=0.4}
				\addplot[mark=diamond,mark size=1.3pt,black]             table [y=errorhDeC   , x=timeDeC]{vortex_time/CFL07_comparison.txt};
				\addlegendentry{DeC5-WENO5 CFL=0.7}
				\addplot[mark=diamond*, mark size=1.3pt,red,dashed]             table [y=errorhmPDeC   , x=timemPDeC]{vortex_time/CFL07_comparison.txt};
				\addlegendentry{mPDeC5-WENO5 CFL=0.7}
			\end{axis}
	\end{tikzpicture}}
\subfigure[Average of Jacobi iterations and confidence interval\label{fig:vortex_jacobi}]{\begin{tikzpicture}
	\begin{axis}[
		xmode=log,
		grid=major,
		xlabel={Mesh elements in $x$},
		ylabel={Jacobi iterations},
		xlabel shift = 1 pt,
		ylabel shift = 1 pt,
		legend pos= outer north east,
		legend style={nodes={scale=0.6, transform shape}},
		width=.45\textwidth
		]		
		\addplot[mark=star,dashed, mark size=1.3pt,blue]             table [y=jacAve   , x=N]{vortex_time/CFL07_comparison.txt};
		\addlegendentry{CFL=0.7}
		\addplot[mark=diamond*,dotted,mark size=1.3pt,red]table [y=jacAve   , x=N]{vortex_time/CFL04_comparison.txt};
		\addlegendentry{CFL=0.4}			
		\addplot[mark=square*,dotted,mark size=1.3pt,darkspringgreen]table [y=jacAve   , x=N]{vortex_time/CFL02_comparison.txt};
		\addlegendentry{CFL=0.2}	
		\addplot[mark=triangle*,dashdotted,mark size=1.3pt,magenta]table [y=jacAve   , x=N]{vortex_time/CFL01_comparison.txt};
		\addlegendentry{CFL=0.1}
		\addplot[name path=us_top,dashed,blue!50]             table [y expr=\thisrow{jacAve}+0.5*\thisrow{jacStd}   , x=N]{vortex_time/CFL07_comparison.txt};
		\addplot[name path=us_bot,dashed, blue!50]             table [y expr=\thisrow{jacAve}-0.5*\thisrow{jacStd}   , x=N]{vortex_time/CFL07_comparison.txt};
		\addplot[blue!30,fill opacity=0.3] fill between[of=us_top and us_bot];
		\addplot[name path=us_top,dotted,red!50]             table [y expr=\thisrow{jacAve}+0.5*\thisrow{jacStd}   , x=N]{vortex_time/CFL04_comparison.txt};
		\addplot[name path=us_bot,dotted, red!50]             table [y expr=\thisrow{jacAve}-0.5*\thisrow{jacStd}   , x=N]{vortex_time/CFL04_comparison.txt};
		\addplot[red!30,fill opacity=0.3] fill between[of=us_top and us_bot];		
		\addplot[name path=us_top,dotted,darkspringgreen!50]             table [y expr=\thisrow{jacAve}+0.5*\thisrow{jacStd}   , x=N]{vortex_time/CFL02_comparison.txt};
		\addplot[name path=us_bot,dotted, darkspringgreen!50]             table [y expr=\thisrow{jacAve}-0.5*\thisrow{jacStd}   , x=N]{vortex_time/CFL02_comparison.txt};
		\addplot[darkspringgreen!30,fill opacity=0.3] fill between[of=us_top and us_bot];		
		\addplot[name path=us_top,dashdotted,magenta!50]             table [y expr=\thisrow{jacAve}+0.5*\thisrow{jacStd}   , x=N]{vortex_time/CFL01_comparison.txt};
		\addplot[name path=us_bot,dashdotted, magenta!50]             table [y expr=\thisrow{jacAve}-0.5*\thisrow{jacStd}   , x=N]{vortex_time/CFL01_comparison.txt};
		\addplot[magenta!30,fill opacity=0.3] fill between[of=us_top and us_bot];
	\end{axis}
\end{tikzpicture}}
\caption{Unsteady vortex test: computational time and Jacobi iterations \label{fig:vortex_time}}
\end{figure}
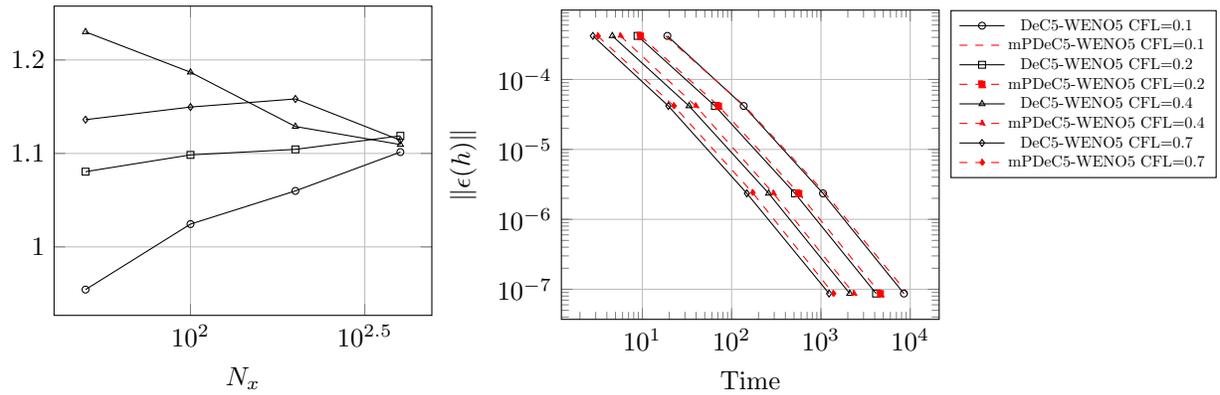
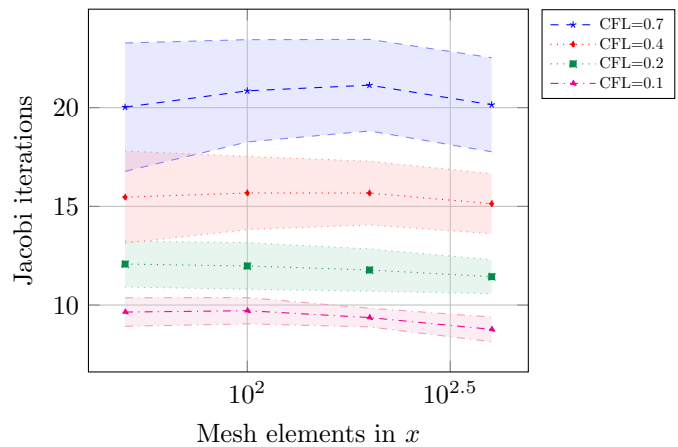

In Figure~\ref{fig:vortex_time} we study the computational costs of the two methods (mPDeC5 and DeC5) with respect to different CFL numbers and mesh refinements. In Figure~\ref{fig:vortex_time_ratio} we plot the ratio of computational time of mPDeC5 over the computational time of DeC5 needed to finish the simulation for the same CFL and mesh. We see that for fine meshes, when the computational times are more reliable, all ratios are close to 1.1. This means that the overhead that mPDeC5 requires, with respect to an explicit method, is of about 10\%. In Figure~\ref{fig:vortex_time_error} we plot errors with respect to computational times and we see a very small difference between mPDeC5 and DeC5, while there is a huge difference in computational costs when changing the CFL. The mPDeC5 is guaranteed to run at any $\text{CFL}^{\text{FE}}<1$, while the positivity for other SSPRK methods is guaranteed only for $\text{CFL}^{\text{FE}}<1/12$ with WENO5. Hence, there is huge advantage with mPDeC.
Finally, in Figure~\ref{fig:vortex_jacobi} we can observe on average how many iterations are needed to solve the linear system and a confidence interval defined by the average $\pm$ $\frac12$ standard deviation. It is clear that Jacobi iterations are mainly driven by the CFL number, which explicitly appear as a coefficient of the mass matrix minus the identity. 
Indeed, all factors $\dt \frac{d}{h}$ or $\dt \frac{p}{h}$ are linearly dependent on the CFL as production and destruction terms are proportional to $\frac{1}{\Delta x}$.

\subsection{Lake at rest}

\begin{figure}
\centering
\subfigure[Computational mesh]{\includegraphics[width=0.45\textwidth]{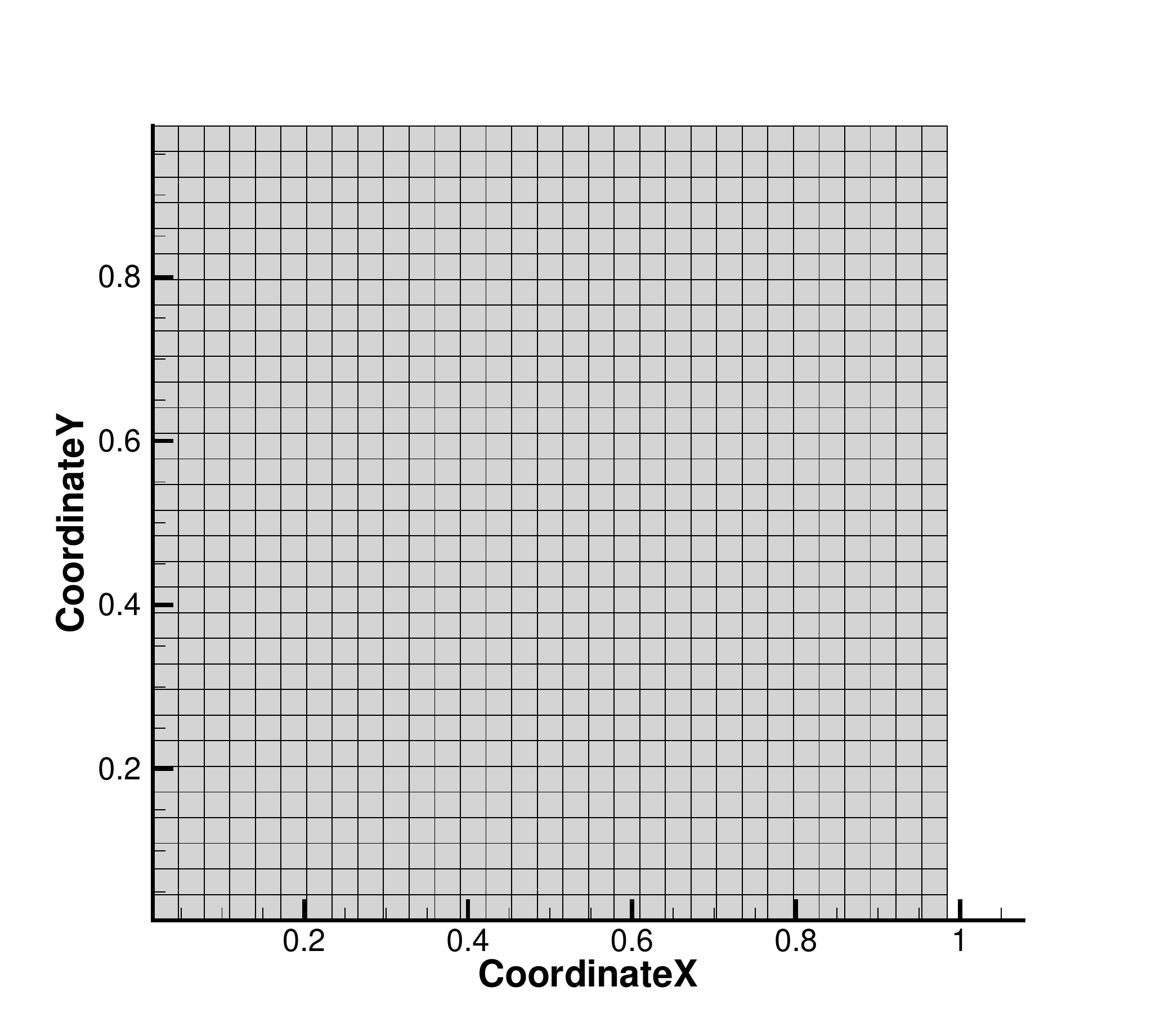}}
\subfigure[Exact solution of $h$]{\includegraphics[width=0.45\textwidth]{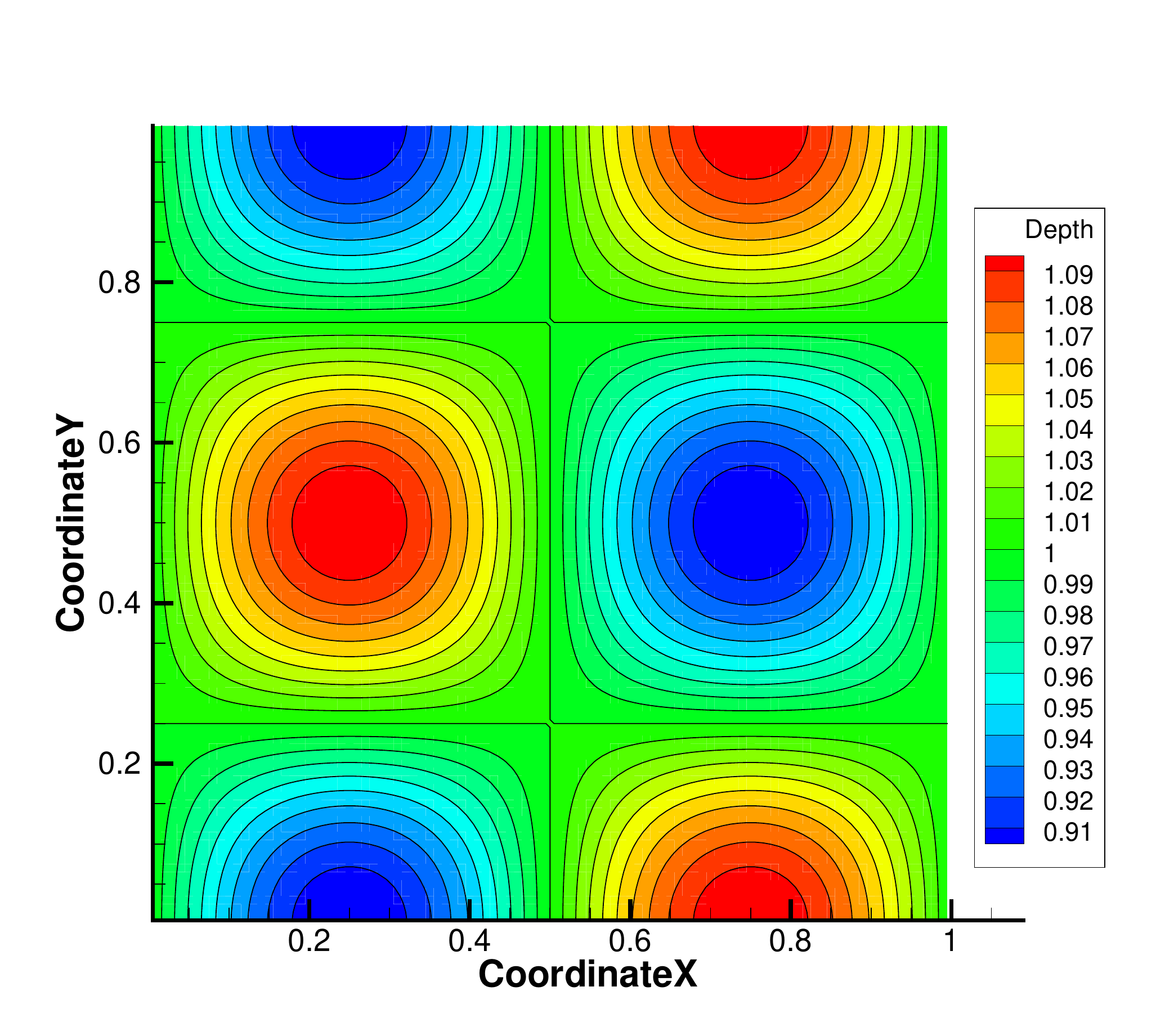}}
\caption{Lake at rest: test case setting.}\label{fig:lakeatrest_set}
\end{figure}

\begin{figure}
\centering
\subfigure[DeC5--WENO5]{
\begin{tikzpicture}
    \begin{axis}[
     xmode=log, ymode=log,
     xmin=0.003,xmax=0.3,
     grid=major,
     xlabel={\textit{mesh size}},
     ylabel={$\|\epsilon_h(\mathbf{u})\|$},
     xlabel shift = 1 pt,
     ylabel shift = 1 pt,
     legend pos= south east,
     legend style={nodes={scale=0.6, transform shape}},
     width=.45\textwidth
     ]
     \addplot[mark=square*,mark size=1.3pt,black] table [y=h, x=size]{Lake_Dec5NoWB.dat};
     \addlegendentry{$h$}

     \addplot[mark=square*,mark size=1.3pt,black,dashed] table [y=hu, x=size]{Lake_Dec5NoWB.dat};
     \addlegendentry{$hu$}

     \addplot[mark=square*,mark size=1.3pt,black,dashdotted] table [y=hv, x=size]{Lake_Dec5NoWB.dat};
     \addlegendentry{$hv$}
     \addplot[magenta,domain=0.2:0.0075, dashdotted]{x*x*x*x*(200)};
     \addlegendentry{fourth order}
     \addplot[magenta,domain=0.2:0.0075]{x*x*x*x*x*(1000)};
     \addlegendentry{fifth order}

     \end{axis}
\end{tikzpicture}}
\subfigure[mPDeC5--WENO5]{
\begin{tikzpicture}
    \begin{axis}[
     xmode=log, ymode=log,
     xmin=0.003,xmax=0.3,
     grid=major,
     xlabel={\textit{mesh size}},
     ylabel={$\|\epsilon_h(\mathbf{u})\|$},
     xlabel shift = 1 pt,
     ylabel shift = 1 pt,
     legend pos= south east,
     legend style={nodes={scale=0.6, transform shape}},
     width=.45\textwidth
     ]
     \addplot[mark=otimes*,mark size=1.3pt,blue] table [y=h, x=size]{Lake_mPDec5NoWB.dat};
     \addlegendentry{$h$}

     \addplot[mark=otimes*,mark size=1.3pt,blue,dashed] table [y=hu, x=size]{Lake_mPDec5NoWB.dat};
     \addlegendentry{$hu$}

     \addplot[mark=otimes*,mark size=1.3pt,blue,dashdotted] table [y=hv, x=size]{Lake_mPDec5NoWB.dat};
     \addlegendentry{$hv$}
     \addplot[magenta,domain=0.2:0.0075, dashdotted]{x*x*x*x*(200)};
     \addlegendentry{fourth order}
     \addplot[magenta,domain=0.2:0.0075]{x*x*x*x*x*(1000)};
     \addlegendentry{fifth order}

     \end{axis}
\end{tikzpicture}}
\caption{Lake at rest: convergence tests without preserving the exact solution.}\label{fig:lakeatrest_conv}
\end{figure}
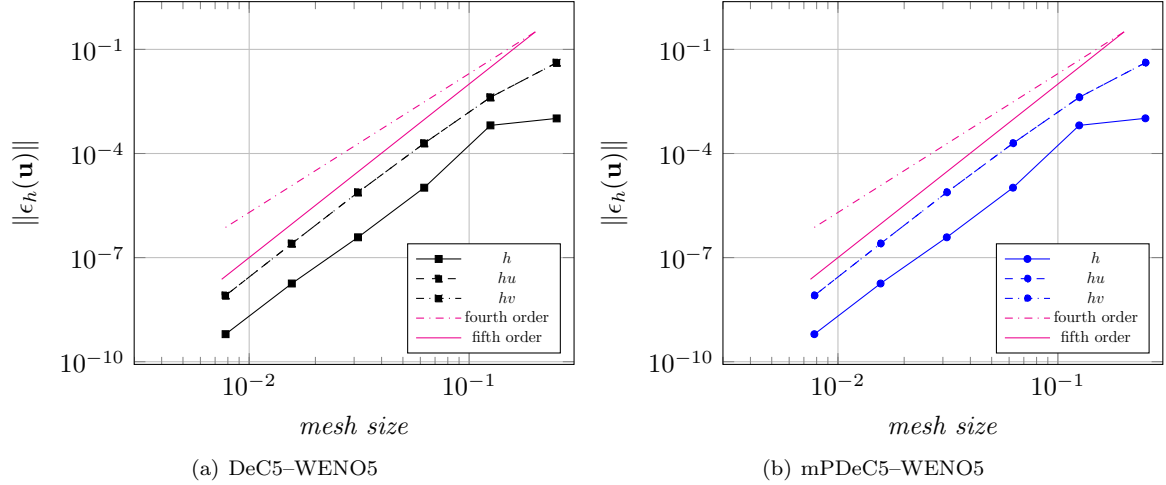

As already introduced in the theoretical part, this test case is needed to prove the presented scheme is well-balanced. The computational domain is the square $[0,1]\times[0,1]$ and the steady solution 
of this problem and bathymetry are briefly summarized below
\begin{equation}
b(x,y) = 0.1\,\sin\left(2\,\pi\,x\right)\,\cos\left(2\,\pi\,y\right),\qquad h(x,y) = 1\,-\,b(x,y),\qquad u=v=0.
\end{equation}
This benchmark can also be used to test once again the order of accuracy of our discretization. Indeed, we expect the method 
to converge with a fifth order slope when not well-balanced and we expect machine precision errors for all the well-balanced tests. This simulation has been performed with four different settings: 
DeC5 and mPDeC5, well-balanced and not well-balanced. For all cases, we employed a fifth order WENO discretization for the spatial derivatives and CFL=0.9. 
Also for this test, we chose $\varepsilon=10^{-30}$ to check the accuracy of the scheme with errors of the order of $10^{-10}$.
As expected, the error computed for the well-balanced simulations is exactly zero therefore we did not plot them along with the other results.
The interested reader can run the simulations and test the properties of our method by downloading the code available at the reproducibility repository \cite{ourrepo}. 
The Cartesian mesh employed for this convergence test are $4\times4$, $8\times8$, $16\times16$,  $32\times32$,  $64\times64$ and $128\times128$.
The exact solution is presented in Figure~\ref{fig:lakeatrest_set}, along with the $32\times32$ mesh.  
As can be noticed from Figure~\ref{fig:lakeatrest_conv}, mPDeC5 allows a fifth order convergence rate as theoretically proved with
results almost identical to those of DeC5.

\subsection{Wet-dry lake at rest}
Now we test the capability of dealing with wet and dry regions of the scheme in a very simple context. We consider a bathymetry given by a bump 
\begin{equation}\label{eq:bathymetry_pert}
	b(x,y) = \begin{cases} e^{1-\frac{1}{1-r^2}}, & \text{ if } r^2<1, \\ 0,  & \text{  else  }, \end{cases}\qquad\text{ where } r^2=x^2+y^2,
\end{equation}
on the domain $[-5,5]\times[-2,2]$.
The lake at rest solution is 
\begin{equation}
	h(x,y) = \max\lbrace 0.7 - b(x,y),0 \rbrace ,\qquad u=v=0.
\end{equation}
The maximum of the bathymetry is 1, hence, there is a dry island at the center of the domain.
For practical purposed, we set the initial conditions to be
\begin{equation}
	h_0(x,y) = \max\lbrace 0.7 - b(x,y),\varepsilon \rbrace ,\qquad u=v=0,
\end{equation} 
with $\varepsilon=10^{-6}$. We consider final time $T=1$.\\
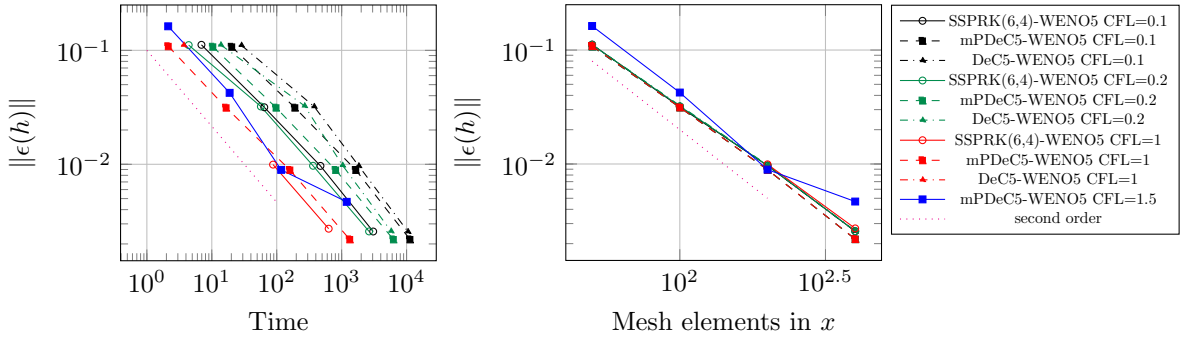
\begin{figure}
	\centering
	\begin{tikzpicture}
		\begin{axis}[
			xmode=log, ymode=log,
			grid=major,
			xlabel={Time},
			ylabel={$\|\epsilon(h)\|$},
			xlabel shift = 1 pt,
			ylabel shift = 1 pt,
			width=.35\textwidth
			]
			\addplot[mark=o,mark size=1.3pt,black]             table [y=errh   , x=time]{IslandAtRest/SSPRK/CFL01/convergence.txt};
			\addplot[mark=square*,dashed,mark size=1.3pt,black]table [y=errh   , x=time]{IslandAtRest/mPDeC/CFL01/convergence.txt};
			\addplot[mark=triangle*,dashdotted,mark size=1.3pt,black]table [y=errh   , x=time]{IslandAtRest/DeC/CFL01/convergence.txt};
			\addplot[mark=o,mark size=1.3pt,darkspringgreen]             table [y=errh   , x=time]{IslandAtRest/SSPRK/CFL02/convergence.txt};
			\addplot[mark=square*,dashed,mark size=1.3pt,darkspringgreen]table [y=errh   , x=time]{IslandAtRest/mPDeC/CFL02/convergence.txt};
			\addplot[mark=triangle*,dashdotted,mark size=1.3pt,darkspringgreen]table [y=errh   , x=time]{IslandAtRest/DeC/CFL02/convergence.txt};
			\addplot[mark=o,mark size=1.3pt,red]             table [y=errh   , x=time]{IslandAtRest/SSPRK/CFL1/convergence.txt};
			\addplot[mark=square*,dashed,mark size=1.3pt,red]table [y=errh   , x=time]{IslandAtRest/mPDeC/CFL1/convergence.txt};
			\addplot[mark=triangle*,dashdotted,mark size=1.3pt,red]table [y=errh   , x=time]{IslandAtRest/DeC/CFL1/convergence.txt};
			\addplot[mark=square*,mark size=1.3pt,blue]             table [y=errh   , x=time]{IslandAtRest/mPDeC/CFL15/convergence.txt};
			\addplot[magenta,domain=1:100, dotted]{0.1*x^(-2./3.)};
		\end{axis}
	\end{tikzpicture}
	\begin{tikzpicture}
		\begin{axis}[
			xmode=log, ymode=log,
			grid=major,
			xlabel={Mesh elements in $x$},
			ylabel={$\|\epsilon(h)\|$},
			xlabel shift = 1 pt,
			ylabel shift = 1 pt,
			legend pos= outer north east,
			legend style={nodes={scale=0.6, transform shape}},
			width=.35\textwidth
			]
			\addplot[mark=o,mark size=1.3pt,black]             table [y=errh   , x=Nx]{IslandAtRest/SSPRK/CFL01/convergence.txt};
			\addlegendentry{SSPRK(6,4)-WENO5 CFL=0.1}
			\addplot[mark=square*,dashed,mark size=1.3pt,black]table [y=errh   , x=Nx]{IslandAtRest/mPDeC/CFL01/convergence.txt};
			\addlegendentry{mPDeC5-WENO5 CFL=0.1}				
			\addplot[mark=triangle*,dashdotted,mark size=1.3pt,black]table [y=errh   , x=Nx]{IslandAtRest/DeC/CFL01/convergence.txt};
			\addlegendentry{DeC5-WENO5 CFL=0.1}				
			\addplot[mark=o,mark size=1.3pt,darkspringgreen]             table [y=errh   , x=Nx]{IslandAtRest/SSPRK/CFL02/convergence.txt};
			\addlegendentry{SSPRK(6,4)-WENO5 CFL=0.2}
			\addplot[mark=square*,dashed,mark size=1.3pt,darkspringgreen]table [y=errh   , x=Nx]{IslandAtRest/mPDeC/CFL02/convergence.txt};
			\addlegendentry{mPDeC5-WENO5 CFL=0.2}
			\addplot[mark=triangle*,dashdotted,mark size=1.3pt,darkspringgreen]table [y=errh   , x=Nx]{IslandAtRest/DeC/CFL02/convergence.txt};
			\addlegendentry{DeC5-WENO5 CFL=0.2}
			\addplot[mark=o,mark size=1.3pt,red]             table [y=errh   , x=Nx]{IslandAtRest/SSPRK/CFL1/convergence.txt};
			\addlegendentry{SSPRK(6,4)-WENO5 CFL=1}
			\addplot[mark=square*,dashed,mark size=1.3pt,red]table [y=errh   , x=Nx]{IslandAtRest/mPDeC/CFL1/convergence.txt};
			\addlegendentry{mPDeC5-WENO5 CFL=1}
			\addplot[mark=triangle*,dashdotted,mark size=1.3pt,red]table [y=errh   , x=Nx]{IslandAtRest/DeC/CFL1/convergence.txt};
			\addlegendentry{DeC5-WENO5 CFL=1}
			\addplot[mark=square*,mark size=1.3pt,blue]             table [y=errh   , x=Nx]{IslandAtRest/mPDeC/CFL15/convergence.txt};
			\addlegendentry{mPDeC5-WENO5 CFL=1.5}
			\addplot[magenta,domain=50:200, dotted]{200*x^(-2)};
			\addlegendentry{second order}
		\end{axis}
	\end{tikzpicture}
	\caption{Wet and dry lake at rest test: error and computational time.}\label{fig:WDLAR}
\end{figure}
First, we test the non-well-balanced schemes, to asses the capability of preserving the water height positivity and the accuracy of such methods. The positivity of the classical schemes is not preserved, even with the positivity limiter of remark~\ref{rem:positivity_limiter}. In the dry region, the SW model \eqref{eq:SWE} does not hold and for the time integration schemes it is hard to verify hypotheses that guarantee the positivity of the solution. Hence, for the classical schemes, we force the positivity of the water height every time we need to compute the flux or to convert the variables from conservative to primitive and \textit{vice versa}. On the other hand, the mPDeC scheme always preserve the positivity of the solution and none of these tricks is required. \\
In figure~\ref{fig:WDLAR} we observe that for similar computational times, the mPDeC5-WENO5 gives the same accuracy of all classical schemes. For all schemes the accuracy is 2, as the solution is not $\mathcal{C}^1$ everywhere. The difference is that mPDeC5-WENO5 can be run up to CFL=1, without any problem, while for CFL=1 the SSPRK(6,4)-WENO5 scheme, even with the checks on the positivity, can have problems and might have exploding velocities and water heights. Since the reconstruction does not guarantee the positivity for such high CFLs, it is not safe to run such simulations.

Adding the well-balanced technique we obtain machine precision errors for all schemes.

\subsection{Almost dry lake at rest}
Now we modify the previous test, in order to have a smooth solution and to be able to obtain a fifth order accuracy in the schemes. We consider again the bathymetry \eqref{eq:bathymetry_pert} 
on the domain $[-5,5]\times[-2,2]$.
The lake at rest solution is defined, this time, as
\begin{equation}
	h(x,y) = \max\lbrace 0.999 - b(x,y),0 \rbrace ,\qquad u=v=0.
\end{equation}
\begin{table}
\centering
\caption{Almost dry lake at rest: mPDeC5-WENO5\label{tab:LAR_wet_mPDeC}}\vspace*{-3mm}
\begin{tabular}{|c||c|c||c|c||c|c|}	\hline
	  $N_x$   & Error $h$  &  Order $h$ & Error $u$  &  Order $u$ & Error $v$  &  Order $v$ \\ \hline
  600  &   8.346e-05  &  3.899  &  7.759e-04 & 3.057  &  6.911e-04 & 3.091 \\ 
  700  &   4.166e-05  &  4.508  &  4.526e-04 & 3.498  &  3.648e-04 & 4.144 \\ 
  800  &   2.134e-05  &  5.007  &  2.533e-04 & 4.346  &  1.886e-04 & 4.941 \\ 
 1000  &   6.185e-06  &  5.551  &  7.406e-05 & 5.511  &  5.225e-05 & 5.753 \\ 
 1200  &   2.219e-06  &  5.621  &  2.478e-05 & 6.006  &  1.774e-05 & 5.924 \\ 
 1400  &   1.120e-06  &  4.438  &  1.034e-05 & 5.669  &  7.561e-06 & 5.532 \\ 
 1600  &   5.896e-07  &  4.804  &  5.200e-06 & 5.148  &  3.799e-06 & 5.155 \\ 
 \hline
\end{tabular} \vspace*{3mm}

\caption{Almost dry lake at rest: SSPRK(6,4)-WENO5 \label{tab:LAR_wet_SSPRK}}\vspace*{-3mm}
\begin{tabular}{|c||c|c||c|c||c|c|}	\hline
	  $N_x$   & Error $h$  &  Order $h$ & Error $u$  &  Order $u$ & Error $v$  &  Order $v$ \\ \hline
  600  &   8.378e-05  &  3.891  &  7.765e-04 & 3.057  &  6.907e-04 & 3.091 \\ 
  700  &   4.191e-05  &  4.493  &  4.532e-04 & 3.493  &  3.644e-04 & 4.147 \\ 
  800  &   2.160e-05  &  4.963  &  2.538e-04 & 4.342  &  1.884e-04 & 4.942 \\ 
 1000  &   6.426e-06  &  5.434  &  7.432e-05 & 5.504  &  5.228e-05 & 5.745 \\ 
 1200  &   2.569e-06  &  5.029  &  2.502e-05 & 5.972  &  1.797e-05 & 5.857 \\ 
 1400  &   1.376e-06  &  4.051  &  1.058e-05 & 5.586  &  7.805e-06 & 5.410 \\ 
 1600  &   9.062e-07  &  3.127  &  5.496e-06 & 4.901  &  4.096e-06 & 4.829 \\ 
 \hline
\end{tabular}
\end{table}
Notice that the bathymetry has a peak with value 1, but, depending on the mesh refinement, the peak will be lower. In most of the simulations, this test will be completely wet and $\mathcal{C}^\infty$. If we discretize the mesh with more than $600 \times 180$ elements, the water level will be below the threshold $\varepsilon=10^{-6}$. As before, we initialize the water heigth at least equal to $\varepsilon=10^{-6}$ and we let the schemes evolve up to final time $T=1$.\\
We compare the mPDeC5-WENO5 and the SSPRK(6,4)-WENO5 schemes. The test is steady, hence, we expect only the spatial discretization to be the only responsible of the order of accuracy. For the SSPRK(6,4)-WENO5 to be run, we need to introduce some extra checks on the water height and computations of the flux, so that it does not become negative, while for the mPDeC5 time integration we do not need this type of extra corrections.
In both cases we expect to have a small perturbation of the solution while computing the $L_1$ error, indeed, the initial condition set with minimum level at $10^{-6}$ should introduce an error, in very few cells, of this order.\\
In table~\ref{tab:LAR_wet_mPDeC} there is the error analysis for the mPDeC5-WENO5 method, while in table~\ref{tab:LAR_wet_SSPRK} there is the one referred to the SSPRK(6,4)-WENO5 method. Despite expecting the error of the initialization to become evident around error of $10^{-6}$, for the mPDeC5 simulation, this small perturbation confined to very few cells does not propagate much and lead to very accurate results also for errors $\approx 5\cdot 10^{-7}$. Even the correction in remark~\ref{rem:division_zero} does not seem to affect the accuracy of the solution, probably because the water height never reaches values much lower than $10^{-6}$.	So, the order of accuracy stays very close to five, the expected one, even for almost dry solutions.\\
On the other side, in the SSPRK simulation the need of extra corrections in the flux every time the solution is below $10^{-6}$ adds further errors that are visible at level of $L_1$ error around $10^{-6}$. Indeed, it seems that the error of the water height starts plateauing close to that value. And there it loses the expected fifth order of accuracy.

\subsection{Perturbation of the lake at rest}
In order to better highlight the improvements one gets with the well-balanced implementation, a perturbation analysis is run
on a problem with both wet and dry areas. This test case is run on the rectangular domain $[-5,5]\times[-2,2]$.
The bathymetry is given by \eqref{eq:bathymetry_pert}.
The lake at rest solution that we are going to conserve with the well-balanced implementation is
\begin{equation}
h(x,y) = \max\lbrace 0.7 - b(x,y),\varepsilon \rbrace ,\qquad u=v=0,
\end{equation}
with $\varepsilon=10^{-6}$, whereas the perturbation shape that we want to study is 
\begin{equation}
\tilde{h}(x,y) = h(x,y) + \begin{cases} 0.05\,e^{1-\frac{1}{(1-\rho^2)^2}}, & \text{ if } \rho^2<1,  \\ 0,  & \text{  else  },\end{cases} \qquad\text{ where } \rho^2=9((x+2)^2+(x-0.5)^2).
\end{equation}

\begin{figure}
\centering
\subfigure[Computational mesh]{\includegraphics[width=0.45\textwidth,trim={1cm 0cm 1cm 8cm},clip]{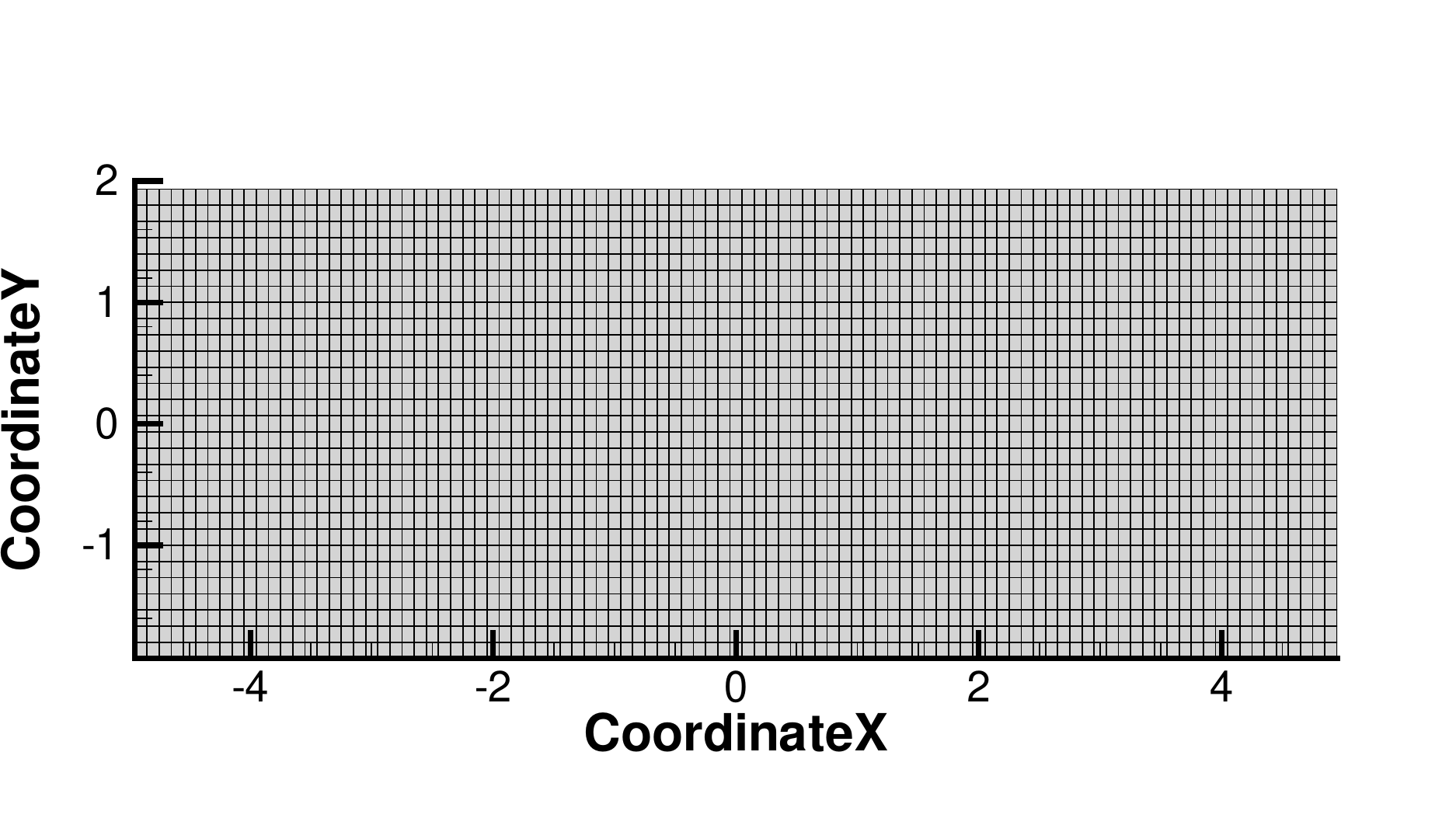}}
\subfigure[Bathymetry]{\includegraphics[width=0.45\textwidth]{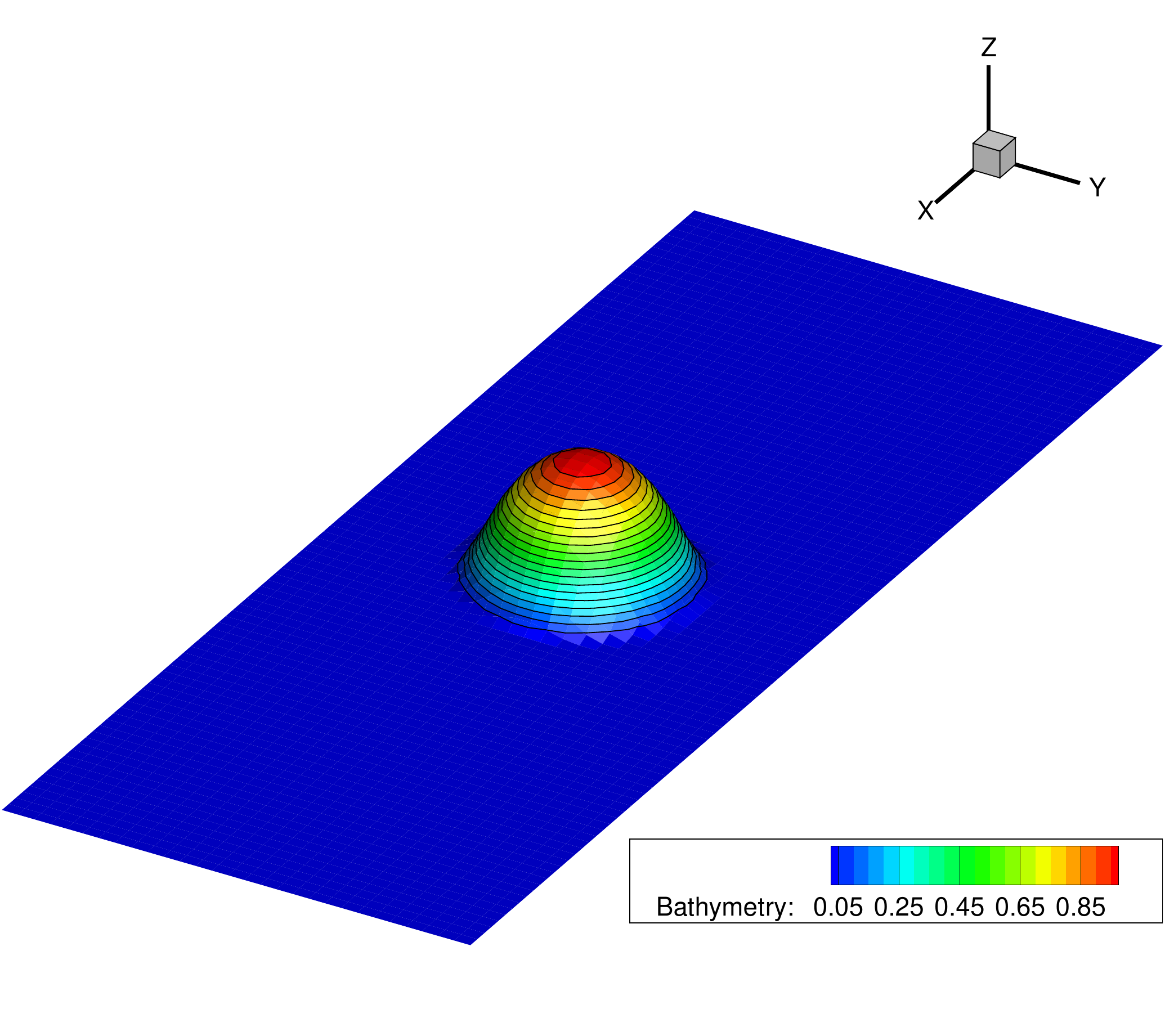}}
\caption{Perturbation analysis over a steady solution: test case setting}\label{fig:perturbation_set}
\end{figure}

Two simulations have been performed: one with the well-balanced correction and one without. Both simulations use WENO5 as spatial discretization and mPDeC5 for integrating the ODEs coming from the semi-discrete system.
From this test on the tolerances of the positivity limiter and of the mPDeC5 divisions is set to $\varepsilon=10^{-6}$.
The computational mesh and bathymetry plot are shown in Figure~\ref{fig:perturbation_set}.

\begin{figure}
\centering
\subfigure[WB, $t=0$]{\includegraphics[width=0.45\textwidth,trim={1cm 0cm 1cm 10cm},clip]{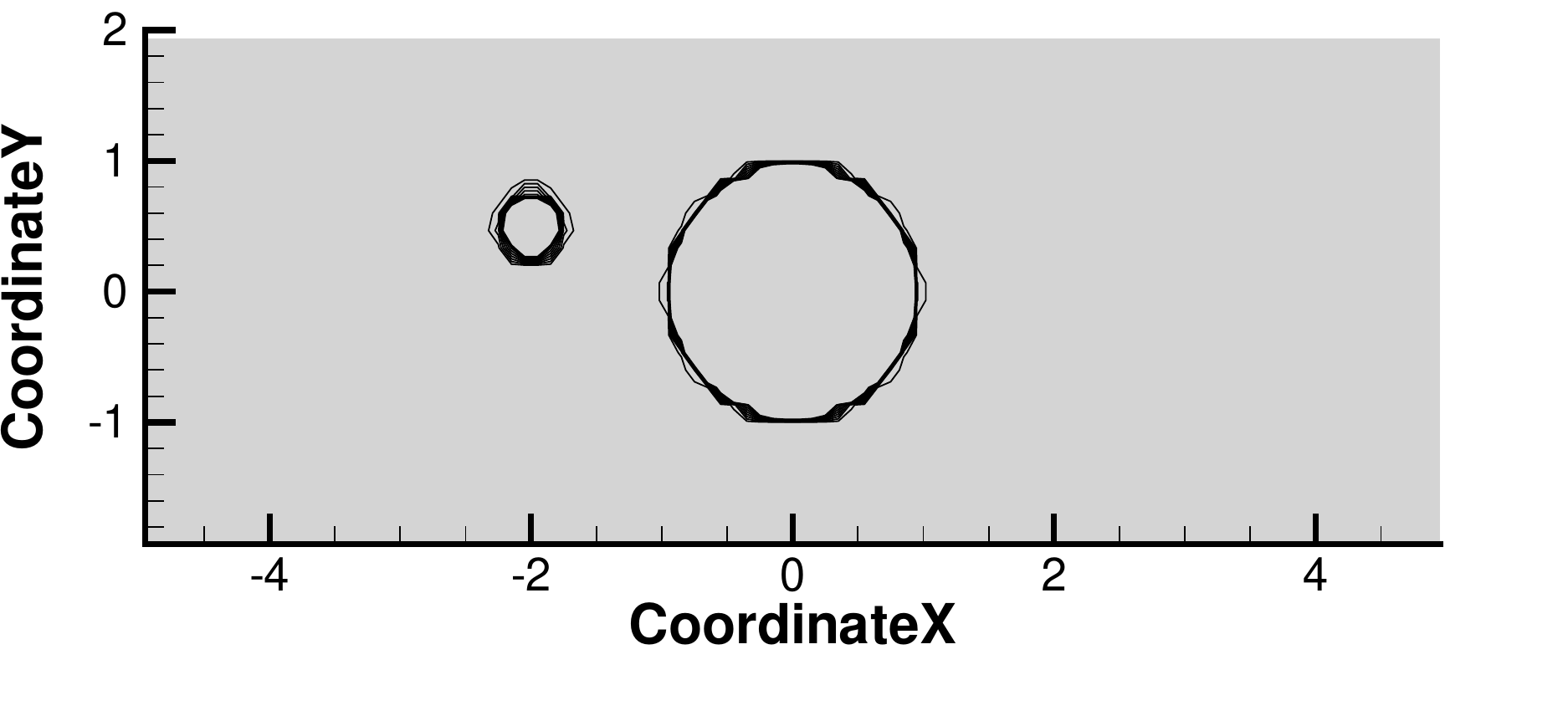}}
\subfigure[non WB, $t=0$]{\includegraphics[width=0.45\textwidth,trim={1cm 0cm 1cm 10cm},clip]{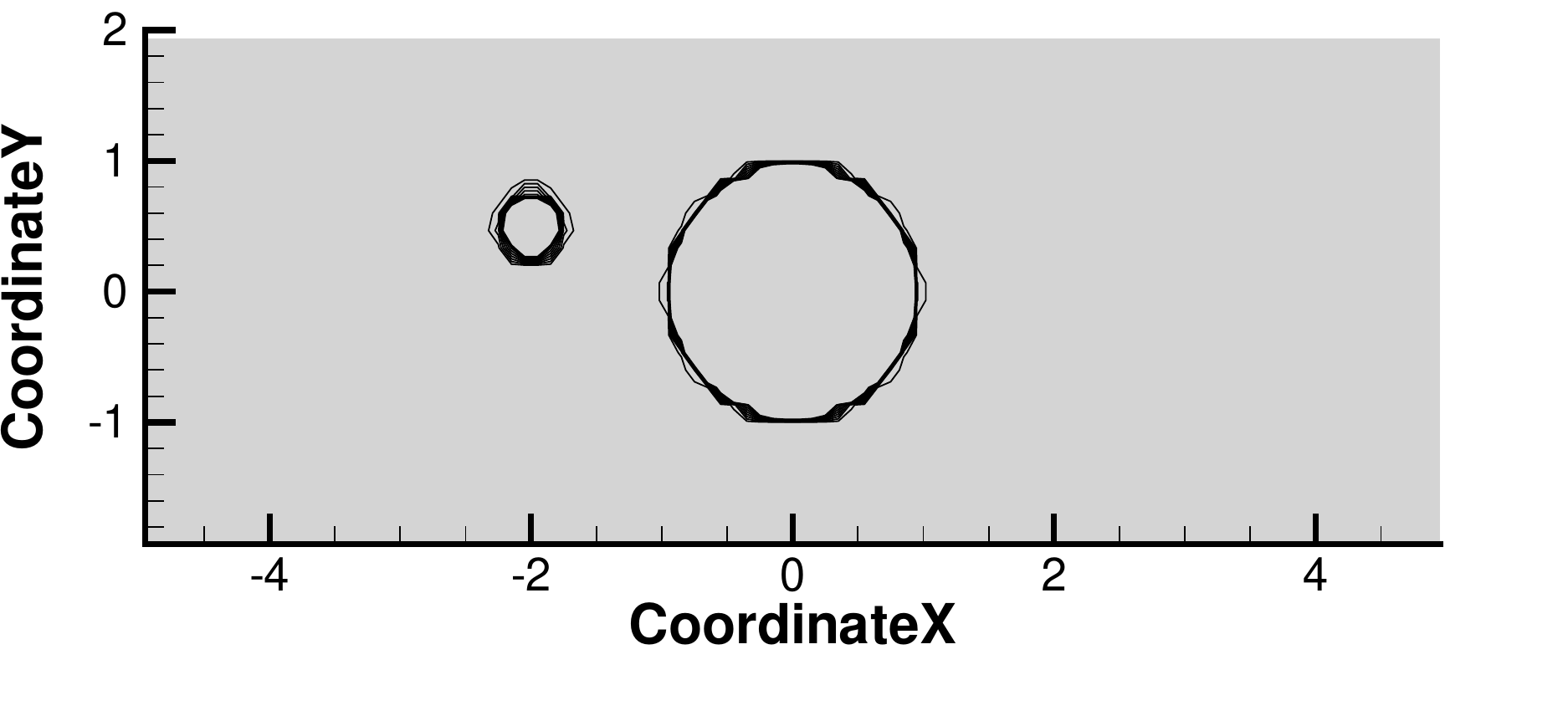}}
\qquad
\subfigure[WB, $t=0.25$]{\includegraphics[width=0.45\textwidth,trim={1cm 0cm 1cm 10cm},clip]{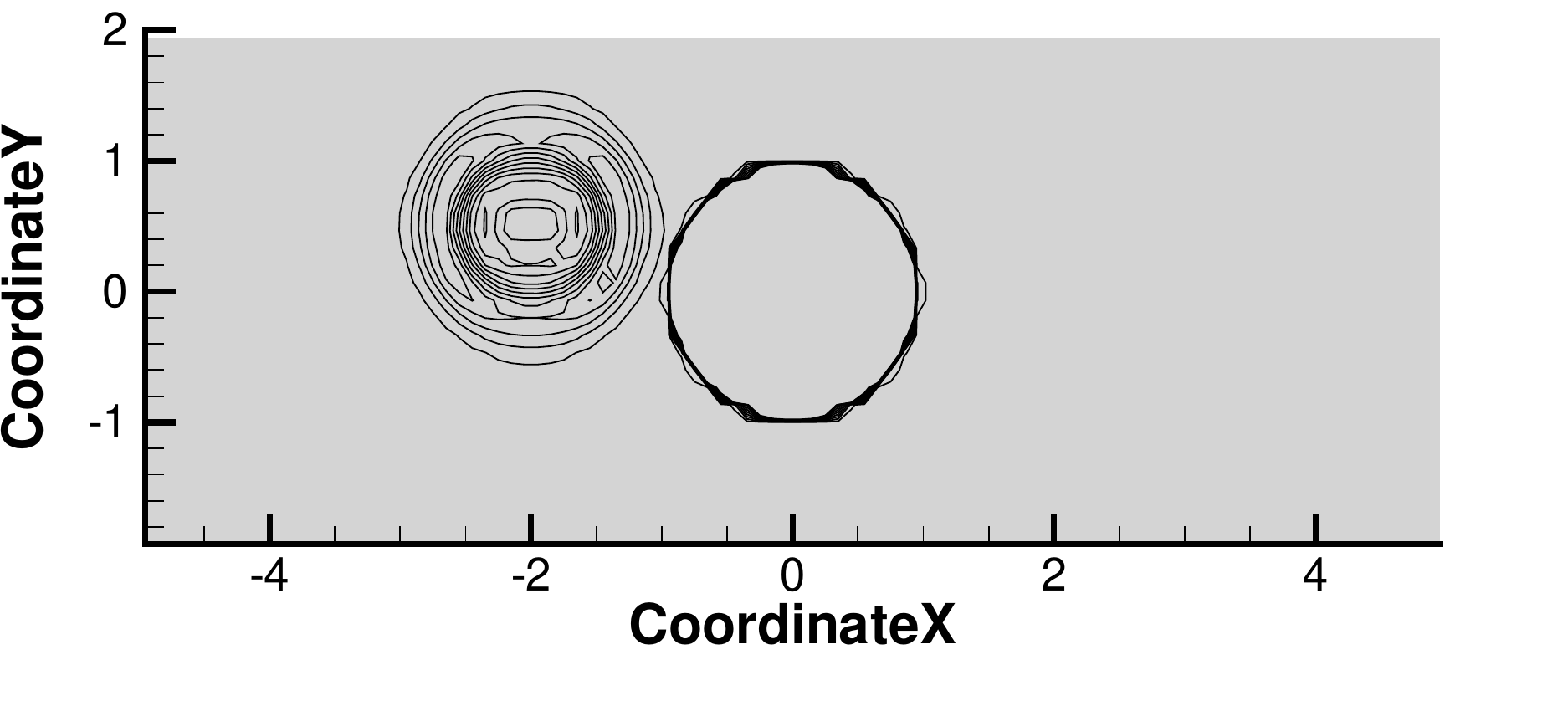}}
\subfigure[non WB, $t=0.25$]{\includegraphics[width=0.45\textwidth,trim={1cm 0cm 1cm 10cm},clip]{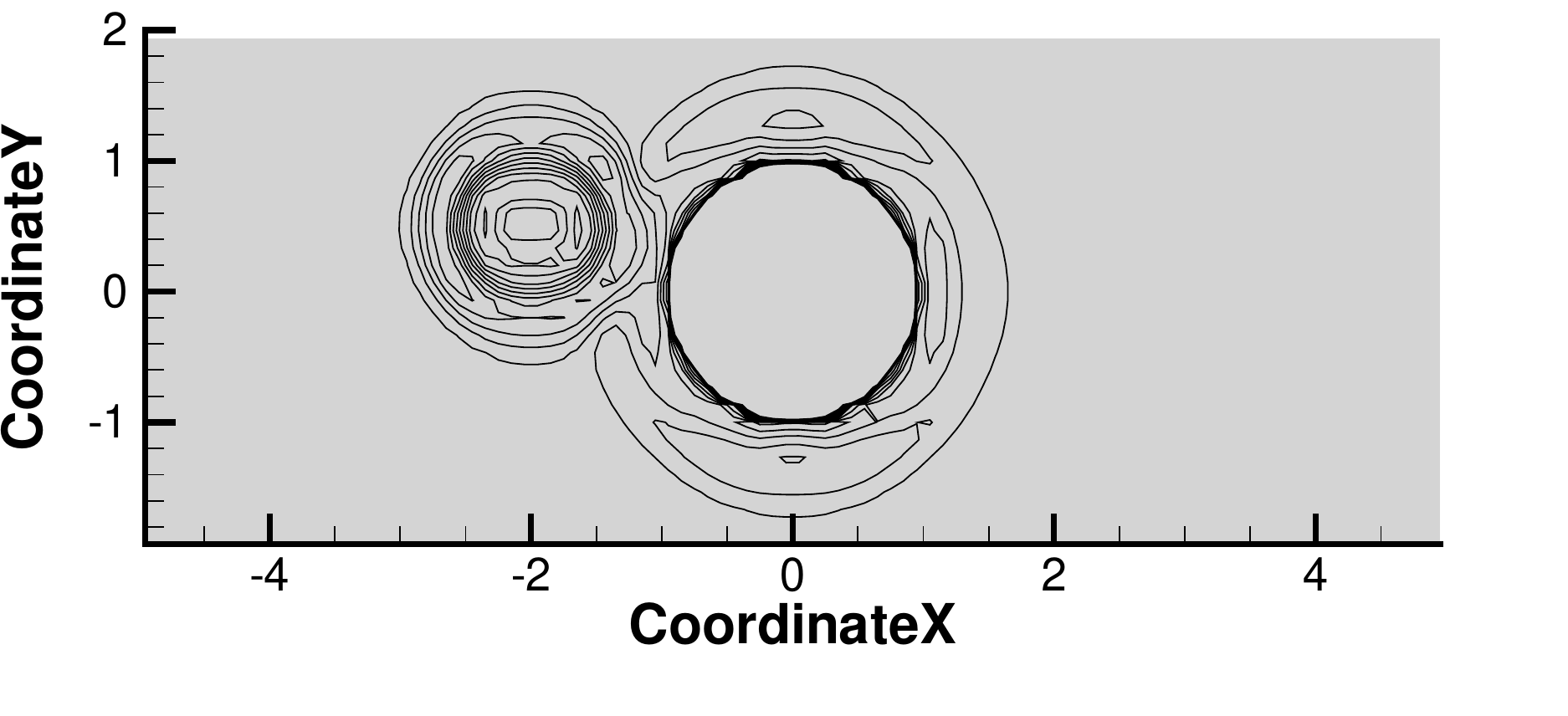}}
\qquad
\subfigure[WB, $t=0.5$]{\includegraphics[width=0.45\textwidth,trim={1cm 0cm 1cm 10cm},clip]{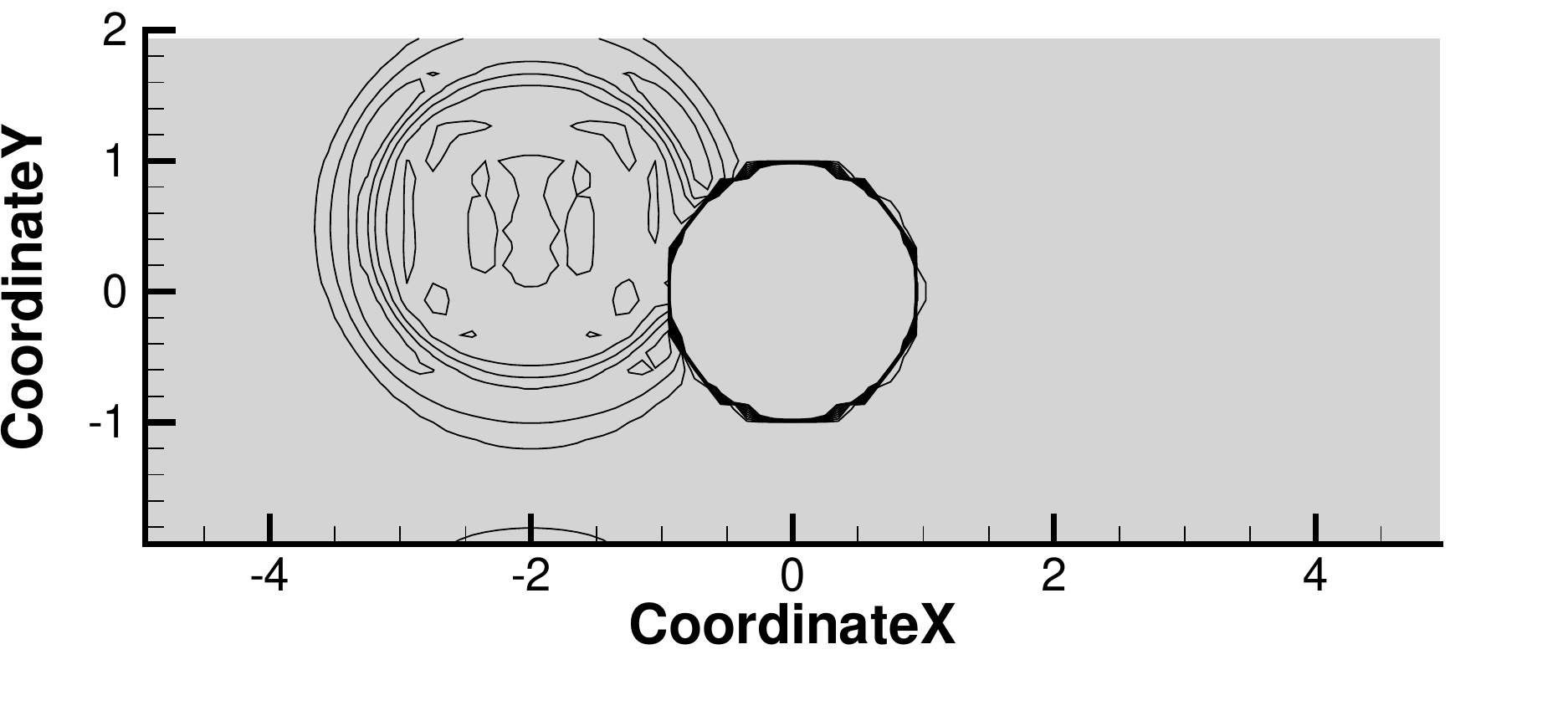}}
\subfigure[non WB, $t=0.5$]{\includegraphics[width=0.45\textwidth,trim={1cm 0cm 1cm 10cm},clip]{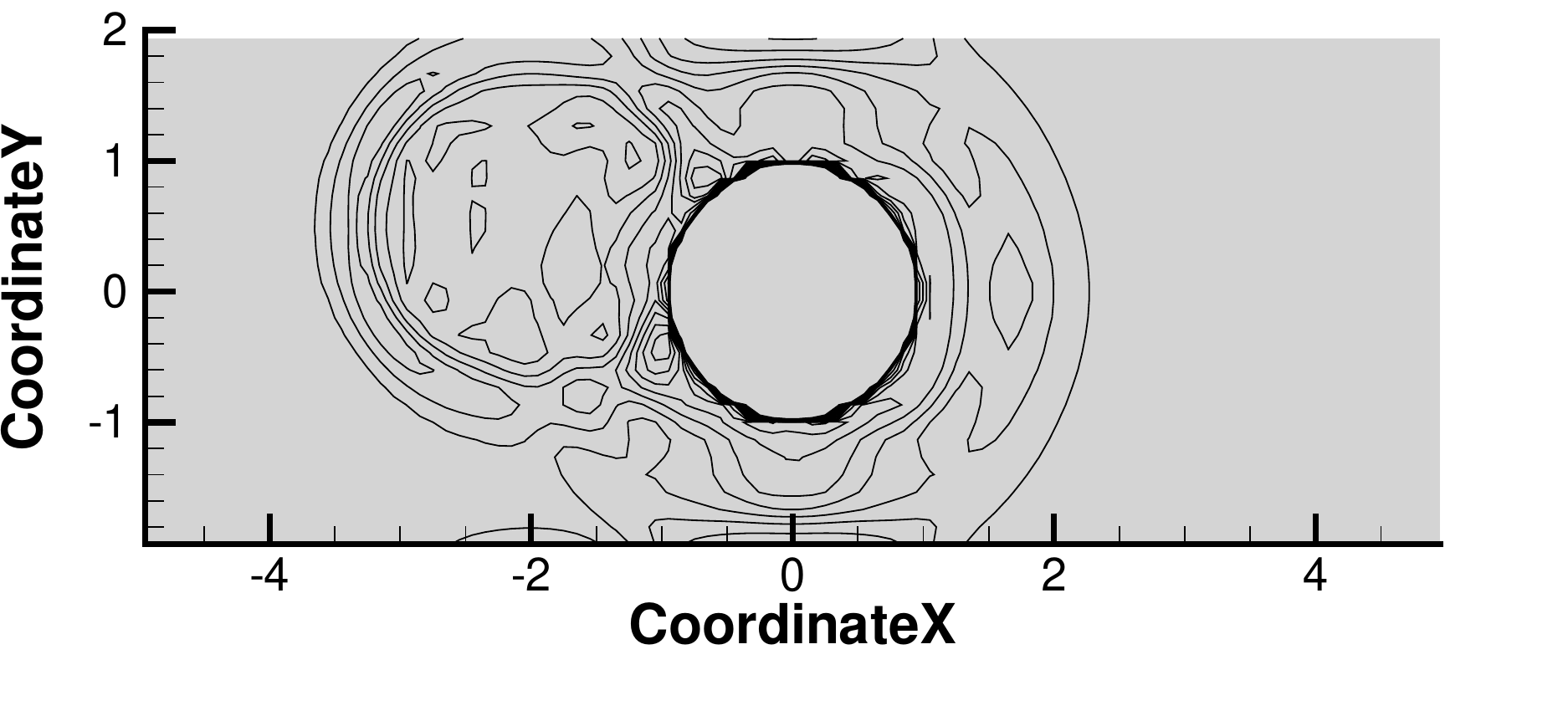}}
\qquad
\subfigure[WB, $t=1$]{\includegraphics[width=0.45\textwidth,trim={1cm 0cm 1cm 10cm},clip]{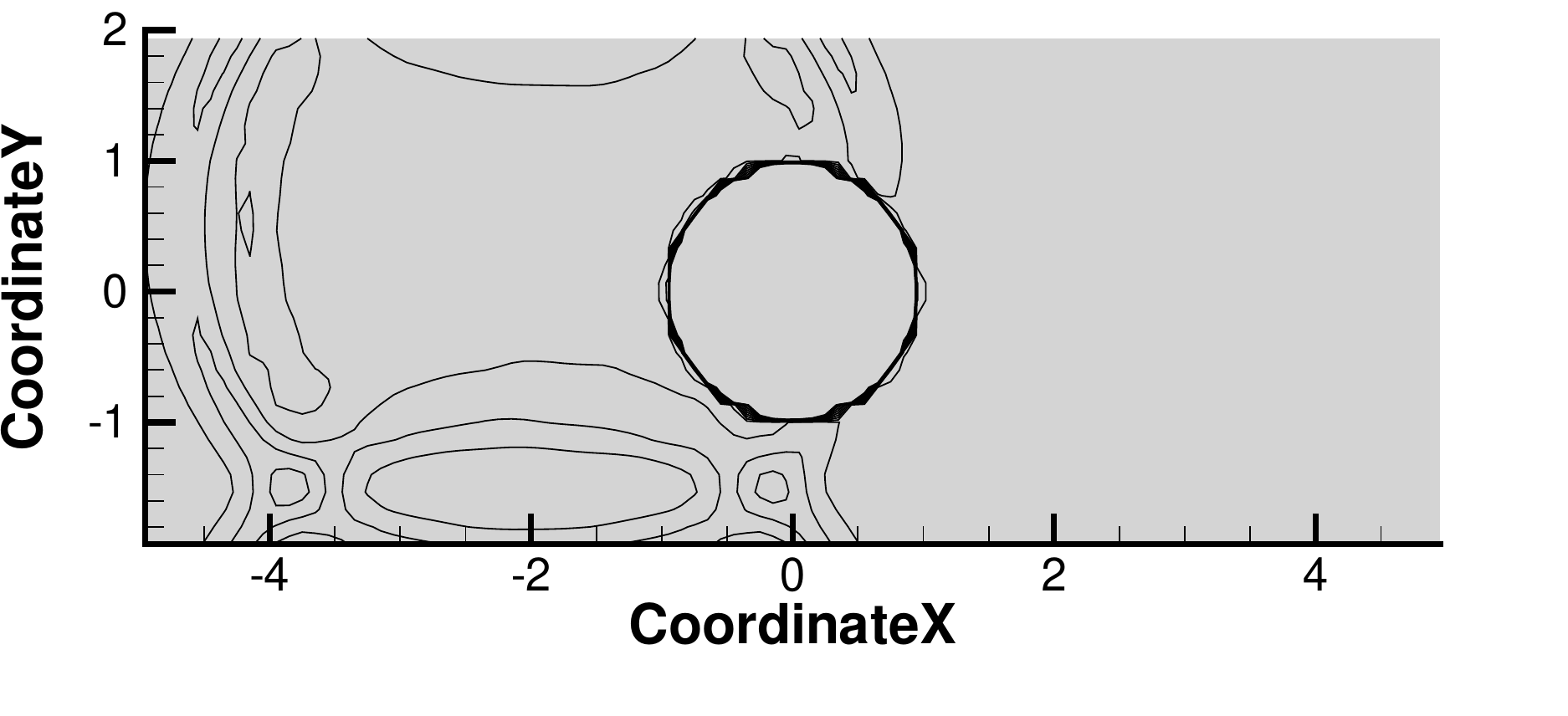}}
\subfigure[non WB, $t=1$]{\includegraphics[width=0.45\textwidth,trim={1cm 0cm 1cm 10cm},clip]{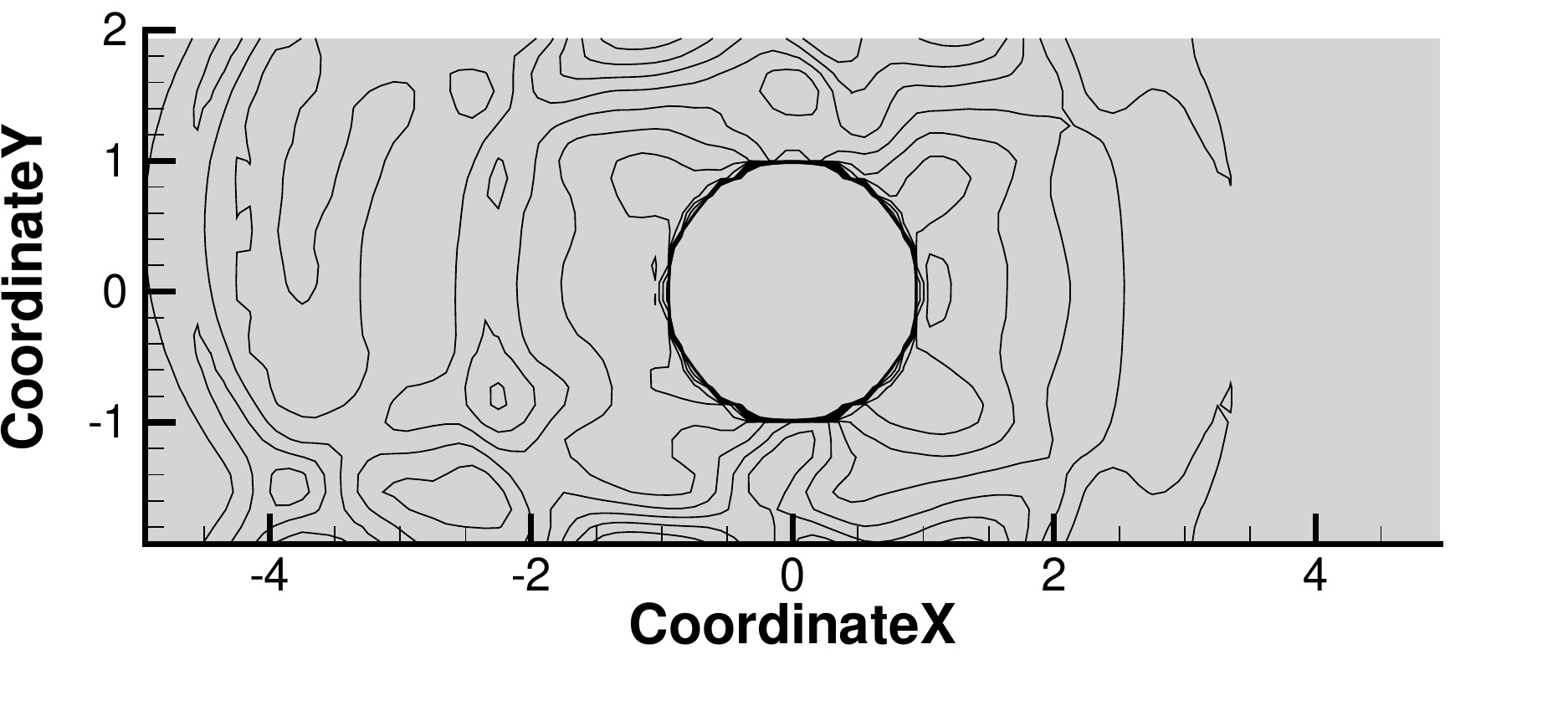}}
\caption{Perturbation analysis over a steady solution: water height $h$ isolines for well-balanced (left-hand side) and non well-balanced (right-hand side) results.}\label{fig:perturbation_results}
\end{figure}

The results obtained for the two implementations are displayed in Figure~\ref{fig:perturbation_results} where only the isolines of the water height $h$
are shown. Four snapshots are presented at different times of the simulation, $t=0,0.25,0.5,1$. The results on the right-hand side of Figure~\ref{fig:perturbation_results} 
are those computed without the well-balancing correction. As can be noticed, in this case, the numerical error propagates from the nonconstant bathymetry area around the island placed in the 
middle of the domain. This error propagates and interacts with the perturbation, making the perturbation waves indistinguishable from the noise. This test case allows to better assess the well-balanced implementation already tested in the 
previous test case. Indeed, for this case, we have a dry area which involves a jump in the derivative of the water height causing a reduction of the order of accuracy given 
by the WENO method, whose limiters work with the high order derivatives of the solution.
On the other side, the simulation runs with the well-balanced correction allows to exactly preserve the lake at rest solution over which the perturbation analysis is carried out.
This leads to a much better capturing of the perturbation, whose evolution is not influenced by the spurious disturbances coming from the wet-dry area. 

\subsection{Circular dry dam break problem}

\begin{figure}
\centering
\subfigure{\includegraphics[width=0.45\textwidth]{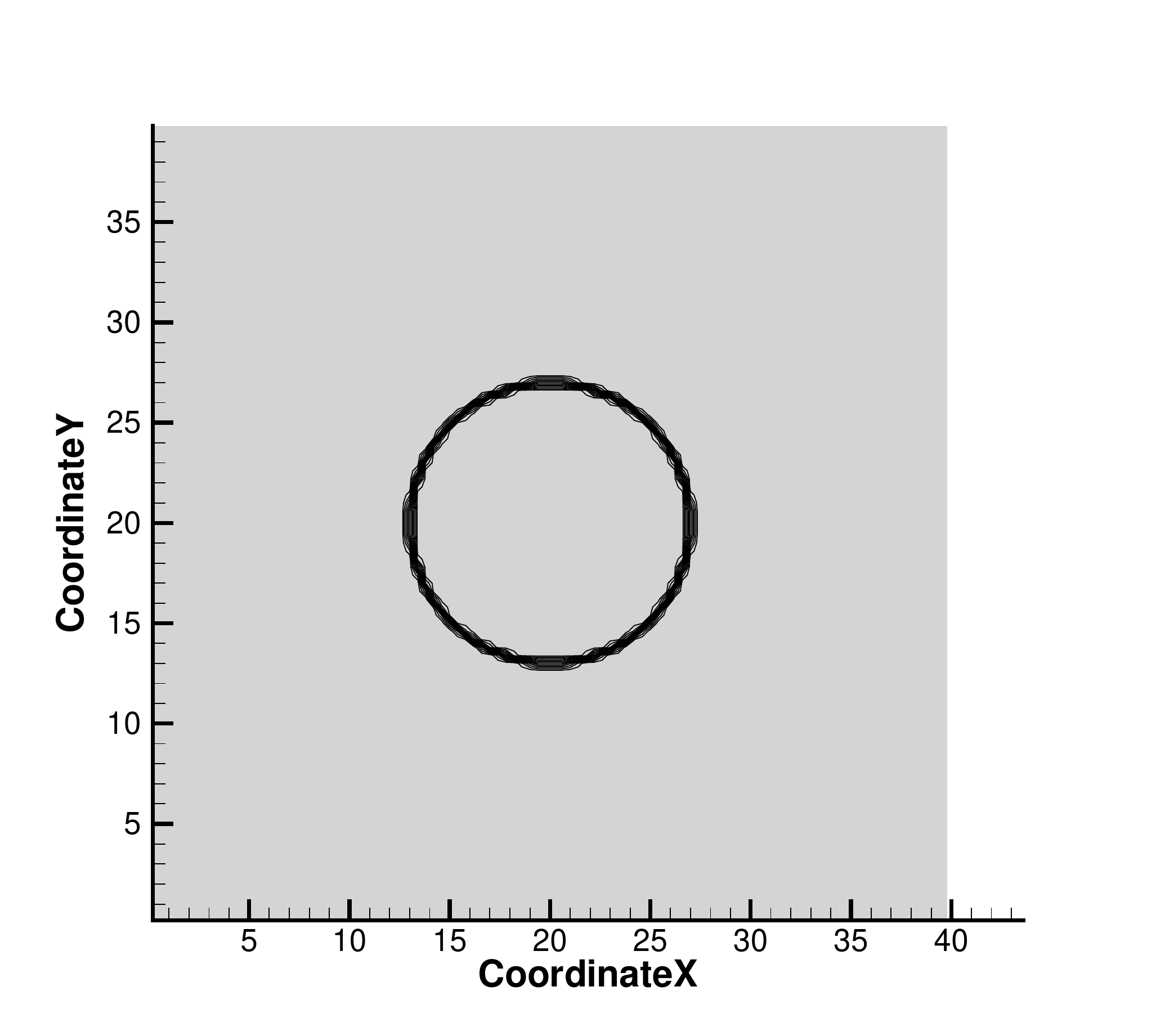}}
\subfigure{\includegraphics[width=0.45\textwidth]{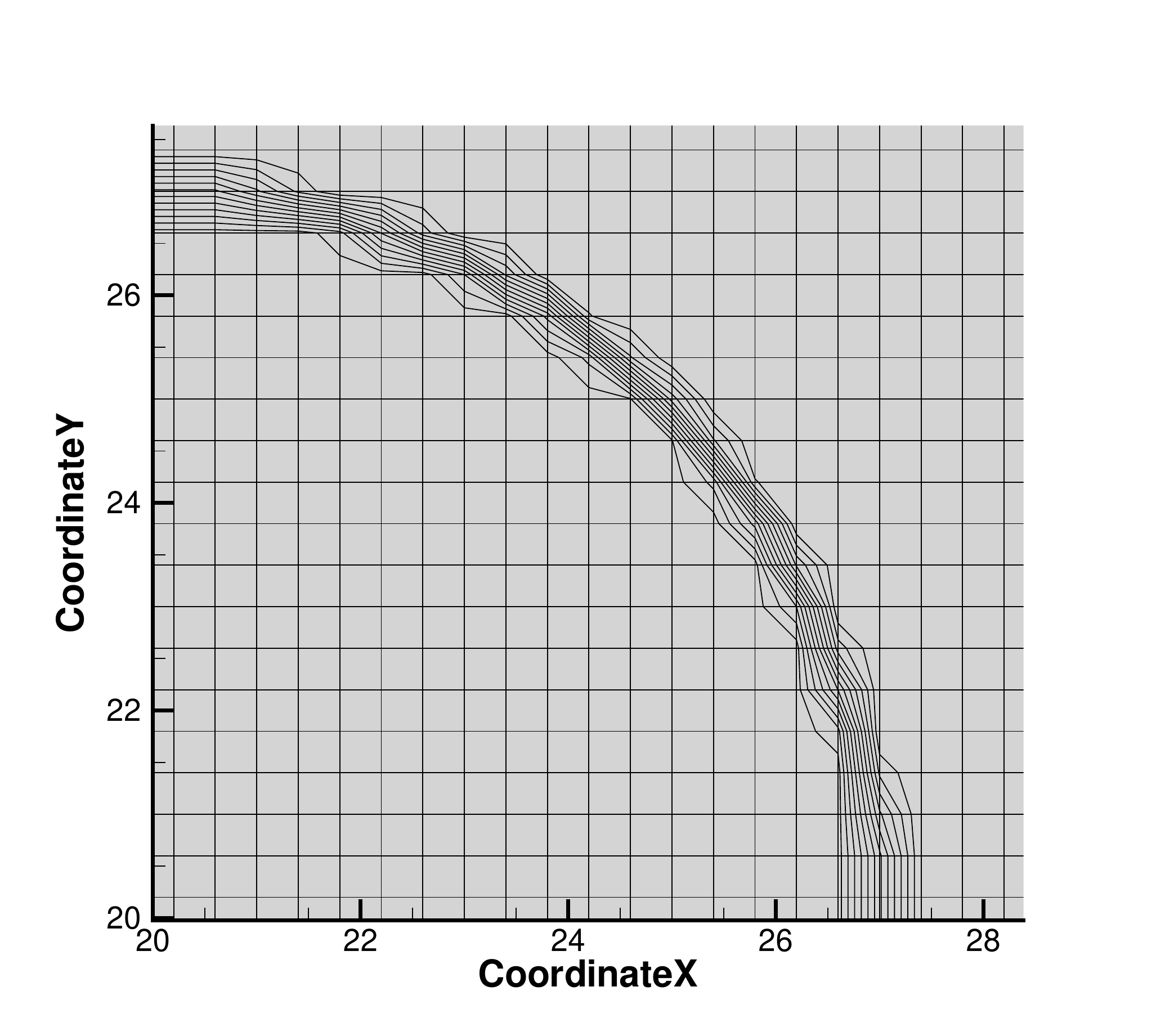}}
\caption{Circular dry dam break: computational domain and initial solution.}\label{fig:dam_set}
\end{figure}

\begin{figure}
\centering
\subfigure[$t=0$]{\includegraphics[width=0.45\textwidth]{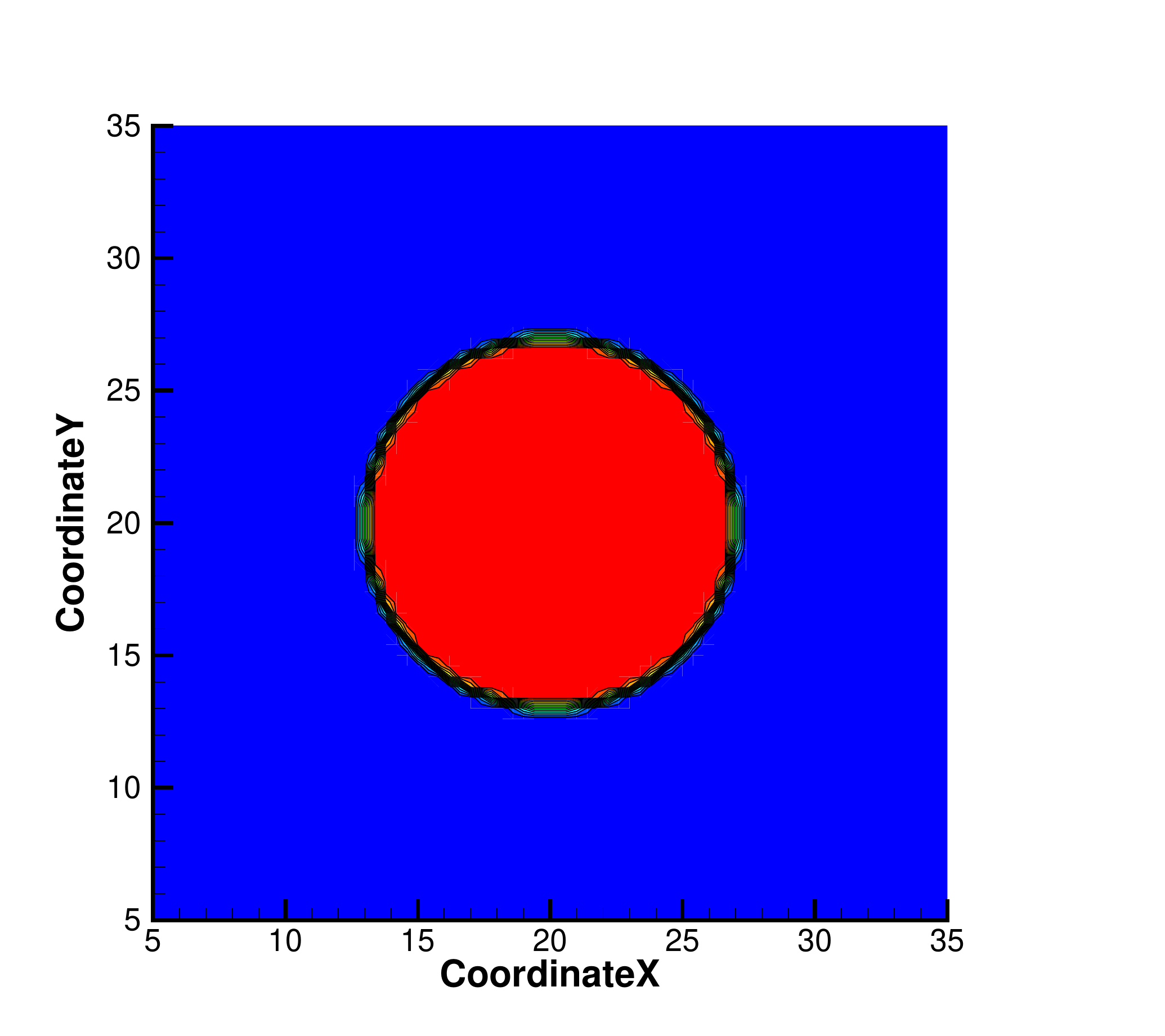}}
\subfigure[$t=0.3$]{\includegraphics[width=0.45\textwidth]{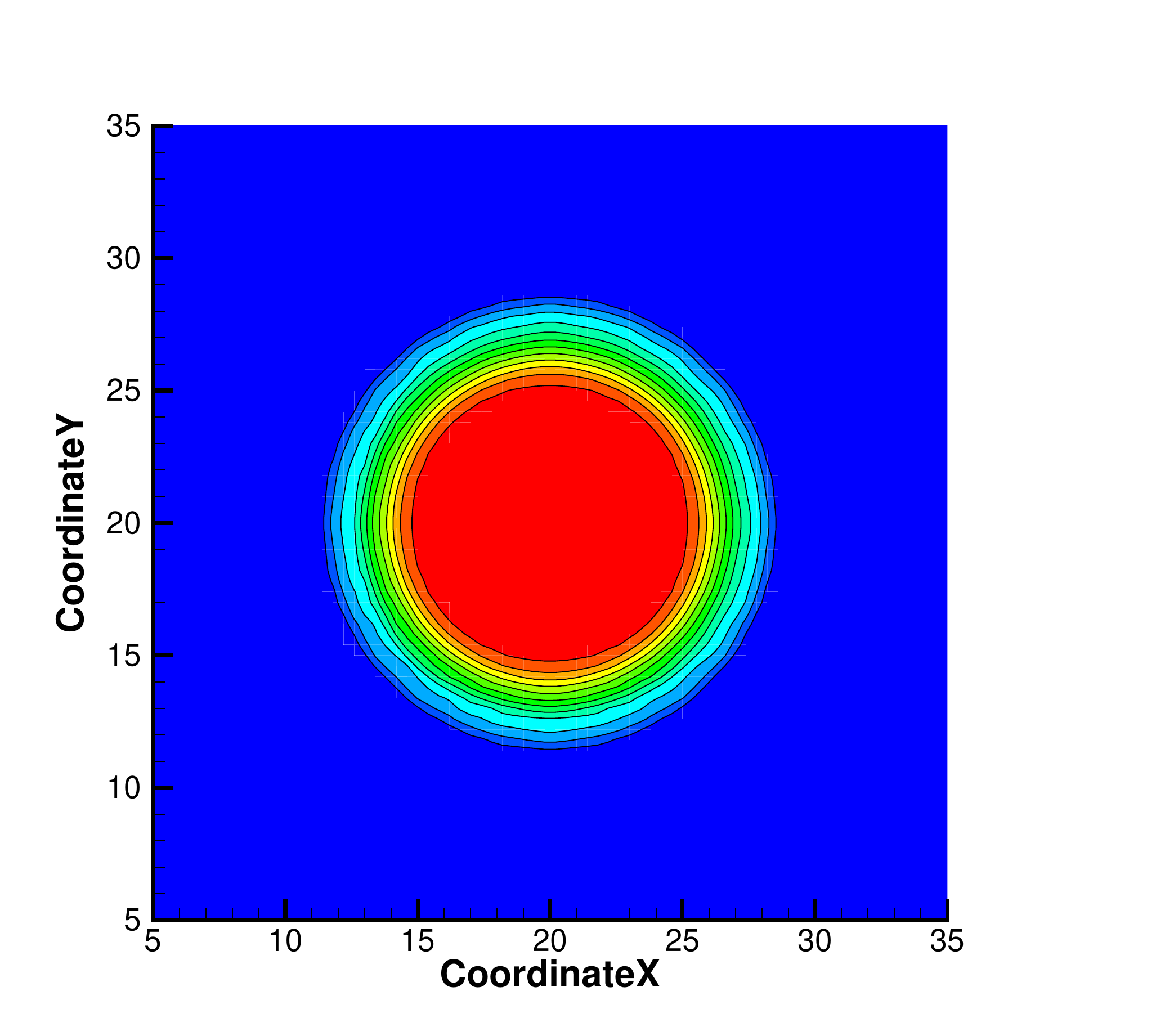}}
\qquad
\subfigure[$t=0.6$]{\includegraphics[width=0.45\textwidth]{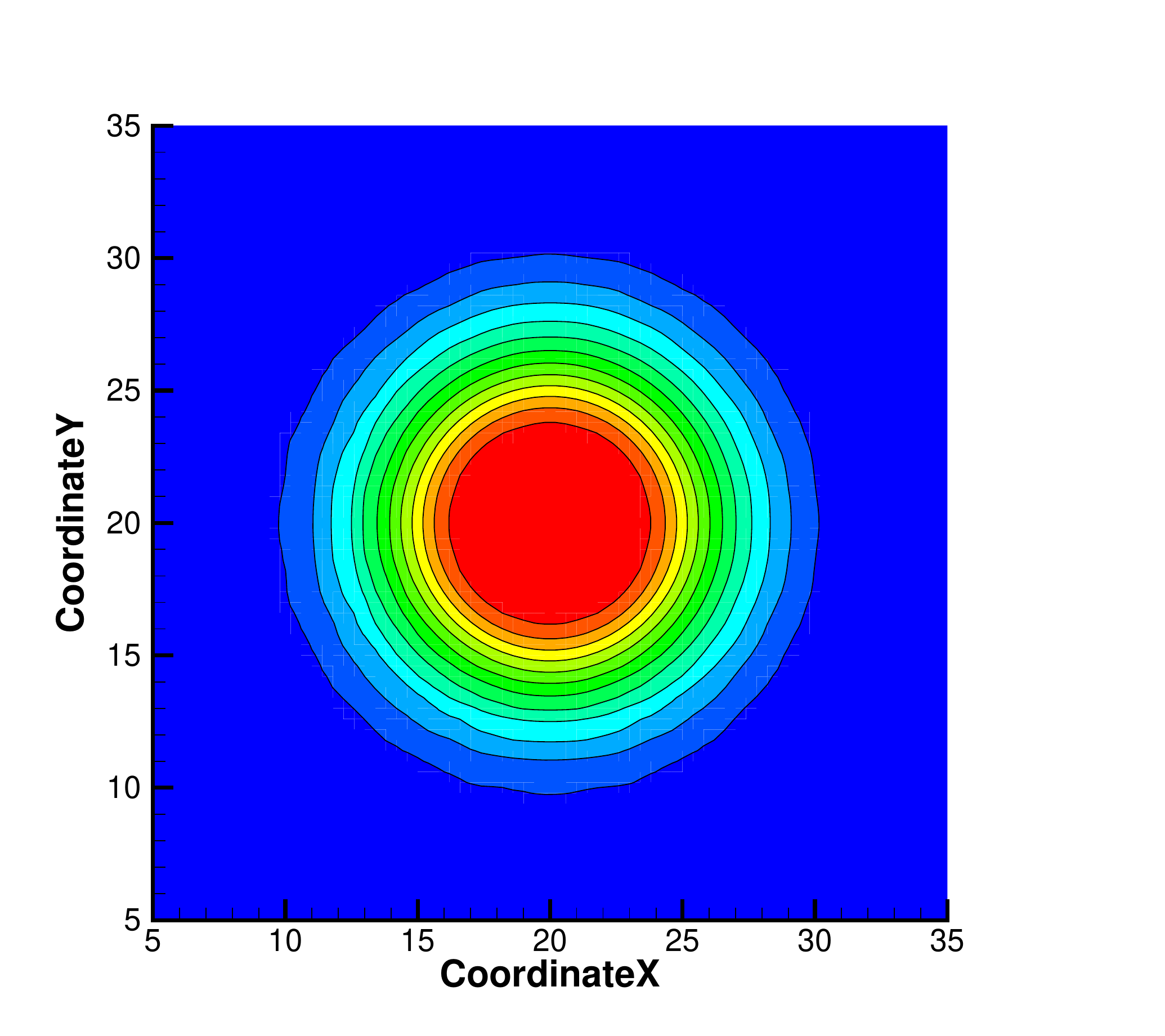}}
\subfigure[$t=0.9$]{\includegraphics[width=0.45\textwidth]{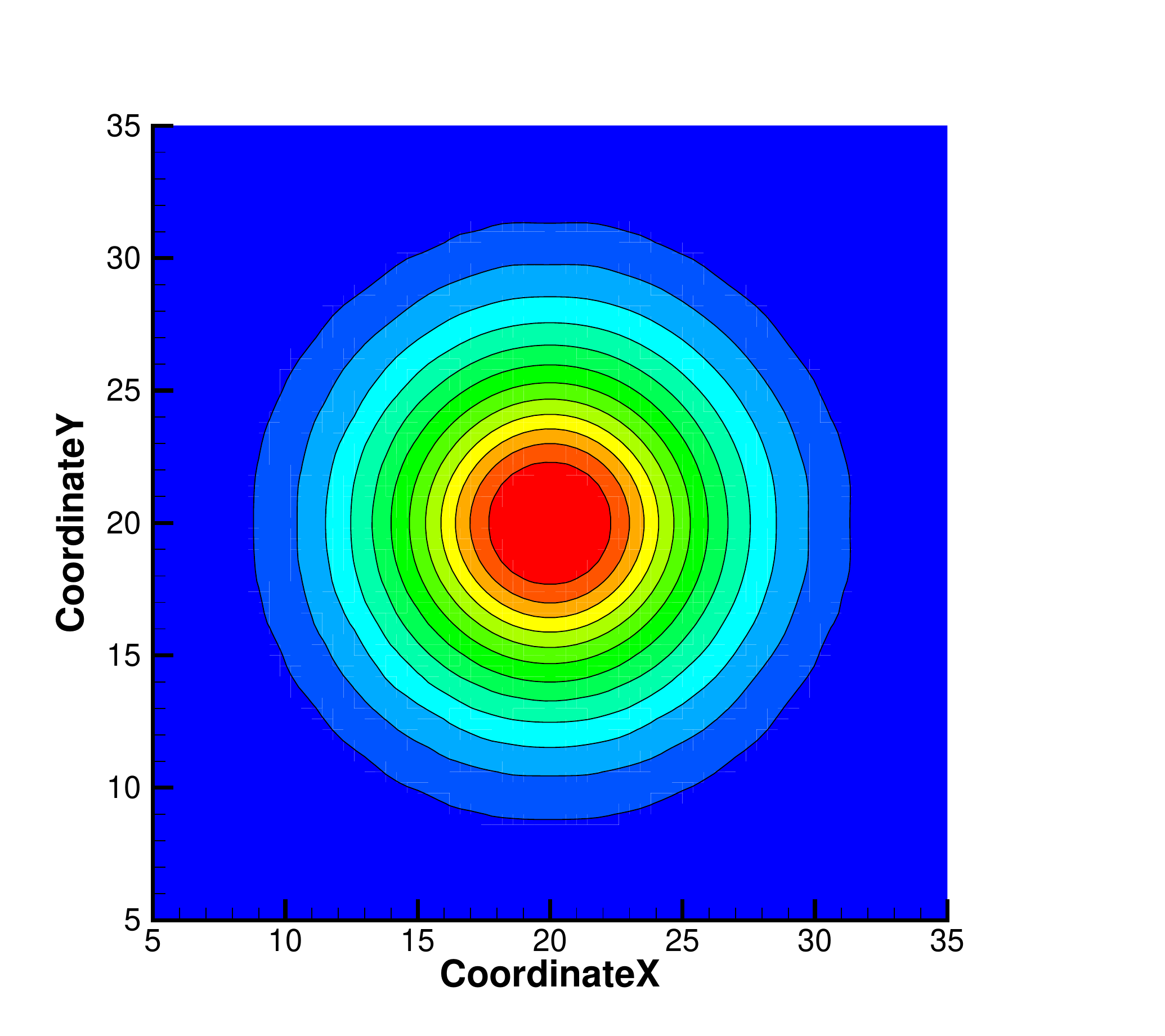}}
\qquad
\caption{Circular dry dam break: water height $h$ isocontours.}\label{fig:dam_sol}
\end{figure}

We simulate the break of a circular dam separating two basins with water heights $h_1=2.5$ and $h_2=\varepsilon=10^{-6}$, meaning that the water in the first basin is falling over a dry area which is all around it.
The radius of the discontinuity is $r=7$. A sketch of the initial condition, along with the computational mesh, is given in Figure~\ref{fig:dam_set}.
The computational domain is the square $[0,40]\times[0,40]$ discretized with $100\times 100$ cells and the simulation is run until a final time $t_{end}=0.9$.
This is a somewhat challenging test for the mPDeC5--WENO5 method that has to face both the capture of a sharp discontinuity and the progressing wetting of a dry area while
always maintaining its appealing properties.  
The results are printed for different times, $t=0,0.3,0.6,0.9$, in Figure~\ref{fig:dam_sol} showing the evolution of the water height. 
The advantages of the mPDeC5 have been clearly proven by running this challenging test case with different CFL conditions. As a matter of fact, 
we managed to run this test case with a CFL up until 1.5 while ensuring positivity.
In particular, the water height equation does not cause any CFL restriction due to the nature of the method used to solve it, which is unconditionally positive.
However, since the momentum equations are solved by means of an explicit DeC scheme, the CFL must be bounded.
It should be noticed that, in order to retain positivity in the spatial reconstruction, a positivity limiter has been implemented.
For explicit SSPRK methods, this limiter has been proven to cause a huge restriction in the CFL condition, which now has to be less than 
$\frac{1}{12}\text{CFL}^{\scriptscriptstyle \text{SSPRK}}$. 
Nevertheless, since arbitrary high order DeC cannot be recast as a convex combination of explicit Euler method, there is no proof that 
the solution would stay positive also under that strict condition.
On the other side, in the modified Patankar DeC framework, the limiter is only imposed on the water height equation which is unconditionally positive for any CFL by definition.
This means that even when the limiter plays a role in the simulation, like in this case, the scheme stays positive with much higher CFL numbers.
Furthermore, it should be underlined the fact that the system is linear implicit, which only requires a simple linear solver, e.g.\ Jacobi method.
Even though implicit methods are much slower than their explicit counterpart, solving a linear problem in this case only yields to a $18\%$ increase in the computational time
with respect to the fully explicit DeC5, which is nothing compared to the CFL restriction imposed by the positivity limiter.

\subsection{Circular wet dam break problem}

\begin{figure}
\centering
\subfigure[Isocontours of $h$ at time $t_{end}=0.8$]{\includegraphics[width=0.32\textwidth]{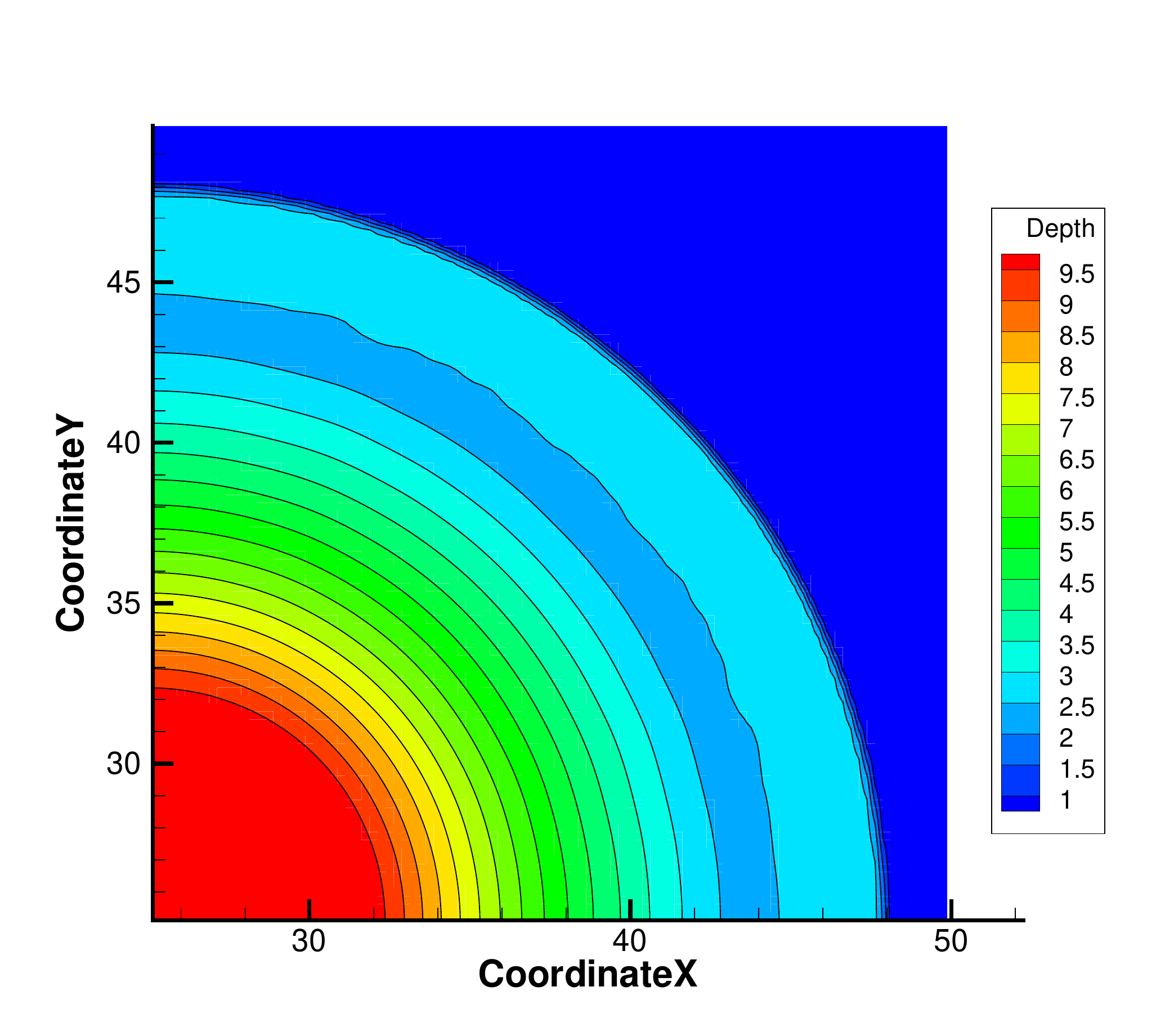}}
\subfigure[Plot of $h$ at time $t_{end}=0.8$]{\includegraphics[width=0.32\textwidth]{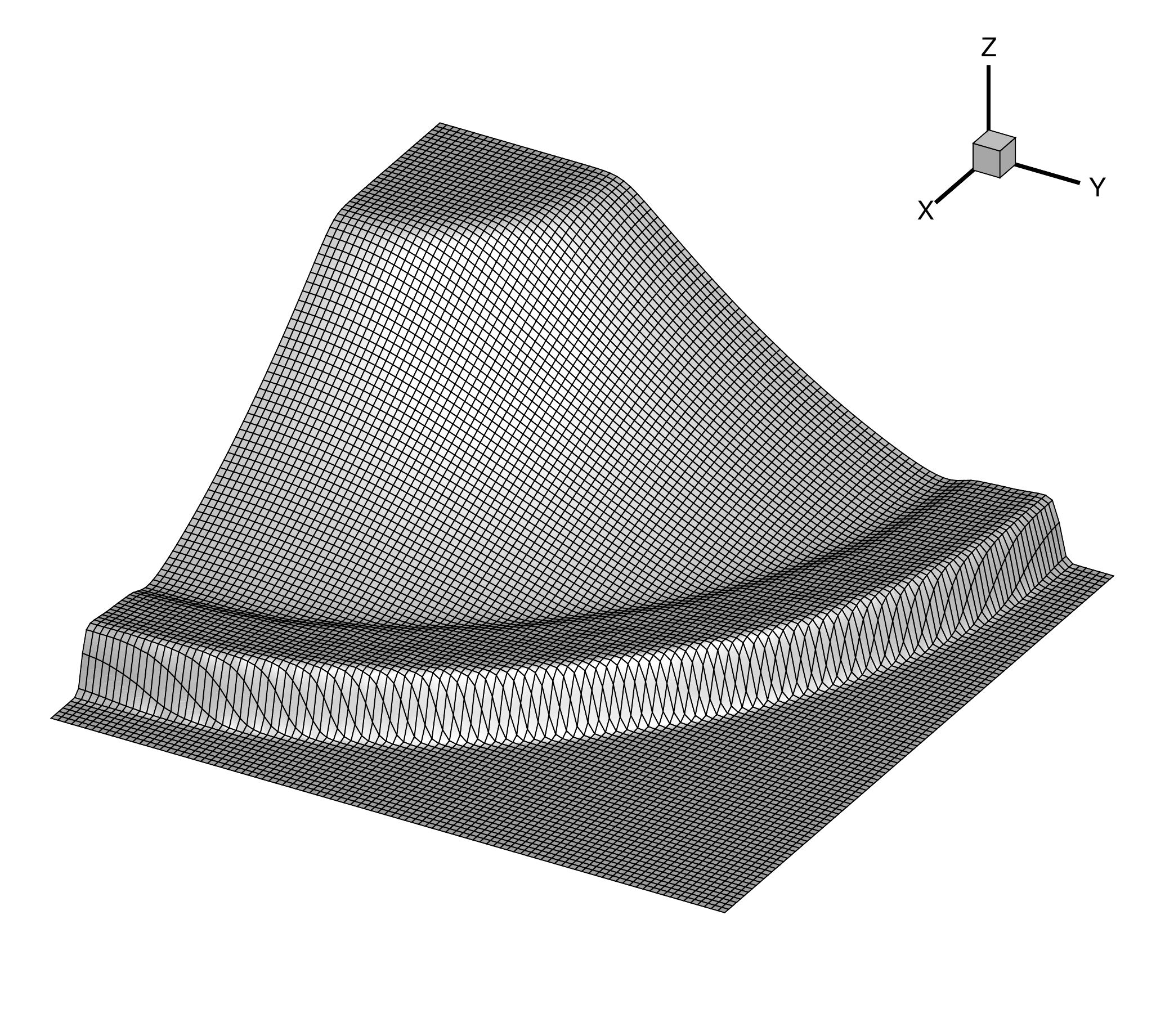}}
\subfigure[$h$ time evolution along the diagonal of the domain]{\includegraphics[width=0.32\textwidth]{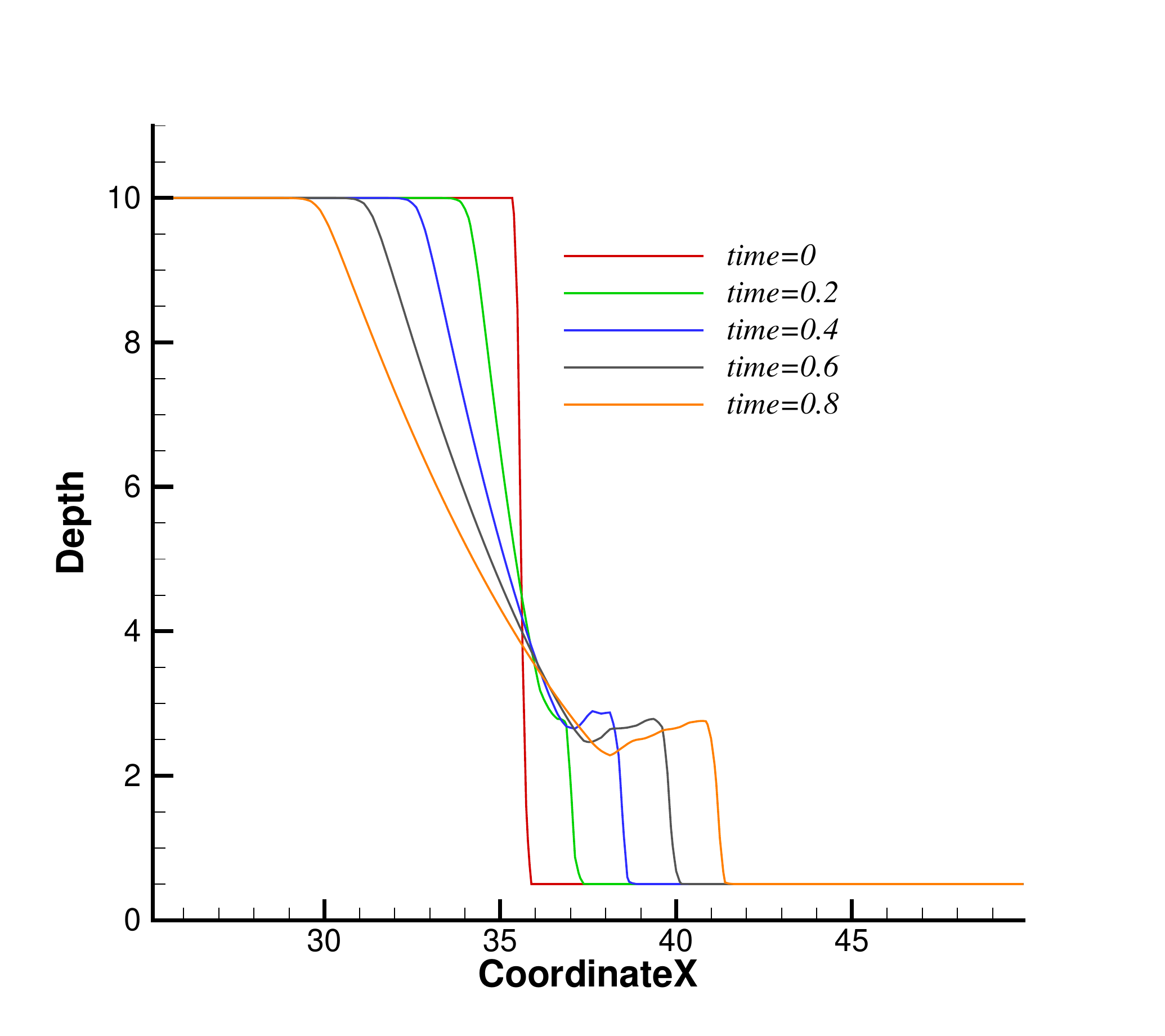}}
\caption{Circular wet dam break}\label{fig:damsod_sol}
\end{figure}

Next, we consider a wet dam break problem with a setup similar to the previous one. In this case the water height of the two basins are $h_1=10$ and $h_2=0.5$, 
like it was done in~\cite{ricchiuto2009stabilized}. In this case, the computational domain is the square $[0,50]\times[0,50]$ discretized by a $200\times200$ mesh.
The simulation is run until $t_{end}=0.8$ with CFL=1 and the solutions is displayed only for one quarter of the domain.
Figure~\ref{fig:damsod_sol} displays the final snapshot of the solution at time $t_{end}=0.8$ plotted both in 2D and 3D displaying the correct evolution of the water height and
the time evolution of the water depth extracted along the diagonal of the portion of the domain studied. The method shows good properties such as discontinuities sharply captured, no oscillations 
and high order approximation of smooth features. 

\subsection{Wave over dry island}

For the last test case, we simulate a wave crashing over a dry island showing the robustness of our method
when facing more realistic simulations. The computational domain is the rectangle $[-5,5]\times[-2,2]$ discretized
by a $400\times120$ mesh. The simulation can be reproduced by taking Eq.~\eqref{eq:bathymetry_pert} for the bathymetry
$b(x,y)$ and the following initial conditions
\begin{equation}
h(x,y) = 0.7 - b(x,y) + \begin{cases} 0.5\,e^{1-\frac{1}{(1-\rho^2)^2}}, & \text{ if } \rho^2<1, \\ 0,  & \text{  else,  }\end{cases} \qquad\text{ where } \rho^2=(x+2)^2 ,\qquad (u,v)=(1,0),
\end{equation}
in case $h(x,y) \leq \varepsilon=10^{-6}$, we set $h(x,y)=\varepsilon$ and $u=0$.
The simulation has been run with the WENO5$-$mPDeC5 scheme until a final time $t=1$ with CFL=0.9 and the results are displayed in Figure~\ref{fig:wodi_sol} for different times $t=0,0.25,0.5,0.75,1$.

In this case, the variable $\eta=h+b$ has been chosen to be plotted since it better represents the underlying physics.
\begin{figure}
	\centering
	\begin{tikzpicture}
		\begin{axis}[
			xmode=log,
			grid=major,
			xlabel={Mesh elements in $x$},
			ylabel={Jacobi iterations},
			xlabel shift = 1 pt,
			ylabel shift = 1 pt,
			legend pos= outer north east,
			legend style={nodes={scale=0.6, transform shape}},
			width=.45\textwidth
			]		
			\addplot[mark=circle,solid, mark size=1.3pt,black]             table [y=jacAve   , x=N]{IslandShort/CFL15_comparison.txt};
			\addlegendentry{CFL=1.5}
			\addplot[mark=star,dashed, mark size=1.3pt,blue]             table [y=jacAve   , x=N]{IslandShort/CFL08_comparison.txt};
			\addlegendentry{CFL=0.8}
			\addplot[mark=diamond*,dotted,mark size=1.3pt,red]table [y=jacAve   , x=N]{IslandShort/CFL04_comparison.txt};
			\addlegendentry{CFL=0.4}			
			\addplot[mark=square*,dotted,mark size=1.3pt,darkspringgreen]table [y=jacAve   , x=N]{IslandShort/CFL02_comparison.txt};
			\addlegendentry{CFL=0.2}	
			\addplot[mark=triangle*,dashdotted,mark size=1.3pt,magenta]table [y=jacAve   , x=N]{IslandShort/CFL01_comparison.txt};
			\addlegendentry{CFL=0.1}
			\addplot[name path=us_top,solid,gray!50]             table [y expr=\thisrow{jacAve}+0.5*\thisrow{jacStd}   , x=N]{IslandShort/CFL15_comparison.txt};
			\addplot[name path=us_bot,solid, gray!50]             table [y expr=\thisrow{jacAve}-0.5*\thisrow{jacStd}   , x=N]{IslandShort/CFL15_comparison.txt};
			\addplot[gray!30,fill opacity=0.3] fill between[of=us_top and us_bot];
			\addplot[name path=us_top,dashed,blue!50]             table [y expr=\thisrow{jacAve}+0.5*\thisrow{jacStd}   , x=N]{IslandShort/CFL08_comparison.txt};
			\addplot[name path=us_bot,dashed, blue!50]          table [y expr=\thisrow{jacAve}-0.5*\thisrow{jacStd}   , x=N]{IslandShort/CFL08_comparison.txt};
			\addplot[blue!30,fill opacity=0.3] fill between[of=us_top and us_bot];
			\addplot[name path=us_top,dotted,red!50]             table [y expr=\thisrow{jacAve}+0.5*\thisrow{jacStd}   , x=N]{IslandShort/CFL04_comparison.txt};
			\addplot[name path=us_bot,dotted, red!50]             table [y expr=\thisrow{jacAve}-0.5*\thisrow{jacStd}   , x=N]{IslandShort/CFL04_comparison.txt};
			\addplot[red!30,fill opacity=0.3] fill between[of=us_top and us_bot];		
			\addplot[name path=us_top,dotted,darkspringgreen!50]             table [y expr=\thisrow{jacAve}+0.5*\thisrow{jacStd}   , x=N]{IslandShort/CFL02_comparison.txt};
			\addplot[name path=us_bot,dotted, darkspringgreen!50]             table [y expr=\thisrow{jacAve}-0.5*\thisrow{jacStd}   , x=N]{IslandShort/CFL02_comparison.txt};
			\addplot[darkspringgreen!30,fill opacity=0.3] fill between[of=us_top and us_bot];		
			\addplot[name path=us_top,dashdotted,magenta!50]             table [y expr=\thisrow{jacAve}+0.5*\thisrow{jacStd}   , x=N]{IslandShort/CFL01_comparison.txt};
			\addplot[name path=us_bot,dashdotted, magenta!50]             table [y expr=\thisrow{jacAve}-0.5*\thisrow{jacStd}   , x=N]{IslandShort/CFL01_comparison.txt};
			\addplot[magenta!30,fill opacity=0.3] fill between[of=us_top and us_bot];
		\end{axis}
	\end{tikzpicture}
	\caption{Wave over dry island: Jacobi iterations and confidence interval. \label{fig:island_jacobi}}
\end{figure}
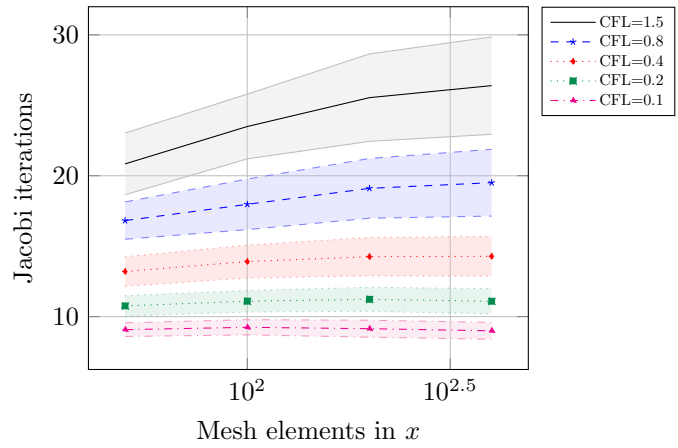
\begin{figure}
	\centering
	\subfigure[$t=0$]{\includegraphics[width=0.45\textwidth,trim={1cm 0cm 1cm 8cm},clip]{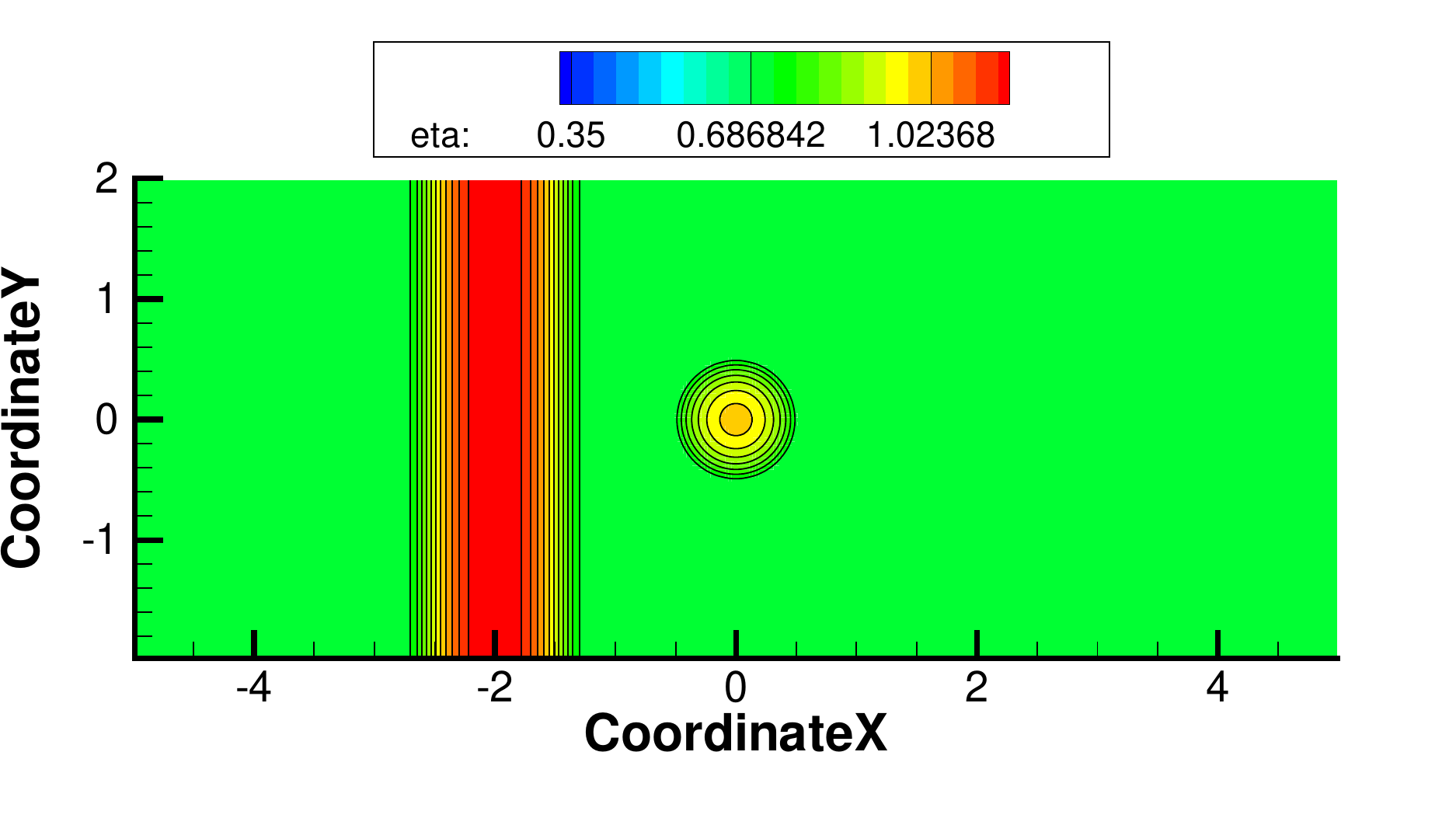}}
	\subfigure[$t=0.25$]{\includegraphics[width=0.45\textwidth,trim={1cm 0cm 1cm 8cm},clip]{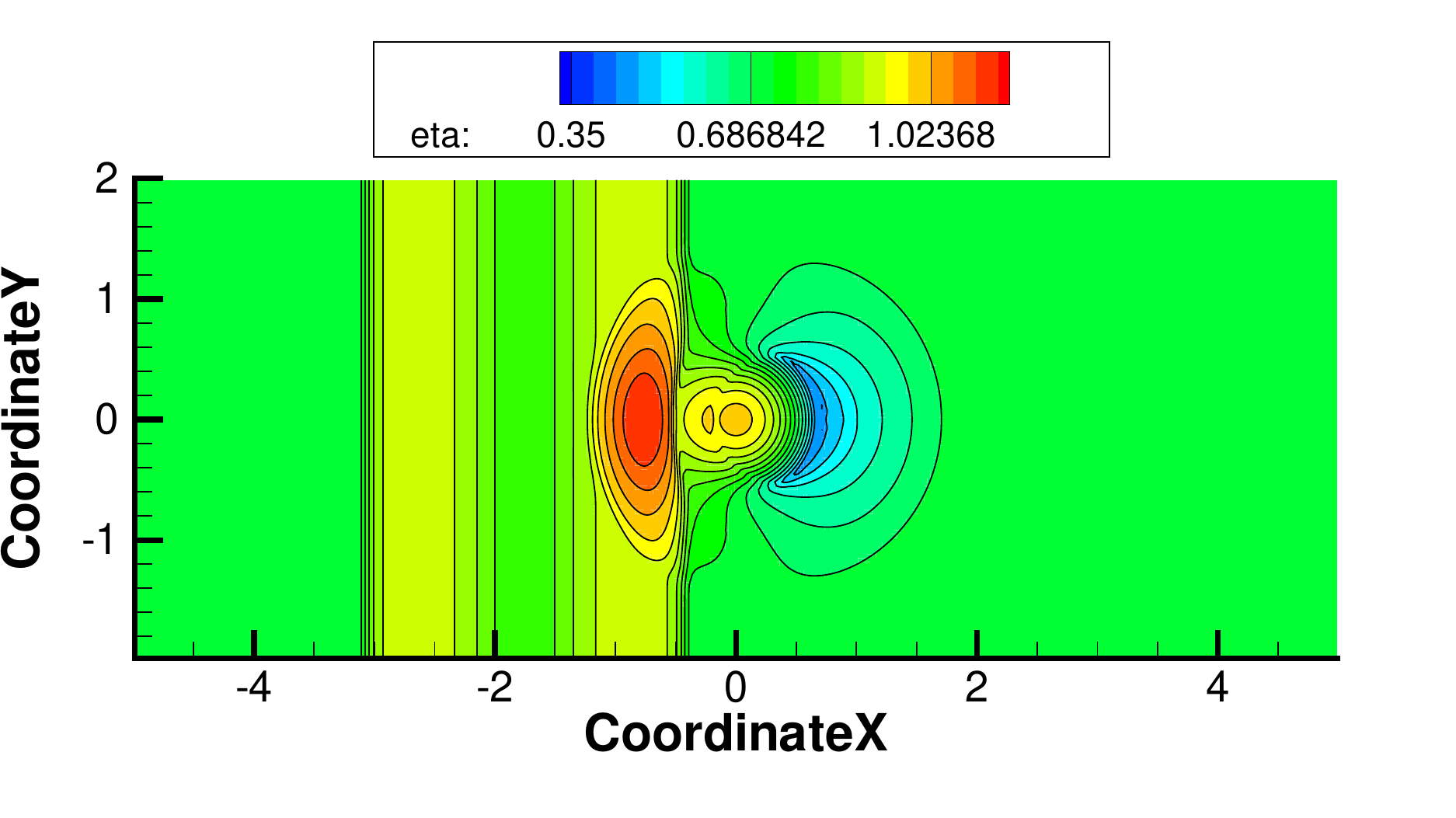}}
	\qquad
	\subfigure[$t=0.5$]{\includegraphics[width=0.45\textwidth,trim={1cm 0cm 1cm 8cm},clip]{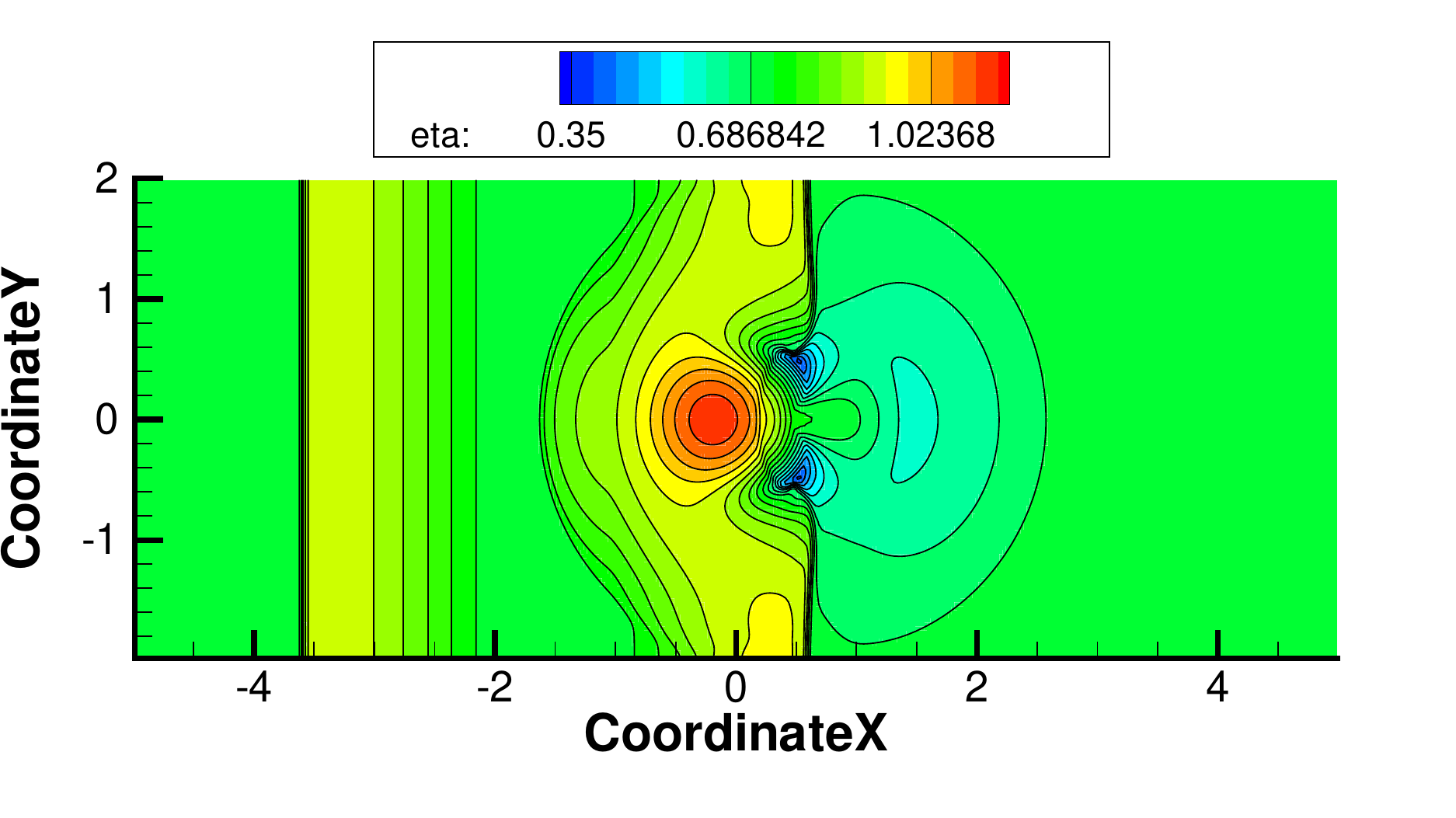}}
	\subfigure[$t=0.75$]{\includegraphics[width=0.45\textwidth,trim={1cm 0cm 1cm 8cm},clip]{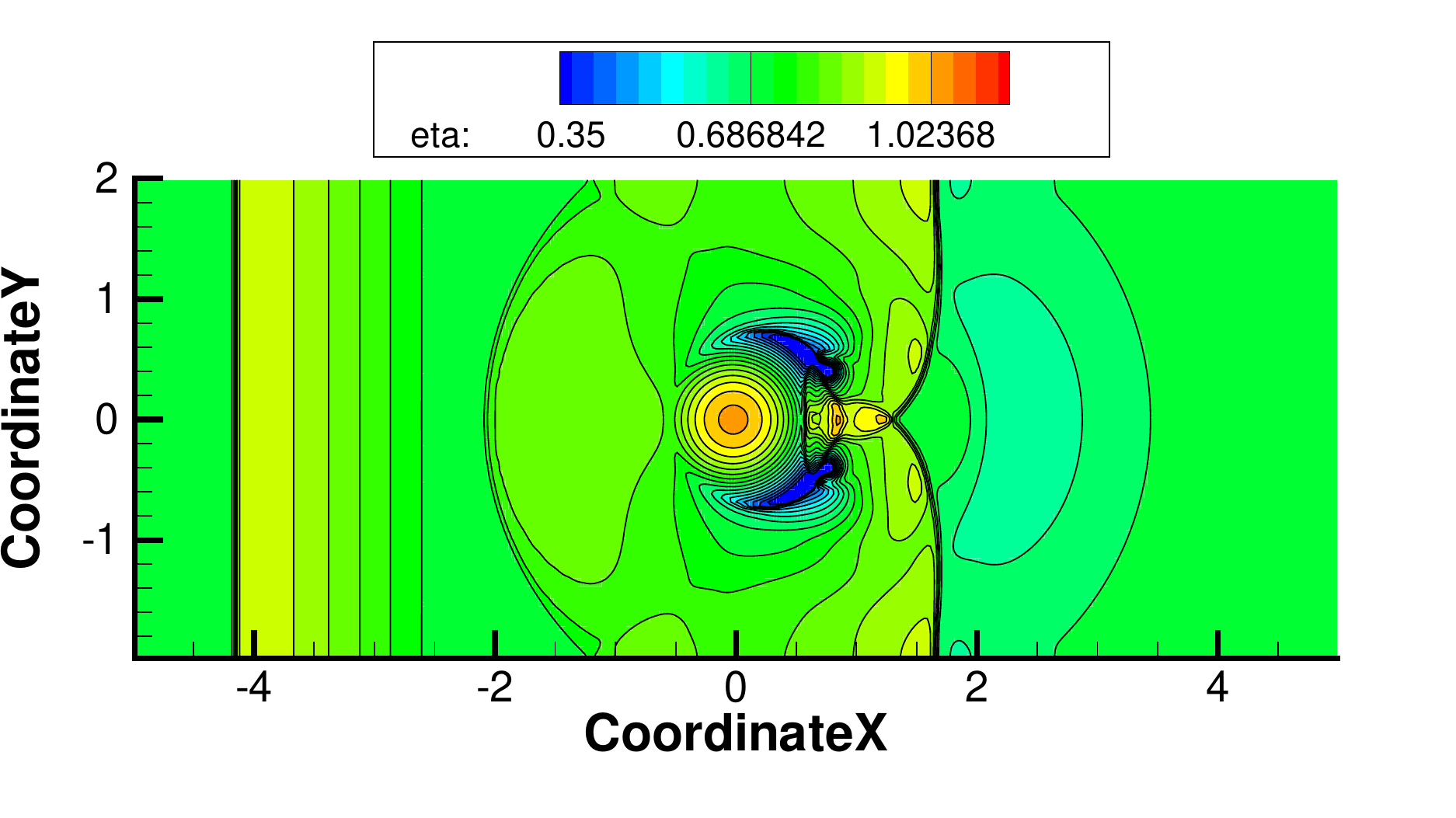}}
	\qquad
	\subfigure[$t=1$]{\includegraphics[width=0.45\textwidth,trim={1cm 0cm 1cm 8cm},clip]{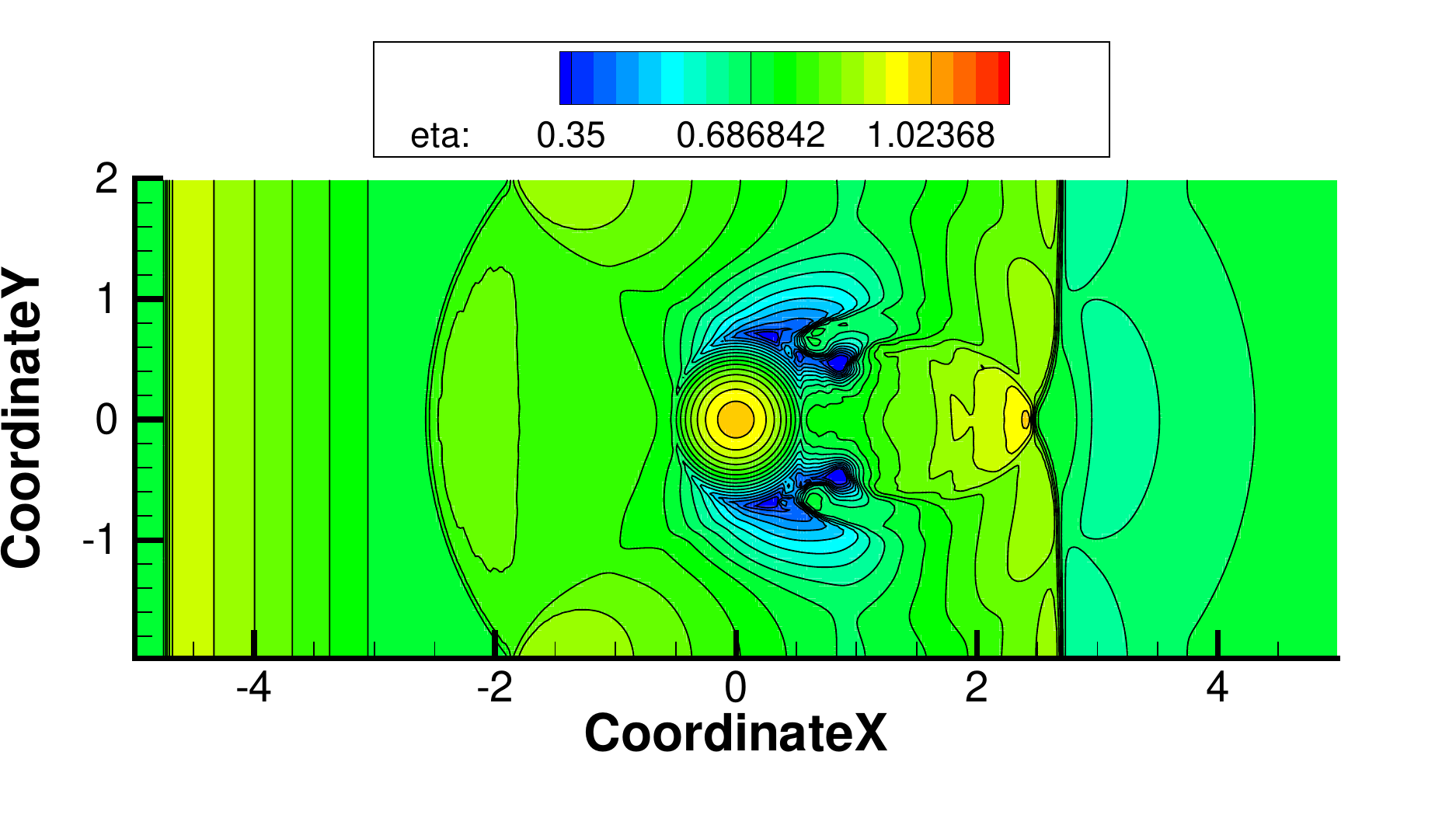}}
	\caption{Wave over dry island: $\eta=h+b$ isocontours at different times.}\label{fig:wodi_sol}
\end{figure}
It should be noticed that the top of the island, which is dry at the beginning, gets wet and dry several times during the simulation while 
never giving rise to problems due to negative water height. This is instead something that we cannot ensure for the DeC5 case. 
The simulation starts with a background state moving with speed $u=1$ which helps the wave traveling towards the island and immediately starts wetting the island from the left and drying it on the right.
The wave breaks into two smaller waves respectively approaching and moving away from the island following the eigenvalues of the flux Jacobian.
At time $t=0.25$ a run up is happening where the traveling wave is trying to submerge the top of the island while, on the other side,
a section of the island is drying. In consecutive times, the wave overtakes the island causing different interacting shock fronts and 
two symmetrical minimum points highlighted in dark blue at time $t=1$.
Several structures could been observed in this simulations and the repeated wetting/drying procedures never results in troubles for the mPDeC method.

In Figure~\ref{fig:island_jacobi} we depict the average number of Jacobi iteration needed for every linear system. The plot compares different CFL numbers and mesh refinements which keep the aspect ratio. The average is plotted with a confidence interval given by $\pm \frac12$ standard deviation. We observe that the CFL number is very incisive in determining how many Jacobi iterations are needed. For low CFL it seems that larger matrices needs more iterations, but this is not uniform with all the tests, as shown in the unsteady vortex test.

\section{Summary and Outlook} \label{se:summary}

We have presented a new well-balanced, positivity preserving high order 
numerical method for solving the shallow water equations. 
By re-writing the WENO semi-discretization in terms of productions-destructions terms, we were able 
to use the modified Patankar Deferred Correction methods of arbitrarily high-order to ensure the unconditionally positivity of water height.  
The restriction on the CFL number comes only from the explicit DeC-solver for the momentum equations. 
The used CFL numbers are of the order of 1 and are much larger than the ones used in explicit SSPRK WENO methods in combination with positivity preserving limiters~\cite{xing2014survey}, where the CFL must be lowered to $1/6$ or $1/12$.
One can relax further the CFL constraint by using implicit DeC or RK methods, though introducing more difficulties.
Even if we  have only presented a fifth order method  in this manuscript,  the approach is actually of arbitrary high-order. 
With classical explicit approaches, one must use SSPRK to guarantee positivity and, as known from literature~\cite{gottlieb2011strong}, explicit (implicit) SSPRK methods exist only up to order four (sixth) for general cases.

By applying mPDeC, we avoid those issues and the price to pay is that of solving a (very sparse) linear system for the water height. 
However, as mentioned before, this increase in computational costs is around $18 \%$ percent in our numerical simulations, but the procedure can still be optimized. 
In the future, a detailed performance test in terms of accuracy and run-time is planned. \\
Additionally, in this work we have only considered positivity preservation and well-balanced properties. In the next step, we would like to extend our investigation to entropy conservative/stable methods. Here, various approaches exist as described \textit{inter alia} in~\cite{abgrall2018general, abgrall2019analysis, chenreview, fisher2013discretely, fjordholm2012arbitrarily, ranocha2017shallow}, 
but from our perspective the convex limiting strategies seems the most promising one~\cite{hajduk2021monolithic, kuzmin2020entropycg}.
In order to do so, the basic stability properties of Patankar methods as ODE solvers have to be first fully understood \cite{izgin2022lyapunov, torlo2021stability}. \\
Finally, herein we have focused on the shallow water system. However, the method can be easily adapted to more complex models, e.g.\ the shallow water equations together 
with biochemical processes like algae bloom in oceans, seas and open water cancels or the Euler equations of gas dynamics, where special treatments for the pressure positivity are necessary \cite{zhang2012positivity,wang2012robust}. Those extensions can and will be also considered in future works.

\appendix
\section{Reconstruction of the primitive variable}\label{sec:PrimitiveVarRecon}

The goal of this section is that of deepening the procedure to compute the coeffiecients of the polynomials and linear weights needed to compute the WENO procedure in Section \ref{sec_Space}. The results and the actual coefficients for WENO5 with 4 Gaussian quadrature points are written in Section~\ref{sec:WENO5GP4}. 

\begin{figure}[ht]
        \centering
\begin{tikzpicture}
\draw [thick]   (0,0) -- (10,0) node[below=2mm] {}; 
\fill[black]    (0,0) circle (1mm) node[below=2mm] {$x_{i-5/2}$}
                (2,0) circle (1mm) node[below=2mm] {$x_{i-3/2}$}
                (4,0) circle (1mm) node[below=2mm] {$x_{i-1/2}$}
                (6,0) circle (1mm) node[below=2mm] {$x_{i+1/2}$} 
                (8,0) circle (1mm) node[below=2mm] {$x_{i+3/2}$}
                (10,0) circle(1mm) node[below=2mm] {$x_{i+5/2}$};
\node at (1,0.5) {$u_{i-2}$};
\node at (3,0.5) {$u_{i-1}$};
\node at (5,0.5) {$u_{i}$};
\node at (7,0.5) {$u_{i+1}$};
\node at (9,0.5) {$u_{i+2}$};
\end{tikzpicture} \caption{Stencil of five cell averages for WENO5 reconstruction.}\label{Fig:Space_Stencil}
\end{figure}
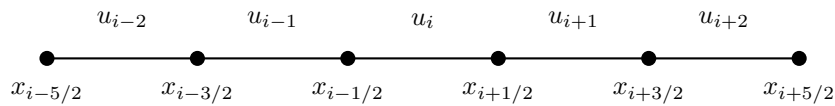

In the following assemble the system that allow to find the coefficients for the WENO polynomial at any quadrature point. The symbolic Matlab script coded for this purpose is also available 
at~\cite{ourrepo}. With few adjustments, the script can be used to compute all the ingredients needed for a 
WENO reconstruction of arbitrary high order with arbitrary high order quadrature formulae.\\
The reconstruction of the primitive variables is needed in a high order finite volume method, as  only cell averages are available from the previous time step. For the sake of simplicity we are going to 
work in a simpler one-dimensional scalar framework since the reconstruction is done dimension-by-dimension.
We consider a scalar function $u(x)$ whose cell averages $u_i$ are known. We aim at reconstructing this variable as a polynomial $v(x)$, where the polynomial may vary in different points. In particular, it will be useful to use a primitive of $v(x)$ which we denote by $\mathcal{P}(x)=\int_{x_0}^x v(s) \diff{s}$. Indeed, we can impose that for every cell average, we have
\begin{equation}\label{eq:sys_coeffWENO_cellave}
u_{i}=\frac{1}{\Delta x}\int_{\xin}^{\xip} u(x)\;\diff{x} \stackrel{!}{=} \int_{\xin}^{\xip} v(x) \diff{x} =\frac{\mathcal{P}(x_{\iip} )- \mathcal{P}(x_{\iin})}{\Delta x}.
\end{equation}
The stencil considered for this example is the one used for the WENO5 reconstruction, which is made up by five cell averages
as denoted in Figure~\ref{Fig:Space_Stencil}.

The next step consists in using interpolating polynomials $\varphi_{j-1/2}$, e.g.\ Lagrange polynomials at cell interfaces, to approximate the primitive $\mathcal{P}(x)$, i.e.,
\begin{equation}\label{eq:lagrange}
\mathcal{P}(x) = \sum_{j=k_1}^{k_2+1} a_{j-1/2}\,\varphi_{j-1/2}(x) \qquad \text{where}\qquad \varphi_{j-1/2}(x) = \prod_{\ell=k_1, i \neq j}^{k_2+1} \frac{x-x_{\ell-1/2}}{x_{j-1/2}-x_{\ell-1/2}}
\end{equation}
where $k_1$ and $k_2$ are the extreme indexes of the considered stencil. So, given then $k_2-k_1+1$ cell averages, the degree of the polynomial $\mathcal{P}$ will be $k_2-k_1+2$.
For instance, when working with WENO5, the reconstruction is composed by a linear combination of three lower order polynomials of degree 3, so with $\mathcal{P} \in \mathbb P_4$, obtained through three cell averages. Then, the 3 polynomials, which then depend on the whole 5 cells stencil, will be be combined with weights which depends on the quadrature points into a fifth order accurate reconstruction. 
In order to have this result, we need to compute the aforementioned three lower order polynomials and a high order polynomial made up
by information coming from the whole stencil.\\
Let us begin with the procedure to compute the high order polynomial with all available cell averages, i.e., $k_1=-2$ and $k_2=2$.
The lower order polynomials can be easily computed following the same approach explained hereafter.
In this case $\mathcal{P}(x)$ is a sixth order polynomial and $v(x)=\mathcal{P}'(x)$ is a fifth order polynomial which gives the right accuracy order. 
Using \eqref{eq:sys_coeffWENO_cellave} for all cell averages in the stencil with the definition of $\mathcal{P}$ given in \eqref{eq:lagrange} and using the property of Lagrangian polynomials $\varphi_{j-1/2}(x_{\ell -1/2}) = \delta_{j,\ell}$, we obtain a system of equations for the coefficients $a_{j-1/2}$ 
with solution
\begin{equation}
a_{j-1/2} = \begin{cases}
 0,  & \text{if }j=k_1, \\ 
 \sum_{\ell=-k_1}^{-k_1+j-1} u_{i+\ell}, & j=k_1+1,\dots, k_2+1. 
\end{cases}
\end{equation}
Finally, the expression for the high order, $ho$, approximation polynomial $v^{ho}(x)$ can be written as 
\begin{equation}
v^{ho}(x)= \sum_{j=k_1}^{k_2+1} a_{j-1/2}\varphi'_{j-1/2}(x) = \sum_{\ell=-k_1}^{k_2} c^{ho}_{\ell}(x) \,u_{i+\ell}= \sum_{\ell=-2}^{2} c^{ho}_{\ell}(x) \,u_{i+\ell},
\end{equation}
where $c^{ho}_\ell$ are obtained collecting all the coefficients and basis functions related to $u_{i+\ell}$. In this way, for any quadrature point $\tilde\xi \in [x_{i-1/2},x_{i+1/2}]$ we can evaluate this high order polynomial $v^{ho}(\tilde \xi)$.
Following the same procedure for three lower order polynomials, $lo$, associated to the 3-cells stencils
$S^0=\{u_{i},u_{i+1},u_{i+2}\}$, $S^1=\{u_{i-1},u_{i},u_{i+1}\}$ and $S^2=\{u_{i-2},u_{i-1},u_{i}\}$,  we obtain an expression for these low order polynomials
\begin{equation}
v^{lo}_j(x) = \sum_{\ell=-j}^{2-j} c^{lo}_{j\ell}(x) \,u_{i+\ell}.
\end{equation}
The last step concerns the computation of the ideal linear weights. These are the weights that allow to recover the high order reconstruction
from a linear combination of the lower order ones for a given quadrature point $\tilde \xi$. Therefore, we need to find the linear weights $d_j$ such that
\begin{equation}\label{eq:ideal_w_problem}
v^{ho}(\tilde \xi) = \sum_{j=0}^2 d_j \,v^{lo}_j(\tilde \xi) \quad\Longleftrightarrow\quad 
\sum_{\ell=-2}^{2} c^{ho}_{\ell}(\tilde \xi) \,u_{i+\ell} = \sum_{j=0}^2 \sum_{\ell=-j}^{2-j} d_j \, c^{lo}_{j\ell}(\tilde \xi) \,u_{i+\ell}, \quad \forall u_{i+\ell}.
\end{equation}
Since \eqref{eq:ideal_w_problem} must hold for any quintuplet $\lbrace u_{i+\ell} \rbrace _{\ell=-2}^{2}$, we can write a system of five equations in the 3 linear weights $d_j$ for each quadrature point $\tilde \xi$, i.e.,
\begin{equation}
\sum_{j=0}^2 d_j \, c^{lo}_{j \ell}(\tilde \xi) = c^{ho}_{\ell}(\tilde \xi), \quad \forall \ell = -2,\dots,2. \\
\end{equation}
This is an overdetermined system with five equations and only three unknowns $d_j$ that can be easily solved by means of a least squares method. Moreover, the solutions found verify exactly all the equations, implying that some of the equations are linearly dependent.

\section{WENO recontruction ($r=3$) with four-point Gaussian quadrature rule}\label{sec:WENO5GP4}
The goal of this section is to present the fifth order WENO reconstruction with four-point Gaussian quadrature rule.
Up to our knowledge, there is no reference in literature that explicitly define linear weights and polynomial coefficients needed for this WENO reconstruction.
Reference~\cite{titarev2004finite} well described the fifth-order WENO5 with two-point Gaussian quadrature rule, which unfortunately does not
allow to go beyond fourth order. Instead, with the four-point Gaussian quadrature rule, one could reach even eighth order.\\ 
Let us consider a one-dimensional cell $[\xi_{\iin},\xi_{\iip}]$, we hereafter provide the expressions for 
\begin{equation}
q(\xi_{\iip}^-)\;,\;\; q(\xi_{\iin}^+)\;,\;\; q\left(\xi_i\pm\frac{\Delta\xi}{2}\sqrt{\frac{3}{7}+\frac{2}{7}\sqrt{\frac{6}{5}}}\right)\;,\;\;  q\left(\xi_i\pm\frac{\Delta\xi}{2}\sqrt{\frac{3}{7}-\frac{2}{7}\sqrt{\frac{6}{5}}}\right),
\end{equation}
which are used for the first sweep (first two terms corresponding to the two boundaries) and the second sweep (the last two terms corresponding to the 4 quadrature points).
For $r=3$ we have only three candidate stencil for the reconstruction
\begin{equation}
S_0 = (i,i+1,i+2)\;,\;\;\;S_1 = (i-1,i,i+1)\;,\;\;\;S_2 = (i-2,i-1,i).
\end{equation}
The corresponding smoothness indicators are given by:
\begin{align*}
\beta_0 &= \frac{13}{12} \left( q_i - 2q_{i+1} + q_{i+2}\right)^2 + \frac{1}{4} \left( 3q_i - 4q_{i+1} + q_{i+2} \right)^2,   \\ 
\beta_1 &= \frac{13}{12} \left( q_{i-1} - 2q_i + q_{i+1}\right)^2 + \frac{1}{4} \left( q_{i-1} - q_{i+1} \right)^2,   \\ 
\beta_2 &= \frac{13}{12} \left( q_{i-2} - 2q_{i-1} + q_i\right)^2 + \frac{1}{4} \left( q_{i-2} - 4q_{i-1} + 3q_i \right)^2.      
\end{align*}
The optimal weights $d_m$ for the left boundary extrapolated value $q_{\iip}^-$ at $\xip$ are 
\begin{equation}
d_0 = \frac{3}{10}  \;,\;\;\;   d_1 = \frac{3}{5}    \;,\;\;\;  d_2 = \frac{1}{10}
\end{equation}
and $q_{\iip}^-$ is given by
\begin{equation}
q_{\iip}^- = \frac{1}{6}\omega_0 \left(-q_{i+2} + 5q_{i+1} + 2q_i\right) + \frac{1}{6}\omega_1 \left(-q_{i-1} + 5q_{i} + 2q_{i+1}\right) + \frac{1}{6}\omega_2 \left(2q_{i-2} - 7q_{i-1} + 11q_i\right).
\end{equation}
The optimal weights $d_m$ for the right boundary extrapolated value $q_{\iin}^+$ at $\xin$ are 
\begin{equation}
d_0 = \frac{1}{10}  \;,\;\;\;   d_1 = \frac{3}{5}    \;,\;\;\;  d_2 = \frac{3}{10}
\end{equation}
and $q_{\iin}^+$ is given by
\begin{equation}
q_{\iin}^+ = \frac{1}{6}\omega_0 \left(2q_{i+2} - 7q_{i+1} + 11q_i\right) + \frac{1}{6}\omega_1 \left(-q_{i+1} + 5q_{i} + 2q_{i-1}\right) + \frac{1}{6}\omega_2 \left(-q_{i-2} + 5q_{i-1} + 2q_i\right).
\end{equation}
For the first Gaussian quadrature point $\xi_1^q = \xi_i - \frac{\Delta\xi}{2}\sqrt{\frac{3}{7}+\frac{2}{7}\sqrt{\frac{6}{5}}}$, 
the optimal weights are:
\begin{align}
\begin{aligned}
d_0 &= \frac{269\,\sqrt{42}\,\sqrt{2\,\sqrt{30}+15}}{50428}-\frac{1751\,\sqrt{35}\,\sqrt{2\,\sqrt{30}+15}}{504280}-\frac{411\,\sqrt{30}}{100856}+\frac{21855}{100856}, \\
d_1 &= \frac{411\,\sqrt{30}}{50428}+\frac{28573}{50428}, \\
d_2 &= \frac{1751\,\sqrt{35}\,\sqrt{2\,\sqrt{30}+15}}{504280}-\frac{269\,\sqrt{42}\,\sqrt{2\,\sqrt{30}+15}}{50428}-\frac{411\,\sqrt{30}}{100856}+\frac{21855}{100856}.
\end{aligned}
\end{align}
The reconstructed value can be computed from the three polynomials associated to each stencil:
{\tiny 
\begin{align*}
p_0(\xi_1^q) &= \left(\frac{\sqrt{5}\,\sqrt{6}}{140}+\frac{3\,\sqrt{\frac{2\,\sqrt{5}\,\sqrt{6}}{35}+\frac{3}{7}}}{4}+\frac{85}{84}\right) \,q_i + \left(-\frac{\sqrt{5}\,\sqrt{6}}{70}-\sqrt{\frac{2\,\sqrt{5}\,\sqrt{6}}{35}+\frac{3}{7}}-\frac{1}{42}\right) \,q_{i+1} + \left(\frac{\sqrt{5}\,\sqrt{6}}{140}+\frac{\sqrt{\frac{2\,\sqrt{5}\,\sqrt{6}}{35}+\frac{3}{7}}}{4}+\frac{1}{84}\right) \,q_{i+2} ,\\ 
p_1(\xi_1^q) &= \left(\frac{\sqrt{5}\,\sqrt{6}}{140}+\frac{\sqrt{\frac{2\,\sqrt{5}\,\sqrt{6}}{35}+\frac{3}{7}}}{4}+\frac{1}{84}\right) \,q_{i-1} + \left(\frac{41}{42}-\frac{\sqrt{5}\,\sqrt{6}}{70}\right) \,q_{i} +  \left(\frac{\sqrt{5}\,\sqrt{6}}{140}-\frac{\sqrt{\frac{2\,\sqrt{5}\,\sqrt{6}}{35}+\frac{3}{7}}}{4}+\frac{1}{84}\right) \,q_{i+1}, \\ 
p_2(\xi_1^q) &= \left(\frac{\sqrt{5}\,\sqrt{6}}{140}-\frac{\sqrt{\frac{2\,\sqrt{5}\,\sqrt{6}}{35}+\frac{3}{7}}}{4}+\frac{1}{84}\right) \,q_{i-2} + \left(\sqrt{\frac{2\,\sqrt{5}\,\sqrt{6}}{35}+\frac{3}{7}}-\frac{\sqrt{5}\,\sqrt{6}}{70}-\frac{1}{42}\right) \,q_{i-1} + \left(\frac{\sqrt{5}\,\sqrt{6}}{140}-\frac{3\,\sqrt{\frac{2\,\sqrt{5}\,\sqrt{6}}{35}+\frac{3}{7}}}{4}+\frac{85}{84}\right) \,q_{i}  .
\end{align*}
}%
For the second Gaussian quadrature point $\xi_2^q = \xi_i - \frac{\Delta\xi}{2}\sqrt{\frac{3}{7}-\frac{2}{7}\sqrt{\frac{6}{5}}}$, 
the optimal weights are:
\begin{align}
\begin{aligned}
d_0 &= \frac{411\,\sqrt{30}}{100856}-\frac{269\,\sqrt{42}\,\sqrt{15-2\,\sqrt{30}}}{50428}-\frac{1751\,\sqrt{35}\,\sqrt{15-2\,\sqrt{30}}}{504280}+\frac{21855}{100856}, \\ 
d_1 &= \frac{28573}{50428}-\frac{411\,\sqrt{30}}{50428}, \\
d_2 &= \frac{1751\,\sqrt{35}\,\sqrt{15-2\,\sqrt{30}}}{504280}+\frac{269\,\sqrt{42}\,\sqrt{15-2\,\sqrt{30}}}{50428}+\frac{411\,\sqrt{30}}{100856}+\frac{21855}{100856}.
\end{aligned}
\end{align}
The reconstructed value can be then computed from the three polynomials associated to the three stencil:
{\tiny 
\begin{align*}
p_0(\xi_2^q) &= \left(\frac{3\,\sqrt{\frac{3}{7}-\frac{2\,\sqrt{5}\,\sqrt{6}}{35}}}{4}-\frac{\sqrt{5}\,\sqrt{6}}{140}+\frac{85}{84}\right) \,q_i     + \left(\frac{\sqrt{5}\,\sqrt{6}}{70}-\sqrt{\frac{3}{7}-\frac{2\,\sqrt{5}\,\sqrt{6}}{35}}-\frac{1}{42}\right) \,q_{i+1} + \left(\frac{\sqrt{\frac{3}{7}-\frac{2\,\sqrt{5}\,\sqrt{6}}{35}}}{4}-\frac{\sqrt{5}\,\sqrt{6}}{140}+\frac{1}{84}\right) \,q_{i+2}, \\ 
p_1(\xi_2^q) &= \left(\frac{\sqrt{\frac{3}{7}-\frac{2\,\sqrt{5}\,\sqrt{6}}{35}}}{4}-\frac{\sqrt{5}\,\sqrt{6}}{140}+\frac{1}{84}\right) \,q_{i-1} + \left(\frac{\sqrt{5}\,\sqrt{6}}{70}+\frac{41}{42}\right) \,q_{i}   + \left(\frac{1}{84}-\frac{\sqrt{\frac{3}{7}-\frac{2\,\sqrt{5}\,\sqrt{6}}{35}}}{4}-\frac{\sqrt{5}\,\sqrt{6}}{140}\right) \,q_{i+1}, \\ 
p_2(\xi_2^q) &= \left(\frac{1}{84}-\frac{\sqrt{\frac{3}{7}-\frac{2\,\sqrt{5}\,\sqrt{6}}{35}}}{4}-\frac{\sqrt{5}\,\sqrt{6}}{140}\right) \,q_{i-2} + \left(\frac{\sqrt{5}\,\sqrt{6}}{70}+\sqrt{\frac{3}{7}-\frac{2\,\sqrt{5}\,\sqrt{6}}{35}}-\frac{1}{42}\right) \,q_{i-1} + \left(\frac{85}{84}-\frac{3\,\sqrt{\frac{3}{7}-\frac{2\,\sqrt{5}\,\sqrt{6}}{35}}}{4}-\frac{\sqrt{5}\,\sqrt{6}}{140}\right) \,q_{i}  .
\end{align*}
}%
For the third Gaussian quadrature point $\xi_3^q = \xi_i + \frac{\Delta\xi}{2}\sqrt{\frac{3}{7}-\frac{2}{7}\sqrt{\frac{6}{5}}}$, 
the optimal weights are:
\begin{align}
\begin{aligned}
d_0 &= \frac{1751\,\sqrt{35}\,\sqrt{15-2\,\sqrt{30}}}{504280}+\frac{269\,\sqrt{42}\,\sqrt{15-2\,\sqrt{30}}}{50428}+\frac{411\,\sqrt{30}}{100856}+\frac{21855}{100856}, \\
d_1 &= \frac{28573}{50428}-\frac{411\,\sqrt{30}}{50428}, \\
d_2 &= \frac{411\,\sqrt{30}}{100856}-\frac{269\,\sqrt{42}\,\sqrt{15-2\,\sqrt{30}}}{50428}-\frac{1751\,\sqrt{35}\,\sqrt{15-2\,\sqrt{30}}}{504280}+\frac{21855}{100856}.
\end{aligned}
\end{align}
The reconstructed value can be then computed from the three polynomials associated to the three stencil:
{\tiny
\begin{align*}
p_0 (\xi_3^q)&= \left(\frac{85}{84}-\frac{3\,\sqrt{\frac{3}{7}-\frac{2\,\sqrt{5}\,\sqrt{6}}{35}}}{4}-\frac{\sqrt{5}\,\sqrt{6}}{140}\right) \,q_i     + \left(\frac{\sqrt{5}\,\sqrt{6}}{70}+\sqrt{\frac{3}{7}-\frac{2\,\sqrt{5}\,\sqrt{6}}{35}}-\frac{1}{42}\right)\,q_{i+1} + \left(\frac{1}{84}-\frac{\sqrt{\frac{3}{7}-\frac{2\,\sqrt{5}\,\sqrt{6}}{35}}}{4}-\frac{\sqrt{5}\,\sqrt{6}}{140}\right)\,q_{i+2}, \\ 
p_1(\xi_3^q) &= \left(\frac{1}{84}-\frac{\sqrt{\frac{3}{7}-\frac{2\,\sqrt{5}\,\sqrt{6}}{35}}}{4}-\frac{\sqrt{5}\,\sqrt{6}}{140}\right) \,q_{i-1} + \left(\frac{\sqrt{5}\,\sqrt{6}}{70}+\frac{41}{42}\right)\,q_{i}   + \left(\frac{\sqrt{\frac{3}{7}-\frac{2\,\sqrt{5}\,\sqrt{6}}{35}}}{4}-\frac{\sqrt{5}\,\sqrt{6}}{140}+\frac{1}{84}\right)\,q_{i+1}, \\ 
p_2(\xi_3^q) &= \left(\frac{\sqrt{\frac{3}{7}-\frac{2\,\sqrt{5}\,\sqrt{6}}{35}}}{4}-\frac{\sqrt{5}\,\sqrt{6}}{140}+\frac{1}{84}\right) \,q_{i-2} + \left(\frac{\sqrt{5}\,\sqrt{6}}{70}-\sqrt{\frac{3}{7}-\frac{2\,\sqrt{5}\,\sqrt{6}}{35}}-\frac{1}{42}\right)\,q_{i-1} + \left(\frac{3\,\sqrt{\frac{3}{7}-\frac{2\,\sqrt{5}\,\sqrt{6}}{35}}}{4}-\frac{\sqrt{5}\,\sqrt{6}}{140}+\frac{85}{84}\right)\,q_{i}  .
\end{align*}
}%
For the fourth Gaussian quadrature point $\xi_4^q = \xi_i + \frac{\Delta\xi}{2}\sqrt{\frac{3}{7}+\frac{2}{7}\sqrt{\frac{6}{5}}}$, 
the optimal weights are:
\begin{align}
\begin{aligned}
d_0 &= \frac{1751\,\sqrt{35}\,\sqrt{2\,\sqrt{30}+15}}{504280}-\frac{269\,\sqrt{42}\,\sqrt{2\,\sqrt{30}+15}}{50428}-\frac{411\,\sqrt{30}}{100856}+\frac{21855}{100856},  \\ 
d_1 &= \frac{411\,\sqrt{30}}{50428}+\frac{28573}{50428} , \\
d_2 &= \frac{269\,\sqrt{42}\,\sqrt{2\,\sqrt{30}+15}}{50428}-\frac{1751\,\sqrt{35}\,\sqrt{2\,\sqrt{30}+15}}{504280}-\frac{411\,\sqrt{30}}{100856}+\frac{21855}{100856}.
\end{aligned}
\end{align}
The reconstructed value can be then computed from the three polynomials associated to the three stencil:
{\tiny
\begin{align*}
p_0 (\xi_4^q)&= \left(\frac{\sqrt{5}\,\sqrt{6}}{140}-\frac{3\,\sqrt{\frac{2\,\sqrt{5}\,\sqrt{6}}{35}+\frac{3}{7}}}{4}+\frac{85}{84}\right) \,q_i     + \left(\sqrt{\frac{2\,\sqrt{5}\,\sqrt{6}}{35}+\frac{3}{7}}-\frac{\sqrt{5}\,\sqrt{6}}{70}-\frac{1}{42}\right) \,q_{i+1} + \left(\frac{\sqrt{5}\,\sqrt{6}}{140}-\frac{\sqrt{\frac{2\,\sqrt{5}\,\sqrt{6}}{35}+\frac{3}{7}}}{4}+\frac{1}{84}\right) \,q_{i+2}, \\ 
p_1(\xi_4^q) &= \left(\frac{\sqrt{5}\,\sqrt{6}}{140}-\frac{\sqrt{\frac{2\,\sqrt{5}\,\sqrt{6}}{35}+\frac{3}{7}}}{4}+\frac{1}{84}\right) \,q_{i-1} + \left(\frac{41}{42}-\frac{\sqrt{5}\,\sqrt{6}}{70}\right) \,q_{i}   + \left(\frac{\sqrt{5}\,\sqrt{6}}{140}+\frac{\sqrt{\frac{2\,\sqrt{5}\,\sqrt{6}}{35}+\frac{3}{7}}}{4}+\frac{1}{84}\right) \,q_{i+1} ,\\ 
p_2 (\xi_4^q)&= \left(\frac{\sqrt{5}\,\sqrt{6}}{140}+\frac{\sqrt{\frac{2\,\sqrt{5}\,\sqrt{6}}{35}+\frac{3}{7}}}{4}+\frac{1}{84}\right) \,q_{i-2} + \left(-\frac{\sqrt{5}\,\sqrt{6}}{70}-\sqrt{\frac{2\,\sqrt{5}\,\sqrt{6}}{35}+\frac{3}{7}}-\frac{1}{42}\right) \,q_{i-1} + \left(\frac{\sqrt{5}\,\sqrt{6}}{140}+\frac{3\,\sqrt{\frac{2\,\sqrt{5}\,\sqrt{6}}{35}+\frac{3}{7}}}{4}+\frac{85}{84}\right) \,q_{i}  .
\end{align*}
}%
%
For all Gaussian quadrature points the solution in $\xi$ can be easily built by assembling 
the three polynomials.
%

{\small
\subsection*{Acknowledgements}
M. Ciallella is funded by an Inria PhD fellowship.
L. Micalizzi is supported by SNF, project number 200020$\_$175784, and by the Forschungskredit grant FK-21-098.
D. Torlo is funded by an Inria Postdoc.
P. \"Offner gratefully acknowledge the support of the Gutenberg Research College and also wants to thank Mario Ricchiuto for his invitation to Inria Bordeaux.  \\
All authors would like to thank Jonatan N\'u\~nez for sharing his high-order FV-WENO code on his repository \cite{Nunezrepo}. We have started our work by adapting his code. 
 \bibliographystyle{abbrv}
\bibliography{literature}
}
\end{document}